%% file: main.tex
\RequirePackage{fix-cm}
\documentclass[smallcondensed]{svjour3}     
\smartqed  
\usepackage{graphicx}
\usepackage{amsmath,amssymb, amsfonts,bm,tikz,url}
\usepackage{xcolor}
\usetikzlibrary{hobby}
\input{CommPreamble}
\usepackage[utf8]{inputenc}
\usetikzlibrary{decorations.pathreplacing,calligraphy}
\usepackage{bm}
\usepackage{mathrsfs}
\usepackage{amsmath}
\usepackage{amsfonts}
\usepackage{color}
\usepackage{setspace}
\usepackage{float}
\usepackage{url}
\usepackage{multicol}
\usepackage{subfigure}
\usepackage{bm}
\usepackage[top=1in,bottom=1in,left=1in,right=1in]{geometry}

\newcommand{\mcL}{\mathcal{L}}

\usepackage{xcolor}

\def\omg{{\Omega}}

\def \phib{{\boldsymbol \phi}}
\def \psib{{\boldsymbol \psi}}

\def \fb{\vf}

\def \ub{\vu}

\def \xb{\vx}

\def \yb{\vy}

\newcommand{\verti}[1]{{\left\vert #1
    \right\vert}}
\newcommand{\vertii}[1]{{\left\vert\left\vert #1
    \right\vert\right\vert}}

\def\PK{\operatorname{P-K}}

\begin{document}

\title{Convergence Analysis and Numerical Studies for Linearly Elastic Peridynamics with Dirichlet-Type Boundary Conditions
\thanks{M. Foss was supported by NSF - DMS 1716790.}
\thanks{P. Radu was supported by NSF - DMS  1716790.}
\thanks{Y. Yu was supported by NSF - DMS 1753031.}}

\author{Mikil Foss \and Petronela Radu \and Yue Yu}

\institute{M. Foss \at
Department of Mathematics, University of Nebraska-Lincoln, NE\\
\email{mfoss@unl.edu}
\and
P. Radu \at
Department of Mathematics, University of Nebraska-Lincoln, NE\\
\email{pradu@unl.edu}
\and
Y. Yu \at
Department of Mathematics, Lehigh University, Bethlehem, PA\\
\email{yuy214@lehigh.edu}
}

\date{Received: date / Accepted: date - Draft, ***, 2021}

\maketitle

\begin{abstract}
The nonlocal models of peridynamics have successfully predicted fractures and deformations for a variety of materials. In contrast to local mechanics, peridynamic boundary  conditions must be defined on a finite volume region outside the body. Therefore, theoretical and numerical challenges arise in order to properly formulate Dirichlet-type nonlocal boundary conditions, while connecting them to the local counterparts. While a careless imposition of local boundary conditions leads to a smaller effective material stiffness close to the boundary and an artificial softening of the material, several strategies were proposed to avoid this unphysical surface effect.

In this work, we study convergence of solutions to nonlocal state-based linear elastic model to their local counterparts as the interaction horizon vanishes, under different formulations and smoothness assumptions for nonlocal Dirichlet-type boundary conditions. Our results provide explicit rates of convergence that are sensitive to the compatibility of the nonlocal boundary data and the extension of the solution for the local model. In particular, under appropriate assumptions, constant extensions yield $\frac{1}{2}$ order convergence rates and linear extensions yield $\frac{3}{2}$ order convergence rates. With smooth extensions, these rates are improved to quadratic convergence. We illustrate the theory for any dimension $d\geq 2$ and numerically verify the convergence rates with a number of two dimensional benchmarks, including linear patch tests, manufactured solutions, and domains with curvilinear surfaces. Numerical results show a first order convergence for constant extensions and second order convergence for linear extensions, which suggests a possible room of improvement in the future convergence analysis.
\end{abstract}

\tableofcontents

\section{Introduction}

The peridynamics theory formulated in \cite{silling_2000} has been successful in modeling deformations, fracture, and predicting behavior in a variety of materials: concrete, metal, viscoelastic, viscoplastic  \cite{bobaru2007influence,foster2010viscoplasticity,gerstle2005concrete,madenci2014peridynamic,madenci2016peridynamic,oterkus2010peridynamic,weckner2013viscoelastic,wu2015stabilizedmetal,yang2018concrete}. As a nonlocal model, peridynamics employs integral, rather than differential, operators which allows to relax the regularity constraints of partial differential equations (PDEs) and to capture effects arising from long range forces at the microscale and mesoscale, not accounted by PDEs \cite{bazant2002nonlocal,diehl2021comparative,diehl2019review,du2013nonlocal,javili2019peridynamics,You2021,You2020Regression}. Consequently, the model (with a variety of formulations) has been of interest even in the absence of damage, as a series of papers showed well-posedness and properties of the solutions, as well as convergence of the nonlocal operators and solutions to classical counterparts in the limit of vanishing horizon of interaction \cite{foss2019nonlocal,FR2,FRW,lipton2014dynamic,mengesha2014bond,mengesha2014nonlocal,RW,tian2014asymptotically}. 

In peridynamics the nonlocality manifests itself through interactions on a finite range. As such, boundary value problems have data supported on a collar that contains the points with whom the domain entertains interactions \cite{du2013analysis,du2013nonlocal,du2021nonlocal,gunzburger2010nonlocal}. The size of the collar (also called the interaction domain) is determined 
by a system parameter, called horizon, and usually denoted by $\delta$. Formulating a boundary value problem in the nonlocal setting brings in an additional level of complexity, as appropriate boundary conditions must be defined. While for a local problem data on the domain's surface can be easily provided by experimentalists through surface measurements, in the nonlocal setting one must provide volumetric data for the boundary collar, which may be difficult (or even impossible) to obtain. Thus, practioners must introduce ad-hoc methods for prescribing both Dirichlet and Neumann-type boundary conditions for nonlocal problem \cite{d2020physically,du2018heterogeneously,du2019uniform,Lee2021,tao2017nonlocal,you2020asymptotically,yu2021asymptotically,zhao2020algorithm}, and in this work we mainly focus on the Dirichlet-type volume constraints.  Most commonly, the convergence of nonlocal solutions to classical counterparts has been studied for homogeneous Dirichlet-type boundary conditions in second order (\cite{mengesha2014bond,mengesha2014nonlocal}), or higher order (\cite{RTY}) problems. In these works the local problem is set on domain $\Omega$, while the nonlocal equation holds on $\Omega \setminus \Gamma$ (where $\Gamma\subset \Omega$ is an interior collar set of positive measure), 
with zero boundary conditions on $\Gamma$. With this setup (where the nonlocal solution has exactly the same values as the local solution on the nonlocal boundary) it is expected to observe quadratic convergence (with respect to $\delta$) for the $L^2$ norm of the difference between the local and nonlocal solutions.

For nonhomogenous Dirichlet-type boundary data a series of approaches have been introduced and investigated for peridynamics and the general nonlocal problems. Currently, popular strategies to enforce local Dirichlet boundary conditions in nonlocal models can be generally classified into two types: by modifying the nonlocal operator near the boundary \cite{du2018heterogeneously}, or by extending or converting the surface (local) data into volumetric data \cite{du2019uniform,Lee2021}. The first approach was introduced in \cite{du2018heterogeneously}, where the authors propose to gradually change the nonlocal operator to a local operator on the boundary, so as to avoid the use of nonlocal boundary conditions. For nonlocal diffusion problems, thanks to the nonlocal trace theorem provided in \cite{tian2017trace}, this approach features well-posedness. However, to the best of our knowledge, no work has yet addressed the well-posedness of this approach in the vectorial framework of peridynamics. On the other hand, the second approach is often developed by means of constant, linear, or higher-order extrapolations of the given boundary data on the surface of codimension-one onto the boundary collar \cite{Lee2021,morris1997modeling}, and it is more commonly employed to prescribe the Dirichlet boundary condition in peridynamics. 

On scalar-valued nonlocal problems such as the nonlocal diffusion problems, in \cite{du2019uniform} the authors proved the second order convergence of the nonlocal solution to the local one with linear extensions on the boundary conditions defined based on the exact derivatives of the local limit. When the exact local derivative is not known a priori, in \cite{Lee2021}, the authors proposed to implicitly compute the appropriate normalization factors and they were able to obtain second order consistency of the nonlocal operators with the local ones. 
For vector-valued systems, in the context of peridynamics a naive fictitious nodes method (FNM) was originally proposed by assigning the same values to all fictitious points corresponding to a boundary point \cite{chen2015selecting}. The naive FNM was later extended to the Taylor FNM by using Taylor expansion (to linear terms) for points in the boundary collar \cite{wang2019modeling} and the mirror-based FNM by reflecting the values on interior collar nodes to their corresponding mirror nodes for domains with simple geometries \cite{bobaru2016handbook,bobaru2011adaptive,le2018surface,oterkus2014fully}. Recently, the mirror-based FNM was further extended to more general domains in \cite{zhao2020algorithm} for nonlocal diffusion problems, where the authors propose a novel approach of extending the solution across the boundary along a ``nonlocal normal" that is computed using the gradient of the index-function that measures the size of the support of interaction. The approach is investigated numerically and it shows better agreement with classical solutions, while also being able to handle non-smooth boundaries, and even crack lines.  Similar as the mirror FNM idea in \cite{zhao2020algorithm} and the vanishing horizon idea in \cite{du2018heterogeneously}, in \cite{prudhomme2020treatment} the authors propose two methods for dealing with boundary conditions in one-dimensional bond-based peridynamics model. The Extended Domain Method (EDM) adds a layer to the domain $\Omega$ on which an odd reflection of the classical solution is imposed. For the Variable Horizon Method (VHM) the boundary is interior to the domain and it is considered of variable depth, so that as the horizon shrinks to zero, the nonlocal solution converges to classical. The numerical studies show quadratic convergence in $\delta$, if some additional corrections are performed.


\subsection{Description of the Results} In this work we present a comprehensive study for convergence of nonlocal solutions to classical counterparts, which holds for the vector-valued linearized peridynamic solid model (LPS) \cite{emmrich2007well}, a prototypical state-based model, with nonhomogeneous Dirichlet-type boundary conditions.  The vector-valued system is notoriusly difficult to study even in classical elasticity due to the cross-interactions between solution components, so well-posedness and regularity results rely heavily on Korn-type inequalities \cite{mengesha2014nonlocal}. Nonlocal versions of these tools are employed here to study how one may optimally impose nonlocal boundary conditions in order to ensure a high degree of compatibility with the local system.

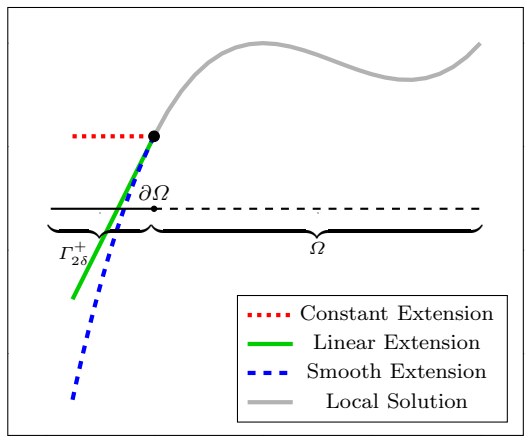
\begin{figure}\label{fig:extensions}
\begin{center}
\begin{tikzpicture}[scale=1]
\begin{axis}[legend pos=south east,
    tickwidth=0,xticklabels={,,},yticklabels={,,}
]
  \addplot[domain=-1.75:-1, red, ultra thick,dotted] {0.5};
  \addplot[domain=-1.75:-1, green!80!black, ultra thick] {0.5+10.5*(x+1)};
   \addplot[domain=-1.75:-1, blue, ultra thick,dashed] {5 -3*x^2+1.5*x^3};
   
    \addplot[domain=-1:2, black!30, ultra thick] {5 -3*x^2+1.5*x^3};
      \addlegendentry{\small{Constant Extension}};
  \addplot [fill] (-1,0.5) circle [radius=2pt];
  \addlegendentry{\small{Linear Extension}}
        \addlegendentry{\small{Smooth Extension}};
    \addlegendentry{\small{Local Solution}}
     \addplot[domain=-1.95:-1, black, thick] {-3};
      \addplot[domain=-1:2, black, thick, dashed] {-3};
         \addplot[fill] (-1,-3) circle (1pt)node[above,font=\small] {$\partial \Omega$};
            \addplot[fill] (-1.5,-3.2)node[below] {$\underbrace{\hspace{39pt}}_{\Gamma_{2\delta}^+\quad\quad}$};
         \addplot[fill] (0.5,-3.2)node[below] {$\underbrace{\hspace{125pt}}_{\Omega}$};
\end{axis}

\end{tikzpicture}
\end{center}
\caption{Different types of extensions of a local solution to the nonlocal boundary layer 
$ \Gamma_{2\delta}^+$: continuous (constant), $C^1$ (linear), $C^2$ (quadratic/smooth)}
\end{figure}

The analysis performed in this paper does not prescribe a particular method for imposing boundary conditions for the nonlocal system. Instead, we establish bounds for the difference between the local and nonlocal solutions, in terms of the difference (on an interior collar) between the nonlocal data and extensions of different degrees of regularity (even fractional) for the classical solution. The construction of different possible extensions is illustrated in Figure \ref{fig:extensions}, where $\Gamma_{2\delta}^+$ is an exterior collar (see the definition \eqref{notationdomain1}) that nonlocal boundary conditions will be prescribed to guarantee the well-posedness of the nonlocal equation. 

To provide a more precise description, suppose that $\{\vu_\delta\}_{\delta>0}\subseteq L^2(\dom;\re^d)$ is a family of solutions to the nonlocal problems parameterized by $\delta$. Suppose that $\ov{\vu}_0\in W^{2,2}(\re^d;\re^d)$ is an extension of the strong solution $\vu_0\in W^{2,2}(\dom;\re^d)$ to the corresponding local problem. The theoretical bound for the rate of convergence of $\|\vu_\delta-\vu_0\|_{L^2(\dom;\re^d)}$ is decomposed into two components. One component captures the rate of convergence away from the local boundary and depends on the regularity above $W^{2,2}$ that the local solution possesses. The other component captures the rate of convergence near the local boundary $\partial\dom$. The convergence near the boundary, in turn, depends upon an exterior and interior term. More specifically, the exterior term is provided by a priori rates of convergence of the prescribed collar values on $\bnd^+_{2\delta}$ to the extended local solution as $\partial\dom$ is approached. Knowledge, if any, about the convergence of $\vu_\delta$ to $\vu_0$ on $\bnd^-_{2\delta}$ can also be incorporated to improve the theoretical rate of convergence.

In Table~\ref{fig:rates}, we assume the extension of the local solution has bounded fourth-order derivatives (quantified by $\alpha=2$, see Theorem \ref{thm:conv}). Thus, the convergence away from the boundary is of order-two, and it is the convergence rate near the boundary that dominates the theoretical rate. The value of $\beta\ge0$ identifies the rate at which the collar-values approach $\ov{\vu}_0$ as $\partial\dom$ is approached. If $\beta=0$, then the function values have order-one convergence, which can be ensured, for example, with a constant extension of the local boundary values. If $\beta=1$ then the first-order derivatives of the collar-values and $\ov{\vu}_0$ have order-one convergence, which can be provided by a linear extension of the local boundary values. If the prescribed collar-values coincide with $\ov{\vu}_0$ on $\bnd^+_{2\delta}$, then $\beta$ is effectively infinite-- there is no boundary component for the lower bound on the theoretical rate of convergence (RC). We use $\gamma$ to reflect a priori information, if any, about the rate of convergence of $\vu_\delta$ to $\vu_0$ as $\partial\dom$ is approached from the interior. If it is only known that $\{\vu_\delta\}_{\delta>0}$ are uniformly bounded, then $\gamma=\frac{1}{2}$. In the table we use $\gamma=\frac{1}{2}$ for $\beta=0$ and $\gamma=\frac{5}{2}$ for $\beta=1$. These choices are based on numerical observations in Section \ref{sec:numerics}.  If $\|\vu_\delta-\vu\|_{L^2(\bnd^-_{2\delta})}$ has order-one convergence to zero, then $\gamma=\frac{3}{2}$. See Theorem \ref{thm:conv} for the result in its full generality and notation.

\begin{table}\label{fig:rates}
\begin{center}
\begin{tabular}{|c|c|c|c|c|c|c|}
\hline
Type of extension & $\alpha$& $\beta$ & $\gamma$ & Theoretical RC& Numerical RC\\
\hline\hline
Constant &$2$&$0$& $3/2$ & $\delta^{1/2}$ & $\delta$ \\
\hline
Linear $C^1$ &$2$&$1$&$5/2$ & $\delta^{3/2}$& $\delta^{2}$\\
\hline
Smooth $C^\infty$ & $2$&$\infty$&-- & $\delta^{2}$& $\delta^{2}$\\
\hline
\end{tabular}
\end{center}
\caption{Theoretical and numerical rates of convergence (RC) for the nonlocal solution to its classical counterpart. $\alpha$, $\beta$, and $\gamma$ are as defined in Theorem \ref{thm:conv}.}
\end{table}

\subsection{Significance/Contributions}

The study of convergence of nonlocal solutions to classical counterparts is performed here by taking into consideration 
\begin{enumerate}[(i)] 
\item how well the nonlocal boundary condition matches with a (smooth) extension of the classical solution;  
\item the degree of smoothness for the extension. 
\end{enumerate}
Note that, in contrast with other works, the analysis allows one to consider {\it general} nonlocal boundary conditions (no a priori construction needs to be made), then based on the compatibility between the nonlocal data and the extension to classical problem, and the availability of additional bounds, an explicit rate of convergence is derived. This flexibility in implementations allows considerations of different types of nonlocal boundary data, which may be selected to best fit other desired features of the model (e.g. to conserve some physical property of the boundary data, or to take into account interactions with other dynamics in a coupled system).

The explicit decay rates obtained in Theorem \ref{thm:conv} are based on carefully analyzing bounds on different sets (interior, exterior collars etc.), thus also identifying the bottleneck in obtaining faster convergence. Thus, the availability of improved bounds on the interior collar will increase the rate of convergence. However, a final verdict on whether these (assumed) bounds are expected to hold for irregular boundaries or solutions, has not been delivered, so investigations into the optimality of these results are forthcoming.

As it was mentioned above, the results apply to vector valued systems, a setting for which a comparison principle is not available. This is in  contrast to diffusion problems, where the availability of comparison principles allows a more direct method of capturing the difference between classical and nonlocal problems. 

The results of this paper are proven for general kernels (of different singularities), so they are applicable for a a variety of nonlocal systems where the profiles of interaction can model different media or behaviors. Additionally, the main argument does {\it not} use any information on 
\begin{enumerate}[(a)]
    \item the shape of the domain; the regularity of the boundary may only (indirectly) affect the smoothness of the extensions. 
    \item the regularity/bounds on the nonlocal solution inside the domain, beyond $L^2$ integrability. 
\end{enumerate}

Finally, the convergence bounds transfer to the dynamic case, with the only significant difference that the constants may increase exponentially in time (but are bounded for finite times).


\subsection{Paper Organization} The background material together with the notation are introduced in Section \ref{sec:intro}. The linear peridynamic solid (LPS) model is introduced in Section \ref{sec:LPS}, together with the notation for the nonlocal operators and domains used throughout. Also in Section \ref{sec:LPS}, the nonlocal boundary value problem with Dirichlet-type constraints is introduced, for which wellposedness follows by employing a Poincar\'e-Korn type inequality. The convergence analysis for this static case is obtained in Section \ref{sec:conv}, under several different scenarios distinguishing between the types of extension (continuous, linear, or smooth) of the boundary value data, and other conditions. The dynamic system is studied in Section \ref{sec:dynamic}, where similar rates are shown. The numerical experiments are presented for both, the static and dynamic cases in Section \ref{sec:numerics}. The final Section \ref{sec:concl}  presents conclusions of this work, as well as future directions.


\section{Notation and Preliminaries}\label{sec:intro}

The results of the paper hold in the vectorial framework, for which vector and tensor operations will be required.
Given a vector space $\mathbb{V}$, we denote $\mathbb{V}^{\otimes \ell}=\underbrace{\mathbb{V}\otimes\cdots\otimes\mathbb{V}}_{\ell\text{ times}}\cong\underbrace{\mathbb{V}\times\cdots\times\mathbb{V}}_{\ell\text{ times}}$, which can be identified with $\mathbb{V}^\ell$. In the sequel, $\mathbb{V}$ will be a space of $\ell$-th order tensors, with components in $\re$, so $\mathbb{V}=\re^{d\otimes\ell}\cong\re^{d^\ell}$, for some $\ell=0,1,2,\dots$.
We use $\bdot$ for the usual vector dot product and $\re^{d\otimes 0}:=\re^1=\re^{d^0}$ and $\re^{d\otimes 1}=\re^d$. 

Let $\{\ve{e}_i\}_{i=1}^d\subseteq\re^d$ be the unit coordinate vectors in $\re^d$. Given $\vv\in\re^d$, we may write
\[
    \vv=\lp\vv\rp_{i}=\sum_{i=1}^d v_i\ve{e}_i,
    \text{ with }
    v_i=\vv\bdot\ve{e}_i,\quad\text{ for }i=1,\dots,d.
\]
More generally, given $\tA\in\re^{d\otimes\ell}$,
\[
    \tA=\lp\tA\rp_{i_1\dots i_\ell}
    =\sum_{i_1=1}^d\cdots\sum_{i_\ell=1}^d
    A_{i_1\cdots i_\ell}\ve{e}_{i_1}\otimes\cdots\otimes\ve{e}_{i_\ell},
    \quad\text{ for }i_1,\dots,i_\ell=1,\dots,d.
\]
We will use an extension of inner products and contractions to tensors of higher-orders. Suppose that $\tA\in\re^{d\otimes n}$ and $\tB\in\re^{d\otimes m}$, with $n\ge m$. The inner product of $\tA$ and $\tB$ is
\begin{align*}
    \tA\bdot\tB
    =
    \lp\tC\rp_{i_1\dots i_{n-m}=1}^d
    =
    \sum_{i_1=1}^d\dots\sum_{i_{m-n}=1}^d
        C_{i_1\dots i_{n-m}}
        \ve{e}_{i_1}\otimes\cdots\otimes\ve{e}_{i_{n-m}}
    \in\re^{d\otimes(n-m)},
\intertext{and}
    \tA\tB
    =
    \lp\te{D}\rp_{i_1\dots i_{n-m}=1}^d
    =
    \sum_{i_1=1}^d\dots\sum_{i_{m-n}}
        D_{i_1\dots i_{n-m}}
        \ve{e}_{i_1}\otimes\cdots\otimes\ve{e}_{i_{n-m}}
    \in\re^{d\otimes(n-m)},
\end{align*}
with
\[
    C_{i_1\dots i_{n-m}}
    =
    \sum_{j_1,\dots,j_m=1}^dA_{j_1\dots j_{m}i_1\dots i_{n-m}}
        B_{j_1\dots j_m}
    \:\text{ and }\:
    D_{i_1\dots i_{n-m}}
    =
    \sum_{j_1,\dots,j_m=1}^d A_{i_1\dots i_{n-m}j_1\dots j_m}B_{j_m\dots j_1}
\]
The norm of $\tA$, we define
\[
    |\tA|=\lp\tA\bdot\tA\rp^\frac{1}{2}.
\]
Cauchy's inequality implies
\[
    |\tA\tB|\le|\tA||\tB|
    \nd
    |\tA\otimes\tB|\le|\tA||\tB|.
\]

Suppose that $\tw=(\tw)_{i_1\dots i_m}\in W^{k,2}(\dom^+_{\delta_0};\re^{d\otimes m})$, so the $(k+m)$-order tensor of $k$-th-order derivatives satisfies $\partial^k\tw\in L^2(\dom^+_{\delta};\re^{d\otimes m}\otimes\re^{\otimes k})$. 
In components, 
\[
    \lp\partial^k\tw\rp_{i_1\dots i_m j_1\dots j_k}
    =\sum_{i_1,\dots,i_m}^d\sum_{j_1,\dots,j_k}^d
        \frac{\partial^k}{\partial x_{j_k}\cdots\partial x_{j_1}}w_{i_1\dots i_m}
    \ve{e}_{i_1}\otimes\cdots\otimes\ve{e}_{i_k}\otimes
    \ve{e}_{j_1}\otimes\cdots\otimes\ve{e}_{j_m}
\]
Finally, we use the notation $\ve{a}^{\otimes m}=\underbrace{\ve{a}\otimes\cdots\otimes\ve{a}}_{m\text{ times}}\in\re^{d\otimes m}$.

We denote by $\tI$ the second-order (matrix) identity tensor and $\tI_{\sym}$ is the fourth-order symmetric identity tensor. In components
\begin{align*}
    (\tI)_{ij}
    =&
    \sum_{i,j=1}^d=\lp\ve{e}_i\bdot\ve{e}_j\rp\ve{e}_i\otimes\ve{e}_j
\intertext{and}
    \lp\tI_{\sym}\rp_{i,j,k,\ell}
    =&
    \frac{1}{2}\sum_{i,j,k,\ell=1}^d
    \ls\lp\ve{e}_i\bdot\ve{e}_k\rp\lp\ve{e}_j\bdot\ve{e}_\ell\rp
        +\lp\ve{e}_i\bdot\ve\ve{e}_\ell\rp\lp\ve{e}_j\bdot\ve{e}_k\rp\rs
        \ve{e}_i\otimes\ve{e}_j\otimes\ve{e}_k\otimes\ve{e}_\ell.
\end{align*}

{\bf Nonlocal boundaries and domains.} The results will involve nonlocal boundaries with interior and exterior layers, for which we will adopt the following notation. For each $\rho>0$ define
\begin{align}
\label{notationdomain1}
    &\bnd^+_{\rho}
    :=\{\vx\in\ren\setminus\dom:0<\dist(\vx,\partial\dom)<\rho\}
    &&\nd&&
    \dom^+_{\rho}:=\dom\cup\bnd^+_\rho\\
\label{notationdomain}
    &\bnd^-_\rho
    :=\{\vx\in\dom:\dist(\vx,\partial\dom)<\rho\}
    &&\nd&&
    \dom^-_\rho  :=\dom\setminus\bnd^-_\rho.
\end{align}
Additional notation regarding operators and domains will be introduced throughout the statements and proofs below.


\section{A Linear State-Based Peridynamic Model}\label{sec:LPS}

We consider the state-based linear peridynamic solid (LPS) model in a body occupying the domain $\Omega\subset\mathbb{R}^d$. While the main results hold for any dimension $d\ge2$, we are primarily interested in $d = 2$ or $3$. Before presenting more details, we introduce some operators and bilinear forms. For each $\delta>0$, let $B_\delta(\vx)\subseteq\re^d$ denote the $\delta$-ball centered at $\vx$. Let $K_\delta:(0,\infty)\to[0,\infty)$ be a Borel-measurable function, with support contained in $[0,\delta]$, and define $\vK_\delta:\re^d\to\re^d$ by
\begin{equation}\label{eq:kernelradial}
\vK_\delta(\vz):=K_\delta(|\vz|)\vz.
\end{equation}
Following~\cite{mengesha2014nonlocal}, the weighted volume for the $\delta$-ball $B_\delta=B_\delta(\ve{0})$ is
\begin{equation}\label{E:WeightedVol}
    m(\delta)=\int_{B_\delta}\vK_\delta(\vz)\bdot\vz\dd\vz
    =\int_{B_\delta}K_\delta(|\vz|)|\vz|^2\dd\vz,
\end{equation}

\subsection{Nonlocal operators and bilinear functionals}
For the above kernel $\vK$ and its weighted volume $m(\delta)$ we introduce the nonlocal dilatation of a map $\vu\in L^2(\re^d;\re^d)$ by
\begin{equation}\label{eq:dilop}
    \theta_\delta(\vx;\vu)=\frac{d}{m(\delta)}\int_{B_\delta}\vK_\delta(\vz)\bdot\vu(\vx+\vz)\dd\vz.
\end{equation}
Next we introduce the main component operators for the linear Navier-system. On the space $L^2(\re^d;\re^d)$ we define $\oA_\delta,\oB_\delta$ by 
\begin{equation}\label{eq:opA}
    \oA_{\delta}\vu(\vx)
    :=
    \frac{C_\oA}{m(\delta)}\int_{B_\delta}\vK_\delta(\vz)
        \theta_\delta(\vx+\vz;\vu)\dd\vz
\end{equation}
and
\begin{equation}\label{eq:opB}
    \oB_{\delta}\vu(\vx):=\frac{C_{\oB}}{m(\delta)}
    \int_{\bll_\delta}\ls\frac{\vK_\delta(\vz)\otimes\vz}{|\vz|^2}\rs(\vu(\vx+\vz)-\vu(\vx))\dd\vz.
\end{equation}
The nonlocal Navier-operator $\oL_\delta$ defined on $L^2(\re^d;\re^d)$ is given by
\begin{equation}\label{E:NonlocalEq}
    \oL_\delta\vu(\vx):=-(\lambda-\mu)\oA_\delta\vu(\vx)-\mu\oB_\delta\vu(\vx).
\end{equation}
The ranges of the operators  $ \oA_{\delta},  \oB_{\delta}, \oL_\delta$ introduced above depend on the integrability/smoothness properties of $\vK_\delta$. The scaling constants $C_{\oA},C_{\oB}>0$ are chosen such that
\begin{equation}\label{E:KernelAss}
    \frac{C_{\oA}}{m(\delta)}\int_{B_\delta}\vK_\delta(\vz)\otimes\vz\dd\vz=\tI
    \nd
    \frac{C_{\oB}}{m(\delta)}
    \int_{B_\delta}\frac{\vK_\delta(\vz)\otimes\vz^{\otimes3}}{|\vz|^2}\dd\vz
    =\tI\otimes\tI+2\tI_{\sym}.
\end{equation}
For $\vu:\re^d\to\re^d$ sufficiently smooth, it will be shown below in Theorem  \ref{T:thetaConv} that $\oA_\delta\vu\to\oA_0 \vu$ and $\oB_\delta \vu\to\oB_0\vu$, with the local Navier-operator $\oL_0$ given by
\begin{equation}\label{E:LocalEq}
    \oA_0\vu(\vx)=\partial\diverge\vu(\vx)
    \nd
   \oB_0\vu(\vx)=\partial\diverge\vu(\vx)+\diverge\partial\vu(\vx)
\end{equation}
and thus $\oL_\delta\vu\to\oL_0\vu$, with
\begin{equation}\label{E:LocalEq1}
    \oL_0\vu(\vx)
    =
    -(\lambda-\mu)\oA_0\vu(\vx)-\mu\oB_0\vu(\vx)
    =
    -\lambda\partial\diverge\vu(\vx)-\mu\diverge\partial\vu(\vx).
\end{equation}

Two popular choices for the kernels are $K_\delta^1(r)=\Lambda_\delta$ and $K_\delta^2(r)=\Lambda_\delta/r$. For these exmaples, the parameters for 3D linear elasticity are $C_{\oA}=3$, $C_{\oB}=30$ and
\begin{equation}\label{eqn:K3}
\dfrac{K^1_\delta(\rho)}{m(\delta)}=
\left\{\begin{array}{cl}
\dfrac{5}{4\pi\delta^5},\; &\text{ for }\rho\leq \delta; \\
0,\;& \text{ for }\rho> \delta; \\
\end{array}
\right.
\quad {\rm and} \quad
\dfrac{K_\delta^2(\rho)}{m(\delta)}=\left\{\begin{array}{cl}
\dfrac{1}{\pi\delta^4\rho},\; &\text{ for }\rho\leq \delta; \\
0,\;& \text{ for }\rho> \delta. \\
\end{array}
\right.
\end{equation}  
For 2D problems, $C_{\oA}=2$, $C_{\oB}=16$.
\begin{equation}\label{eqn:K}
\dfrac{K^1_\delta(\rho)}{m(\delta)}=
\left\{\begin{array}{cl}
\dfrac{2}{\pi\delta^4},\; &\text{ for }\rho\leq \delta; \\
0,\;& \text{ for }\rho> \delta; \\
\end{array}
\right.
\quad {\rm and} \quad
\dfrac{K_\delta^2(\rho)}{m(\delta)}=\left\{\begin{array}{cl}
\dfrac{3}{2\pi\delta^3\rho},\; &\text{ for }\rho\leq \delta; \\
0,\;& \text{ for }\rho> \delta. \\
\end{array}
\right.
\end{equation}

Given $\vv,\vw\in L^2(\re^d;\re^d)$ satisfying $\vv=\vw=\ve{0}$ on $\re^d\setminus\dom$, define the bilinear forms $W_{\alpha,\delta}(\vv,\vw)$ and $W_{\beta,\delta}(\vv,\vw)$ by
\begin{equation}\label{eq:wa1}
    W_{\oA,\delta}(\vv,\vw)
    =
    \int_{\dom}\int_{B_\delta}\ls\oA_\delta\vv(\vx)\rs\vw(\vx)\dd\vx
\end{equation}
and
\begin{equation}\label{eq:wb}
    W_{\oB,\delta}(\vu,\vv)=
    \int_{\dom}\int_{B_\delta}\ls\oB_\delta\vv(\vx)\rs\vw(\vx)\dd\vx.
\end{equation}
Using Fubini's theorem and a change of variables, we may rewrite
\begin{equation}\label{eq:wa2}
    W_{\oA,\delta}(\vv,\vw)
    =
    \frac{1}{d}\int_{\dom}\theta_\delta(\vx;\vv)\theta_\delta(\vx;\vw)\dd\vx
\end{equation}
and
\begin{equation}\label{eq:wb2}
    W_{\oB,\delta}(\vv,\vw)
    =
    \frac{1}{2m(\delta)}\int_{\dom}\int_{B_\delta}
    K_\delta(|\vz|)\lp\frac{\ls\vv(\vx+\vz)-\vv(\vx)\rs\bdot\vz}{\vz}\rp
        \lp\frac{\ls\vw(\vx+\vz)-\vw(\vx)\rs\bdot\vz}{\vz}\rp\dd\vz\dd\vx.
\end{equation}
Then the bilinear form associated with $\oL_\delta$ is
\begin{equation}\label{eq:wl}
    W_\delta(\vv,\vw)
    =C_{\oA}(\lambda-\mu)W_{\oA,\delta}(\vv,\vw)
    +C_{\oB}\mu W_{\oB,\delta}(\vv,\vw),
\end{equation}
and the nonlocal strain energy is $W_{\delta}(\vv,\vv)$.


\subsection{Well-posedness for the Dirichlet-type Constraint Problem}\label{sec:d}

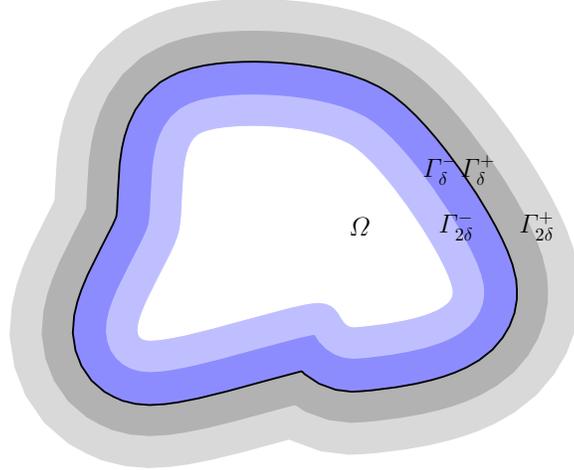
\begin{figure}\label{fig:layers}
\begin{center}
\scalebox{.69}{
\begin{tikzpicture}[scale=1.55]
\draw[line width=140pt,black!15] plot [smooth cycle] coordinates{(-0.7,-0.7) (1.45,-.25) (1.95, -0.55) (3.15,-0.122) (1.9,1.7) (0.02,1.85) (-0.25,0.6) };
\draw[line width=105pt,black!30] plot [smooth cycle] coordinates{(-0.7,-0.7) (1.45,-.25) (1.95, -0.55) (3.15,-0.122) (1.9,1.7) (0.02,1.85) (-0.25,0.6) };
\draw[line width=72pt,black] plot [smooth cycle] coordinates
 {(-0.7,-0.7) (1.45,-.25) (1.95, -0.55) (3.15,-0.122) (1.9,1.7) (0.02,1.85) (-0.25,0.6) };
\draw[line width=70pt,blue!45] plot [smooth cycle] coordinates{(-0.7,-0.7) (1.45,-.25) (1.95, -0.55) (3.15,-0.122) (1.9,1.7) (0.02,1.85) (-0.25,0.6) };
\draw[line width=35pt,blue!25] plot [smooth cycle] coordinates
 {(-0.7,-0.7) (1.45,-.25) (1.95, -0.55) (3.15,-0.122) (1.9,1.7) (0.02,1.85) (-0.25,0.6) };
 \draw[line width=0pt,fill=white,white] plot [smooth cycle] coordinates
 {(-0.7,-0.7) (1.45,-.25) (1.95, -0.55) (3.15,-0.122) (1.9,1.7) (0.02,1.85) (-0.25,0.6) };
  \node at (2,0.7) {\Large$\Omega$};
     \node at (3.0,1.4) {\Large$\Gamma_{\delta}^-$};
        \node at (3.47,1.4) {\Large$\Gamma_{\delta}^+$};
    \node at (3.2,0.7) {\Large$\Gamma_{2\delta}^-$};
        \node at (4.2,0.7) {\Large$\Gamma_{2\delta}^+$};
  \end{tikzpicture}}
\caption{Domain $\Omega$ with adjacent layers. Interior layers: $\Gamma^-_{\delta}, \Gamma^-_{2\delta}\subset \Omega$;  Exterior layers: $\Gamma^+_{\delta}, \Gamma^+_{2\delta}\subset \re^d\setminus \Omega$}
\end{center}
\end{figure}

In order to define the nonlocal system with a Dirichlet-type boundary condition, we use the notation introduced in \eqref{notationdomain1}-\eqref{notationdomain} for the domain $\Omega$ with two interior collars and two exterior collars, as illustrated in Figure \ref{fig:layers}.
We now consider a state-based peridynamic problem with Dirichlet-type boundary condition: 

\begin{equation}\label{eqn:probd}
\left\{\begin{array}{ll}
    \ds- \frac{C_{\oA}}{m(\delta)}  \int_{B_\delta (\vx)}\left(\lambda - \mu\right) K_\delta(\left|\vy-\vx\right|) \left(\vy-\vx \right)\left(\theta(\vx) + \theta(\vy) \right) d\vy&\\
    \ds\qquad-  \frac{C_\beta}{m(\delta)}\int_{B_\delta (\vx)} \mu K_\delta(\left|\vy-\vx\right|)\frac{\left(\vy-\vx\right)\otimes\left(\vy-\vx\right)}{\left|\vy-\vx\right|^2} \cdot \left(\vu(\vy) - \vu(\vx) \right) d\vy = \vf(\vx),\quad &\text{ for $\vx$ in }\Omega,\\
\ds\theta(\vx)=\dfrac{d}{m(\delta)}\int_{B_\delta (\vx)} K_\delta(\left|\vy-\vx\right|) (\vy-\vx)\cdot \left(\vu(\vy) - \vu(\vx) \right)d\vy,\quad &\text{ for $\vx$ in }\omg^+_{\delta},\\
\vu(\vx)=\vu_D(\vx),\quad &\text{ for $\vx$ in }\bnd^+_{2\delta}.
\end{array}\right.
\end{equation}

The well-posedness of the above system follows from the Korn-Poincar\'{e} type inequality below which is an application of Proposition 5 in \cite{mengesha2014nonlocal}.
\begin{theorem}\label{T:PKorn}
Assume $\Omega\subset\re^d$ is a bounded domain with a suffiicently smooth boundary (in particular, it satisfies the interior cone condition). Then there exists a $\delta_0>0$ and a constant $C_{\PK}<\infty$ such that for each $0<\delta\le\delta_0$,
\[
    \|\vv\|_{L^2(\dom;\re^d)}
    \le
    C_{\PK}W_{\oB,\delta}(\vv,\vv),
    \frl
    \vv\in\{\vw\in L^2(\re^d):\vw(\vx)=\ve{0}\text{ a.e. }\vx\in\re^d\setminus\dom\}.
\]
\end{theorem}






With existence of solutions to the nonlocal problem \eqref{eqn:probd} established we will study in the next section convergence of nonlocal solutions to the classical counterparts. 

\section{Convergence Analysis}\label{sec:conv}

We will now show the truncation error of the Dirichlet-type volume constraint formulation. We denote by $\vu_\delta$ the solution of the nonlocal problem \eqref{eqn:probd} and $\vu_0$ as the solution of the classical elastic problem:
 \begin{equation}\label{eqn:local}
 \left\{\begin{array}{ll}
 -\nabla\cdot(\lambda tr(\mathbf{E})\mathbf{I}+2\mu \mathbf{E})=\vf,\quad \text{where }\mathbf{E}=\dfrac{1}{2}(\nabla \vu+(\nabla \vu)^T),&\quad \text{in }\Omega,\\
 \vu= \vu_D,&\quad \text{on }\partial\Omega.
 \end{array}\right.
 \end{equation}

\subsection{Operator Convergence}

For an integer $k>0$ let $\tw\in C^k(\re^d;\re^{d\otimes m})$. Then 
Taylor's formula, with remainder, can be written as
\begin{equation}\label{E:TaylorExp}
    \tw(\vx+\vz)
    =\tw(\vx)+
    \sum_{j=1}^k\frac{1}{j!}\partial^j\tw(\vx)\vz^{\otimes j}
    +\set{R}_k(\tw(\vx);\vz)\vz^{\otimes k}.
\end{equation}
Here
\[
    \set{R}_k(\tw(\vx);\vz):=
    \frac{1}{(k-1)!}\lp\int_{0}^1(1-s)^{k-1}
        \ls\partial^k\tw(\vx+s\vz)-\partial^k\tw(\vx)\rs\dd s\rp.
\]
We note that $j$ is not a multi-index. Observe that $\oR_k$ is linear and
\begin{align}
\nonumber
    \oR_k(\tw(\vx+\vz_1);\vz_2)-\oR_k(\tw(\vx);\vz_2)
    &=
    \oR_k(\partial_{\vx}\tw(\vx);\vz_2)\vz_1
        +\oR_k(\oR_1(\tw(\vx);\vz_1)\vz_1;\vz_2)\\
\label{E:RemComp}
    &=
    \ls\oR_k(\partial_{\vx}\tw(\vx);\vz_2)
        +\oR_1\lp\oR_k(\tw(\vx);\vz_2);\vz_1\rp\rs\vz_1
\end{align}

We will need the following pointwise representation for mappings in fractional Sobolev functions. It provides a generalization of the Hajłasz-type pointwise representation for functions in a 1st-order Sobolev space~\cite{Haj:96a}. The argument is nearly identical to the proof for Theorem 1.1, part (2), in~\cite{Hu:03a} and is included for the sake of completeness.
\begin{theorem}
Let $0<\alpha<1$ and $0<\delta_0\le1$ and $0<\delta'<\delta_0$ be given. Suppose that $\tw\in W^{\alpha,2}(\dom^+_{\delta_0};\re^{d\otimes m})$, so
\[
    |\tw|^2_{W^{\alpha,2}(\dom^+_{\delta_0})}
    :=
    \int_{\dom^+_{\delta_0}}\int_{\dom^+_{\delta_0}}\frac{|\tw(\vy)-\tw(\vx)|^2}
        {|\vy-\vx|^{d+2\alpha}}\dd\vy\dd\vx<\infty.
\]
Put $c_\alpha=c_\alpha(d):=2^{d+1}\cdot 4^{\alpha}/|B_1|$ and $\delta_1:=\delta_0-\delta'$, and define $w_\alpha:\dom^+_\delta\to\re$ by
\begin{equation}\label{E:walphaDef}
    w_\alpha(\vx):=\lp\int_{\dom^+_{\delta_0}\cap B_2(\vx)}
    \frac{|\tw(\vy)-\tw(\vx)|^2}{|\vy-\vx|^{d+2\alpha}}
        \dd\vy\rp^\frac{1}{2}.
\end{equation}
Then $w_\alpha\in L^2(\dom^+_{\delta_0})$ and
\begin{equation}\label{E:SobolevRep}
    |\tw(\vy)-\tw(\vx)|^2\le c_\alpha|\vy-\vx|^{2\alpha}\lp w_\alpha(\vx)^2+w_\alpha(\vy)^2\rp,
    \quad\text{ for all }\vx\in\dom_{\delta'}\text{ and }\vy\in B_{\delta_1}(\vx).
\end{equation}
\end{theorem}
\begin{proof}
It is clear that $\|w_\alpha\|_{L^2(\dom\cup\bnd^+_{\delta})}\le|\tw|_{W^{\alpha,2}(\dom^+_{\delta_0};\re^{d\otimes m})}$. Let $\vx\in\dom_{\delta'}$ and $\vy\in B_{\delta_1}(\vx)\subseteq\dom^+_{\delta_0}$ be given. Put $r:=|\vy-\vx|<\delta_1$. Observe that $\bll_{r}(\vx) \subseteq\dom^+_{\delta_0}\cap\bll_{2r}(\vy)$, so
\begin{align*}
    |\tw(\vy)-\tw(\vx)|^2
    =&
    \frac{1}{|B_r|}\int_{B_r(\vx)}|\tw(\vy)-\tw(\vx)|^2\dd\vx'\\
    \le&
    \frac{2}{|B_r|}\lp\int_{B_r(\vx)}|\tw(\vy)-\tw(\vx')|^2\dd\vx'
        +\int_{B_r(\vx)}|\tw(\vx')-\tw(\vx)|^2\dd\vx'\rp\\
    \le&
    \frac{2}{|B_r|}\lp\int_{\dom^+_{\delta_0}\cap B_{2r}(\vy)}|\tw(\vy)-\tw(\vx')|^2\dd\vx'
        +\int_{B_r(\vx)}|\tw(\vx')-\tw(\vx)|^2\dd\vx'\rp\\
    \le&
    \frac{2r^{d+2\alpha}}{|B_r|}\lp2^{d+2\alpha}\int_{\dom^+_{\delta_0}\cap B_{2r}(\vy)}
        \frac{|\tw(\vy)-\tw(\vx')|^2}{|\vy-\vx'|^{n+2\alpha}}\dd\vx'
        +\int_{B_r(\vx)}\frac{|\tw(\vx')-\tw(\vx)|^2}{|\vx'-\vx|^{d+2\alpha}}\dd\vx'\rp\\
    \le&
    \lp\frac{2^{d+1+2\alpha}}{|B_1|}\rp r^{2\alpha}\lp w_\alpha(\vy)^2+w_\alpha(\vx)^2\rp\\
    =&
    c_\alpha|\vy-\vx|^{2\alpha}\lp w_\alpha(\vy)^2+w_\alpha(\vx)^2\rp.
\end{align*}
\qed
\end{proof}

\begin{lemma}\label{L:RCont}
Let $k\ge1$ and $0\le\alpha\le1$ be given. Suppose that $\tv\in W^{k+\alpha,2}(\dom^+_{\delta_0};\re^{d\otimes m})$. Let $0<\delta'<\delta_0$ be given, and put $\delta_1:=\delta_0-\delta'$.
\begin{enumerate}[(a)]
\item If $\alpha=0$, then for each $\vep>0$, there exists $0<\delta\le\delta_1$ such that
\[
    \|\oR_k(\tv(\cdot);\vz)\|_{L^2(\dom^+_{\delta'})}<\vep,\quad\text{ for all $\vz\in B_\delta$}.
\]
\item If $0<\alpha<1$, then for each $0<\delta\le\delta_1$, we have
\[
    \|\oR_k(\tv(\cdot);\vz)\|_{L^2(\dom^+_{\delta'})}
    \le
    \frac{2c_\alpha}{(k-1)!}|\vz|^\alpha|\partial^k\tv|
            _{W^{\alpha,2}(\dom_{\delta_0}^+;\re^{d\otimes m}\otimes\re^{\otimes k})},
    \quad\text{ for all $\vz\in B_\delta$}.
\]
\item If $\alpha=1$, then for each $0<\delta\le\delta_1$, we have
\[
    \|\oR_k(\tv(\cdot);\vz)\|_{L^2(\dom^+_{\delta'})}
    \le
    \frac{1}{(k-1)!}|\vz|
        |\partial^{k+1}\tv|
            _{L^2(\dom_{\delta_0}^+;\re^{d\otimes m}\otimes\re^{\otimes (k+1)})},
    \quad\text{ for all $\vz\in B_\delta$}.
\]
\end{enumerate}
\end{lemma}

\begin{proof}
With a standard density argument we may assume $\tv\in C^\infty(\dom_{\delta_0}^+;\re^{d\otimes m})$. Part (a) is an immediate consequence of the definition of $\oR_k(\vv(\cdot);\vz)$ and the continuity for the translation operator in $L^2(\re^d)$. Note that if $\vx\in\dom^+_{\delta'}$ and $\vz\in\bll_{\delta}$, then $\vx+\vz\in\dom_{\delta_0}$

For part (b), let $\vz\in\bll_\delta$. We use Minkowski's integral inequality and~\eqref{E:SobolevRep} with \[
\tw=\partial^k\tv\in W^{\alpha,2}(\dom^+_{\delta_0};\re^{d\otimes m}\otimes\re^k)
\]
as follows:
\begin{align}
\label{E:ukCont1}
    \|\oR_k(\tv;\vz)\|^2_{L^2(\dom^+_{\delta'};\re^{d\otimes m}\otimes\re^{\otimes k})}
    \le&
    \frac{1}{(k-1)!}
    \int_{\dom^+_{\delta'}}\int_0^1
        \left|\partial^k\tv(\vx+s\vz)-\partial^k\tv(\vx)\right|^2\dd s\dd\vx\\
\nonumber
    \le&
    \frac{c_\alpha}{(k-1)!}\int_{\dom^+_{\delta'}}\int_0^1
        s^{2\alpha}|\vz|^{2\alpha}\lp w_\alpha(\vx)^2+w_\alpha(\vx+s\vz)^2\rp\dd s\dd\vx\\
\label{E:ukHoldCont1}
    \le&
    \frac{c_\alpha}{(k-1)!}|\vz|^{2\alpha}\lp
        \int_{\dom^+_{\delta_0}} w_\alpha(\vx)^2\dd\vx+\int_0^1\int_{\dom^+_{\delta'}} w_\alpha(\vx+s\vz)^2\dd\vx\dd s\rp
\end{align}
Recalling the definition of $w_\alpha$ in~\eqref{E:walphaDef}, with $\tw=\partial^k\tv$, the first integral is bounded by $|\partial^k\tv|_{W^{\alpha,2}}$. For the second integral, given $0<s<1$,
\begin{align*}
    \int_\dom w_\alpha(\vx+s\vz)^2\dd\vx
    =&
    \int_{\dom^+_{\delta'}}\int_{\dom_{\delta_0}^+\cap\bll_2(\vx+s\vz)}
        \frac{\left|\partial^k\tv(\vy)-\partial^k\tv(\vx+s\vz)\right|^2}
        {|\vy-(\vx+s\vz)|^{d+2\alpha}}
        \dd\vy\dd\vx\\
    \le&
    \int_{\dom^+_{\delta'}}\int_{\dom_{\delta_0}^+}
        \frac{\left|\partial^k\tv(\vy)-\partial^k\tv(\vx+s\vz)\right|^2}
        {|\vy-(\vx+s\vz)|^{d+2\alpha}}
        \dd\vy\dd\vx\\
    \le&
    \int_{\dom_{\delta_0}^+}\int_{\dom_{\delta_0}^+}
        \frac{\left|\partial^k\tv(\vy)-\partial^k\tv(\vx')\right|^2}{|\vy-\vx'|^{d+2\alpha}}
        \dd\vx'\dd\vy
    =|\partial^k\tv|_{W^{\alpha,2}(\dom_{\delta_0}^+;\re^{d\otimes m}\otimes\re^{\otimes k})}.
\end{align*}
For the last line, we used Fubini's theorem and the change of variables $\vx+s\vz\mapsto\vx'$. Again, we note that $\vx\in\dom^+_{\delta'}$ and $\vz\in\bll_\delta$ implies $\vx+s\vz\in\dom^+_{\delta_0}$. Incorporating these bounds into~\eqref{E:ukHoldCont1} completes the proof.

Part (c) follows from~\eqref{E:ukCont1} and the fundamental theorem of calculus.
\qed
\end{proof}

\subsubsection{Operator Convergence for $\theta_\delta$}

\begin{theorem}\label{T:thetaConv}
    Let $0\le\alpha\le 3$ and $0<\delta'<\delta_0$ be given. Suppose that $\vu\in W^{1+\alpha,2}(\dom^+_{\delta_0};\re^d)$, and put $\delta_1:=\delta_0-\delta'$.
\begin{enumerate}[(a)]
    \item If $\alpha=0,1,2$, then for each $\vep>0$, there exists  $0<\delta\le\delta_1$ such that
\[
    \left\|\theta_\delta(\cdot;\vu)-\diverge\vu\right\|_{L^2(\dom^+_{\delta'})}
    <\vep\delta^{\alpha}.
\]
    \item If $0<\alpha\le3$ and $\alpha\neq 1,2$, then and a constant $\Lambda_\alpha$, independent of $\delta$ and $\vu$, such that
\[
    \left\|\theta_\delta(\cdot;\vu)-\diverge\vu\right\|_{L^2(\dom^+_{\delta'})}
    <\Lambda_\alpha|\partial\vu|_{W^{1+\alpha}(\dom^+_{\delta_0};\re^d\otimes\re)}\delta^\alpha,
    \quad\text{ for all }0<\delta\le\delta_1,
\]
with
\[
    \Lambda_\alpha:=\lc\begin{array}{ll}
        2c_{\alpha-1}, & 0<\alpha<3\text{ and }\alpha\neq 1,2,\\
        1, &\alpha=3.
    \end{array}\rper
\]
\end{enumerate}
\end{theorem}
\begin{proof}
For each $0\le s\le 1$, $\vx\in\dom^+_{\delta'}$, and $\vz\in B_\delta$, we find $\vx+s\vz\in\dom^+_{\delta_0}$. Using~\eqref{E:TaylorExp},
\[
    \vu(\vx+\vz)-\vu(\vx)=\partial\vu(\vx)\vz+\oR_1(\vu(\vx);\vz)\vz.
\]
Incorporating this into the definition of $\theta_\delta$, we may write
\begin{align}
\nonumber
    \lp\frac{m(\delta)}{d}\rp\theta_\delta(\vx;\vu)
    =&
    \int_{B_\delta}\vK_\delta(\vz)\bdot\ls\vu(\vx+\vz)-\vu(\vx)\rs\dd\vz\\
\nonumber
    =&
    \int_{B_\delta}\vK_\delta(\vz)\bdot\ls\partial\vu(\vx)\vz\rs\dd\vz
    +\int_{B_\delta}\vK_\delta(\vz)\bdot\ls\oR_1(\vu(\vx);\vz)\vz\rs\dd\vz\\
\label{E:thetaExp1}
    =&
    \lp\int_{B_\delta}\vK_\delta(\vz)\otimes\vz\dd\vz\rp\bdot\partial\vu(\vx)
    +\int_{B_\delta}\big[\vK_\delta(\vz)\otimes\vz\big]\bdot\oR_1(\vu(\vx);\vz)\dd\vz.
\end{align}
For the first equality, we used the antisymmetry of $\vK_\delta$. Recalling~\eqref{E:WeightedVol} provides
\begin{align*}
    \frac{d}{m(\delta)}
    \lp\int_{B_\delta}\vK_\delta(\vz)\otimes\vz\dd\vz\rp\bdot\partial\vu(\vx)
    =&
    \frac{1}{m(\delta)}\lp\int_{B_\delta}K_\delta(|\vz|)|\vz|^2\dd\vz\rp\tI
        \bdot\partial\vu(\vx)
    =\tI\bdot\partial\vu(\vx)\\
    =&
    \diverge\vu(\vx).
\end{align*}
Thus
\begin{align}
\nonumber
    \left\|\theta_\delta(\cdot;\vu)-\diverge\vu\right\|_{L^2(\dom)}
    \le&
    \frac{d}{m(\delta)}\lp\int_\dom\lp\int_{B_\delta}
        |\vK_\delta(\vz)||\vz||\oR_1(\vu(\vx);\vz)|
        \dd\vz\rp^2\dd\vx\rp^\frac{1}{2}\\
\nonumber
    \le&
    \frac{d}{m(\delta)}\int_{B_\delta}\lp\int_\dom|\vK_\delta(\vz)|^2
        |\vz|^2|\oR_1(\vu(\vx);\vz)|^2\dd\vx
        \rp^\frac{1}{2}\dd\vz\\
\label{E:thetaDivBnd1}
    \le&
    \frac{d}{m(\delta)}\int_{B_\delta}|K_\delta(|\vz|)||\vz|^2
        \|\oR_1(\vu(\cdot);\vz)\|_{L^2(\dom)}\dd\vz.
\end{align}
For the second inequality, we used Minkowski's integral inequality. Recall that $\vu\in W^{1+\alpha,2}(\dom^+_{\delta_0})$. Part (a) of the theorem, with $\alpha=0$, follows from Lemma~\ref{L:RCont}(a), and part (b), with $0<\alpha<1$, follows from Lemma~\ref{L:RCont}(b).

Now suppose that $1\le\alpha<2$. We can then use the antisymmetry of $\vz\mapsto\vK_\delta(\vz)\otimes\vz^{\otimes2}=K_\delta(|\vz|)\vz^{\otimes3}$ and~\eqref{E:TaylorExp} to write
\begin{multline}\label{E:thetaExp2}
    \theta_\delta(\vx;\vu)
    =
    \diverge\vu(\vx)
    +\frac{d}{2m(\delta)}
        \underbrace{
        \lp\int_{B_\delta}\vK_\delta(\vz)\otimes\vz^{\otimes2}\dd\vz\rp}_{=\te{0}}
            \bdot\partial^2\vu(\vx)\\
        +\frac{d}{m(\delta)}
        \int_{B_\delta}\ls\vK_\delta(\vz)\otimes\vz^{\otimes2}\rs
        \bdot\oR_2(\vu(\vx);\vz)\dd\vz.
\end{multline}
Arguing as in~\eqref{E:thetaDivBnd1}, we deduce that
\begin{align}
\nonumber
    \left\|\theta_\delta(\cdot;\vu)-\diverge\vu\right\|_{L^2(\dom^+_{\delta'})}
    \le&
    \frac{2d}{m(\delta)}\int_{B_\delta}|\vK_\delta(\vz)||\vz|^2
        \|\oR_2(\vu(\cdot);\vz)\|_{L^2(\dom)}\dd\vz\\
\label{E:thetaDivBnd2}
    \le&
    \frac{2d\delta}{m(\delta)}\int_{B_\delta}|K_\delta(|\vz|)||\vz|^3
        \|\oR_2(\vu(\cdot);\vz)\|_{L^2(\dom)}\dd\vz.
\end{align}
If $\alpha=1$, then Lemma~\ref{L:RCont}(a) implies part (a) of the theorem. If $1<\alpha<2$, then Lemma~\ref{L:RCont}(b) implies part (b). For $\alpha=2$, part (b) is implied by Lemma~\ref{L:RCont}(c).

As observed above, the kernel $\vz\mapsto\vK_\delta(\vz)\otimes\vz^{\otimes2}$ is antisymmetric, so the corresponding integral in~\eqref{E:thetaExp2} is zero. Hence, if $\alpha\ge 2$, we have
\[
    \theta_\delta(\vx;\vu)
    =
    \diverge\vu(\vx)
    +\frac{d}{m(\delta)}\int_{B_\delta}
        \ls\vK_\delta(\vz)\otimes\vz^{\otimes3}\rs\bdot\oR_3(\vu(\vx);\vz)\dd\vz.
\]
For $\alpha=2$, an argument similar to that used for~\eqref{E:thetaDivBnd1} yields part (a). Arguing as in~\eqref{E:thetaDivBnd2}, we obtain part (b) when $2<\alpha<3$. For the rest of part (b), when $\alpha=3$, we apply part (c) of Lemma~\ref{L:RCont} to complete the proof of the theorem.
\qed
\end{proof}

\begin{theorem}\label{T:LConv}
    Let $0\le\alpha\le 2$ be given. Suppose that $\vu\in W^{2+\alpha,2}(\dom^+_{\delta_0};\re^d)$.
\begin{enumerate}[(a)]
    \item If $\alpha=0$ or $\alpha=1$, then for each $\vep>0$, there exists a $0<\delta\le\delta_0/2$ such that
\[
    \left\|\oA_\delta\vu-\oA_0\vu\right\|_{L^2(\dom)}
    <\vep\delta^{\alpha}
    \nd
    \left\|\oB_\delta\vu-\oB_0\vu\right\|_{L^2(\dom)}
    <\vep\delta^{\alpha}.
\]
    \item If $0<\alpha\le 2$ and $\alpha\neq 1$, then there is a constant $\Lambda'_\alpha$, independent of $\delta$ and $\vu$, such that for all $0<\delta\le\delta_0/2$,
\[
    \left\|\oA_{\delta}\vu-\oA_0\vu\right\|_{L^2(\dom)}
    <\Lambda'_{\alpha}|\vu|_{W^{2+\alpha,2}(\dom_{\delta_0})}\delta^{\alpha}
    \nd
    \left\|\oB_{\delta}\vu-\oB_0\vu\right\|_{L^2(\dom)}
    <\Lambda'_\alpha|\vu|_{W^{2+\alpha,2}(\dom_{\delta_0})}\delta^{\alpha}
\]
\end{enumerate}
\end{theorem}
\begin{proof}
The argument is similar to the one used for Theorem~\ref{T:thetaConv}. First, we prove the convergence rates for $\oA_{\delta}$, then the argument for $\oB_{\delta}$ will be discussed. Since $\vu\in W^{2+\alpha,2}(\dom^+_{\delta_0};\re^d)$, we find $g:=\diverge u\in W^{1+\alpha,2}(\dom^+_{\delta_0})$.

As $\alpha\ge0$, applying~\eqref{E:TaylorExp} provides
\begin{align*}
    \diverge\vu(\vx+\vz)
    =
    g(\vx+\vz)
    =&
    g(\vx)+\partial g(\vx)\vz+\oR_1(g(\vx);\vz)\vz\\
    =&
    g(\vx)+\partial g(\vx)\bdot\vz+\oR_1(g(\vx);\vz)\bdot\vz.
\end{align*}
First, we observe that the antisymmetry of $\vK_\delta$ implies
\[
    \int_{B_\delta}\vK_\delta(\vz_1)g(\vx)\dd\vz_1=\ve{0}
\]
and, by Fubini's theorem,
\begin{multline*}
    \int_{B_\delta}\int_{B_\delta}
        \vK_\delta(\vz_1)\lp\ls\vK_\delta(\vz_2)\otimes\vz_2\rs
            \bdot\oR_1(\vu(\vx);\vz_2)\rp
        \dd\vz_2\dd\vz_1\\
    =\lp\int_{B_\delta}\vK_\delta(\vz_1)\dd\vz_1\rp
        \int_{B_\delta}\ls\vK_\delta(\vz_2)\otimes\vz_2\rs
            \bdot\oR_1(\vu(\vx);\vz_2)
        \dd\vz_2=\ve{0}.
\end{multline*}
Suppose that $0\le\alpha<2$. Using~\eqref{E:thetaExp1},
\begin{multline*}
    \oA_{\delta}\vu(\vx)
    =
    \underbrace{\frac{C_{\oA}}{m(\delta)}\int_{B_\delta}\vK_\delta(\vz_1)g(\vx+\vz_1)\dd\vz_1}
        _{I_1(\vx)}\\
    +
    \underbrace{\frac{dC_{\oA}}{m(\delta)^2}\int_{B_\delta}\int_{B_\delta}
        \vK_\delta(\vz_1)\lp\ls\vK_\delta(\vz_2)\otimes\vz_2\rs
            \bdot\oR_1(\vu(\vx+\vz_1);\vz_2)\rp\dd\vz_2\dd\vz_1}_{=:I_2(\vx)}.
\end{multline*}
For $I_1$, we have
\begin{align*}
    I_1
    =&\frac{C_{\oA}}{m(\delta)}\int_{B_\delta}\vK_\delta(\vz_1)
        \ls g(\vx+\vz_1)-g(\vx)\rs\dd\vz_1
    =
    \frac{C_{\oA}}{m(\delta)}\int_{B_\delta}\ls\vK_\delta(\vz_1)\otimes\vz_1\rs
        \ls\partial g(\vx)+\oR_1(g(\vx);\vz_1)\rs\dd\vz_1\\
    =&
    \partial g(\vx)
    +\frac{C_{\oA}}{m(\delta)}\int_{B_\delta}\ls\vK_\delta(\vz_1)\otimes\vz_1\rs
        \oR_1(g(\vx);\vz_1)\dd\vz_1,
\end{align*}
where~\eqref{E:KernelAss} was used for the last line. Thus
\[
    \left|I_1(\vx)-\partial\diverge\vu(\vx)\right|
    \le
    \frac{C_{\oA}}{m(\delta)}\int_{B_\delta}|\vK_\delta(\vz_1)||\vz_1|
        |\oR_1(g(\vx);\vz_1)|\dd\vz_1,
    \frl\vx\in\dom.
\]
The Minkowski integral inequality yields
\[
    \left\|I_1-\partial\diverge\vu\right\|_{L^2(\dom)}
    \le
    \frac{C_{\oA}}{m(\delta)}\int_{B_\delta}
    K_\delta(\vz_1)|\vz_1|^2\|\oR_1(g(\cdot);\vz_1)\|_{L^2(\dom)}\dd\vz_1
\]
For $I_2$, define $\te{r}_1(\cdot;\vz_2):=\oR_1(\vu(\cdot);\vz_2)$. Applying~\eqref{E:RemComp}, we deduce that
\begin{align*}
    I_2(\vx)
    =&
    \frac{dC_{\oA}}{m(\delta)^2}\int_{B_\delta}\int_{B_\delta}
        \vK_\delta(\vz_1)\lp\ls\vK_\delta(\vz_2)\otimes\vz_2\rs\bdot
            \big[ \oR_1(\vu(\vx+\vz_1);\vz_2)-\oR_1(\vu(\vx);\vz_2)\big]\rp\dd\vz_2\dd\vz_1\\
    =&
    \frac{dC_{\oA}}{m(\delta)^2}\int_{B_\delta}\int_{B_\delta}
        \vK_\delta(\vz_1)\lp\ls\vK_\delta(\vz_2)\otimes\vz_2\rs\bdot
            \big[ \oR_1(\partial\vu(\vx);\vz_2)
            +\oR_1(\te{r}_1(\vx;\vz_2);\vz_1)\big]\vz_1\rp\dd\vz_2\dd\vz_1,
\end{align*}
and so
\begin{align*}
    |I_2(\vx)|
    \le&
    \frac{dC_{\oA}}{m(\delta)^2}\int_{B_\delta}\int_{B_\delta}
        |\vK_\delta(\vz_1)||\vz_1||\vK_\delta(\vz_2)||\vz_2|
            \ls
            |\oR_1(\partial\vu(\vx);\vz_2)|+
            |\oR_1(\te{r}_1(\vx;\vz_2);\vz_1)|\rs\dd\vz_2\dd\vz_1.
\end{align*}
Using, again, Minkowskis's integral inequality,
\begin{multline*}
    \|I_2\|_{L^2(\dom)}
    \le
    \frac{dC_{\oA}}{m(\delta)^2}\int_{B_\delta}\int_{B_\delta}
    K_\delta(|\vz_1|)K_\delta(|\vz_2|)|\vz_1|^2|\vz_2|^2\\
    \times
    \ls\|\oR_1(\partial\vu(\cdot);\vz_2)\|_{L^2(\dom)}+
            \|\oR_1(\te{r}_1(\cdot;\vz_2);\vz_1)\|_{L^2(\dom)}\rs\dd\vz_2\dd\vz_1.
\end{multline*}
Note that in the integrals for both $I_1$ and $I_2$, the domains of integration are $\bll_\delta\subseteq B_{\delta_0/2}$. We find $g\in W^{1+\alpha,2}(\dom_{\delta_0})$, $\partial\vu\in W^{1+\alpha,2}(\dom_{\delta_0};\re^{d\otimes 2})$, and $\te{r}_1(\cdot;\vz_2)\in W^{1+\alpha,2}(\dom_{\delta_0/2};\re^{d\otimes 2})$, for each $\vz_2\in\bll_{\delta_0/2}$.  It follows that
\begin{align*}
    &\left\|\oA_{\delta}\vu-\partial\diverge\vu\right\|_{L^2(\dom)}\\
    &\quad\le
    \|I_1-\partial\diverge\vu\|_{L^2(\dom)}+\|I_2\|_{L^2(\dom)}\\
    &\quad\le
    \frac{C_{\oA}}{m(\delta)}\int_{B_\delta}K_\delta(|\vz_1|)|\vz_1|^2
        \|\oR_1(g(\cdot);\vz_1)\|_{L^2(\dom)}\dd\vz_1\\
    &\qquad+
    \frac{dC_{\oA}}{m(\delta)^2}
    \int_{B_\delta}\int_{B_\delta}
    K_\delta(|\vz_1|)K_\delta(|\vz_2|)|\vz_1|^2|\vz_2|^2
    \ls\|\oR_1(\partial\vu(\cdot);\vz_2)\|_{L^2(\dom)}+
            \|\oR_1(\te{r}(\cdot;\vz_2);\vz_1)\|_{L^2(\dom)}\rs\dd\vz_2\dd\vz_1
\end{align*}

For $0\le\alpha<1$, the convergence for $\oA_{\delta}\vu$ follows from Lemma~\ref{L:RCont} and assumption~\eqref{E:KernelAss}.

If $1\le\alpha\le 2$, we can use the expansion in~\eqref{E:thetaExp2}. Then
\[
    \oA_{\delta}\vu(\vx)
    =
    I_1(\vx)
    +
    \underbrace{\frac{C_{\oA}}{m(\delta)^2}\int_{B_\delta}\int_{B_\delta}
        \ls\vK_\delta(\vz_1)\rs\ls\vK_\delta(\vz_2)\otimes\vz_2^{\otimes2}\rs
        \bdot\oR_2(\vu(\vx+\vz_1);\vz_2)\dd\vz_2\dd\vz_1}_{=:I_3(\vx)}
\]
With the additional differentiability of $\vu$,
\begin{align*}
    I_1(\vx)
    =&\frac{C_{\oA}}{m(\delta)}\int_{B_\delta}\vK_\delta(\vz_1)
        \ls g(\vx+\vz_1)-g(\vx)\rs\dd\vz_1\\
    =&
    \frac{C_{\oA}}{m(\delta)}\int_{B_\delta}\vK_\delta(\vz_1)
        \ls\partial g(\vx)\vz_1+\frac{1}{2}\partial^2g(\vx)\vz_1^{\otimes2}
            +\oR_2(g(\vx);\vz_1)\vz_1^{\otimes2}\rs\dd\vz_1\\
    =&
    \partial g(\vx)+\frac{C_{\oA}}{m(\delta)}\int_{B_\delta}K_\delta(|\vz_1|)\vz_1
        \ls\oR_2(g(\vx);\vz_1)\vz_1^{\otimes2}\rs\dd\vz_1,
\end{align*}
so
\[
    |I_1(\vx)-\partial\diverge\vu(\vx)|
    \le
    \delta\frac{\delta C_{\oA}}{m(\delta)}\int_{B_\delta}K_\delta(|\vz_1|)|\vz_1|^2
        |\oR_2(g(\vx);\vz_1)|\dd\vz_1.
\]
For $I_3$, as with $I_2$, we can use the antisymmetry of $\vz_1\mapsto\vK(\vz_1)$ to argue
\begin{align*}
    |I_3(\vx)|
    =&
    \frac{C_{\oA}}{m(\delta)^2}\int_{B_\delta}\int_{B_\delta}
        \vK_\delta(\vz_1)\lp\ls\vK_\delta(\vz_2)\otimes\vz_2^{\otimes2}\rs
            \ls\oR_2(\partial\vu(\vx);\vz_2)
            +\oR_1(\te{r}_2(\vx;\vz_2);\vz_1)\rs\vz_1\rp\dd\vz_2\dd\vz_1\\
    \le&
    \frac{\delta C_{\oA}}{m(\delta)^2}\int_{B_\delta}\int_{B_\delta}
        |\vK_\delta(\vz_1)||\vz_1||\vK_\delta(\vz_2)||\vz_2|
        \ls|\oR_2(\partial\vu(\vx);\vz_2)|+
            |\oR_1(\te{r}_2(\vx;\vz_2);\vz_1)|\rs\dd\vz_2\dd\vz_1.
\end{align*}
Here $\te{r}_2(\cdot;\vz_2):=\oR_2(\vu(\cdot);\vz_2)\in W^{\alpha-1}(\dom_{\delta_0/2};\re^{d\otimes3})$, for each $\vz_2\in\bll_\delta\subseteq\bll_{\delta_2}$. As before, we apply Minkowski's integral inequality to conclude
\begin{multline*}
    \|I_3\|_{L^2(\dom)}
    \le
    \frac{\delta C_{\oA}}{m(\delta)^2}\int_{B_\delta}\int_{B_\delta}
        K_\delta(|\vz_1|)K_\delta(|\vz_2|)|\vz_1|^2|\vz_2|^2\\
    \times
    \ls\|\oR_2(\partial\vu(\cdot);\vz_2)\|_{L^2(\dom;\re^{d\otimes4})}+
            \|\oR_1(\te{r}_2(\cdot;\vz_2);\vz_1)\|_{L^2(\dom;\re^{d\otimes4})}
            \rs\dd\vz_2\dd\vz_1
\end{multline*}
The convergence for $\oA_{\delta}\vu$, again, follows from Lemma~\ref{L:RCont} and assumption~\eqref{E:KernelAss}.

The arguments for $\oB_{\delta}$ are similar, but more straightforward, than those used above. For the case $0\le\alpha<1$, we use
\[
    \vu(\vx+\vz)=\vu(\vx)+\partial\vu(\vx)\vz+\frac{1}{2}\partial^2\vu(\vx)\vz^{\otimes2}
        +\oR_2(\vu(\vx);\vz)\vz^{\otimes2}.
\]
Hence, with assumption~\eqref{E:KernelAss},
\begin{align*}
    \oB_{\delta}\vu(\vx)
    =&
    \frac{C_{\oB}}{m(\delta)}\underbrace{\lp\int_{B_\delta}
    \frac{\vK_\delta(\vz)\otimes\vz^{\otimes2}}{|\vz|^2}
        \dd\vz\rp}_{=\te{0}}\partial\vu(\vx)
    +
    \underbrace{\frac{C_{\oB}}{m(\delta)}\lp\int_{B_\delta}
    \frac{\vK_\delta(\vz)\otimes\vz^{\otimes3}}{|\vz|^2}
    \dd\vz\rp}_{=\tI\otimes\tI+2\tI_{\sym}}
    \partial^2\vu(\vx)\\
    &\qqqquad+
    \frac{C_{\oB}}{m(\delta)}\int_{B_\delta}
        \lp\frac{\vK_\delta(\vz)\otimes\vz^{\otimes3}}{|\vz|^2}\rp
        \oR_2(\vu(\vx);\vz)\dd\vz.
\end{align*}
For the case $1\le\alpha\le2$, we use
\[
    \vu(\vx+\vz)=\vu(\vx)+\partial\vu(\vx)\vz+\frac{1}{2}\partial^2\vu(\vx)\vz^{\otimes2}
        +\frac{1}{6}\partial^3\vu(\vx)\vz^{\otimes3}+\oR_3(\vu(\vx);\vz)\vz^{\otimes3}.
\]
Since $\vz\mapsto\vK_\delta(\vz)\otimes\vz^{\otimes4}=K_\delta(|\vz|)\vz^{\otimes5}$ is antisymmetric, we find
\begin{align*}
    \oB_{\delta}\vu(\vx)
    =&
    \frac{C_{\oB}}{m(\delta)}\underbrace{\lp\int_{B_\delta}
    \frac{\vK_\delta(\vz)\otimes\vz^{\otimes2}}{|\vz|^2}
        \dd\vz\rp}_{=\te{0}}\partial\vu(\vx)
    +
    \underbrace{\frac{C_{\oB}}{m(\delta)}\lp\int_{B_\delta}
    \frac{\vK_\delta(\vz)\otimes\vz^{\otimes3}}{|\vz|^2}
    \dd\vz\rp}_{=\tI\otimes\tI+2\tI_{\sym}}
    \partial^2\vu(\vx)\\
    &\qqqquad+
    \frac{C_{\oB}}{m(\delta)}
        \int_{B_\delta}\lp\frac{\vK_\delta(\vz)\otimes\vz^{\otimes4}}{|\vz|^2}\rp
        \oR_3(\vu(\vx);\vz)\dd\vz.
\end{align*}
In both cases, the result again follows from Lemma~\ref{L:RCont}.
\qed
\end{proof}

\subsection{Convergence of Solutions}

Assume that $\ve{f}\in L^2(\Omega)$ and $\vu_D\in L^2(\re^d)\cap L^2(\partial\dom)$. Throughout this section, for each $\delta>0$ with  $2\delta\le\ov{\delta}$, we use $\vu_\delta\in L^2(\dom^+_{2\delta};\re^d)$ and $\vu_0\in W^{2,2}(\dom;\re^d)$ to denote the solutions to
\begin{equation}\label{eqn:static}
\lc\begin{array}{ll}
        \set{L}_\delta\vu_\delta(\vx)=\ve{f}(\vx), & \vx\in\dom,\\
        \vu_\delta(\vx)=\vu_D(\vx), & \vx\in\bnd^+_{2\delta}
    \end{array}\rper
    \nd
     \lc\begin{array}{ll}
        \set{L}_0\vu_0(\vx)=\ve{f}(\vx), & \vx\in\dom,\\
        \vu_0(\vx)=\vu_D(\vx), & \vx\in\partial\dom
    \end{array}\rper.
\end{equation}
Our main assumptions are
\begin{itemize}
\item[(A1)] With $0\le\alpha\le 2$, we find $\vu_0\in W^{2+\alpha,2}(\dom)$ and that there is an extension $\ov{\vu}_0\in W^{2+\alpha,2}(\re^d;\re^d)$ for $\vu_0$.
\item[(A2)] There exists a constant $M_1<\infty$ and $\beta\ge0$ such that,
\[
    |\ov{\vu}_0(\vx)-\vu_D(\vx)|\le M_1\dist(\vx,\partial\dom)^{1+\beta},
    \quad\text{ for a.e. }\vx\in\dom^+_{\ov{\delta}}.
\]
\item[(A3)] There exists an $M_2<\infty$ such that
\[
    \frac{\delta}{m(\delta)}\int_{B_\delta}|\vK_\delta(\vz)|
    \dd\vz\le M_2,\frl0<\delta\le\ov{\delta}.
\]
\item[(A4)] There is an $M_3<\infty$ and $\gamma\in \re$ such that
\begin{equation}\label{E:Aiv}
    \|\vu_0-\vu_\delta\|_{L^2(\bnd^-_{2\delta})}\le M_3\delta^\gamma,
    \frl0<\delta\le\ov{\delta}.
\end{equation}
\end{itemize}

\begin{theorem}\label{thm:conv}
Suppose that (A1), (A2), and (A3) hold.
\begin{enumerate}[(a)]
\item Then there exists a $C<\infty$ with the following property:
\begin{itemize}
    \item[(a1)] If $\alpha=0,1$, then for each $\vep>0$, there exists a $0<\delta=\delta(\vep)\le\ov{\delta}$ such that
\[
    \|\vu_\delta-\vu_0\|_{L^2(\dom;\re^d)}
    \le
    \vep\delta^{\alpha}+C\delta^{\beta-\frac{1}{2}}.
\]
    \item[(a2)] If $0<\alpha\le2$ and $\alpha\neq1$, then
\[
    \|\vu_\delta-\vu_0\|_{L^2(\dom;\re^d)}
    \le C\lp\delta^\alpha+\delta^{\beta-\frac{1}{2}}\rp.
\]
\end{itemize}
\item If the assumption (A4) also holds, then there exists a $C<\infty$ such that
\begin{itemize}
    \item[(b1)] If $\alpha=0,1$, then for each $\vep>0$, then there exists a $C<\infty$ with the following property:
\[
    \|\vu_\delta-\vu_0\|_{L^2(\dom;\re^d)}
    \le
    \vep\delta^{\alpha}+C\delta^{\frac{2\beta+2\gamma-1}{4}}.
\]
    \item[(b2)] If $0<\alpha\le2$ and $\alpha\neq1$, then
\[
    \|\vu_\delta-\vu_0\|_{L^2(\dom;\re^d)}
    \le C\lp\delta^\alpha+\delta^{\frac{2\beta+2\gamma-1}{4}}\rp.
\]
\end{itemize}
\end{enumerate}
\end{theorem}

\begin{proof}
Throughout the proof $C'<\infty$ denotes a constant that may change from line to line but is independent of $\delta$. Define $\ov{\vu}_\delta,\ov{\vv}_\delta\in L^2(\re^d;\re^d)$ by
\begin{equation}\label{E:FunDiff}
    \ov{\vu}_\delta(\vx)
    :=
    \lc\begin{array}{ll}
         \vu_\delta(\vx), & \vx\in\dom\\
         \ov{\vu}_0(\vx),& \vx\in\re^d\setminus\dom
    \end{array}\rper
    \nd
    \ov{\vv}_\delta(\vx)=\ov{\vu}_\delta(\vx)-\ov{\vu}_0(\vx).
\end{equation}
Since $\ov{\vv}_\delta$ is identically zero on $\re^d\setminus\dom$, we may use the Poincar\'{e}-Korn inequality in Theorem~\ref{T:PKorn}. Thus, we may select $\delta_0>0$ so that for each $0<\delta\le\delta_0$,
\begin{align}
\nonumber
    \|\vu_\delta-\vu_0\|^2_{L^2(\dom)}
    =&
    \|\ov{\vv}_\delta\|^2_{L^2(\dom)}
    \le
    C_{\PK}W_{\oB,\delta}(\ov{\vv}_\delta,\ov{\vv}_\delta)
    \le
    C_{\PK}W_\delta(\ov{\vv}_\delta,\ov{\vv}_\delta)\\
\nonumber
    =&
    -C_{\PK}\int_\dom\underbrace{\ls\set{L}_\delta\vu_\delta(\vx)\rs}_{=\ve{f}(\vx)}
        \bdot\ov{\vv}_\delta(\vx)\dd\vx
    -C_{\PK}\int_\dom\ls\set{L}_\delta\lp\ov{\vu}_\delta-\vu_\delta\rp(\vx)\rs
            \bdot\ov{\vv}_\delta(\vx)\dd\vx\\
\nonumber
    &\qqqquad
    +C_{\PK}\int_\dom\ls\set{L}_\delta\ov{\vu}_0(\vx)\rs\bdot\ov{\vv}_\delta(\vx)\dd\vx\\
\label{E:MainConvIneq}
    =&
    C_{\PK}\underbrace{\int_\dom\ls\set{L}_\delta\ov{\vu}_0(\vx)-\set{L}_0\ov{\vu}_0(\vx)\rs
        \bdot\ov{\vv}_\delta(\vx)\dd\vx}_{=:I_{1,\delta}}
    +C_{\PK}\underbrace{\int_\dom\ls\set{L}_\delta\lp\vu_\delta-\ov{\vu}_\delta\rp(\vx)\rs
            \bdot\ov{\vv}_\delta(\vx)\dd\vx}_{=:I_{2,\delta}}.
\end{align}
With assumption (A1), we may apply Theorem~\ref{T:LConv} to bound $|I_{1,\delta}|$ in terms of $\delta$.

We use the remaining assumptions to bound $|I_{2,\delta}|$. Define $\vv_\delta\in L^2(\re^d;\re^d)$ by
\[
    \vv_\delta(\vx):=\vu_\delta(\vx)-\ov{\vu}_\delta(\vx)
    =\lc\begin{array}{ll}
        \ve{0}, & \vx\in\dom,\\
        \vu_D(\vx)-\ov{\vu}_0(\vx), & \vx\in\re^d\setminus\dom.
    \end{array}\rper
\]
Since the support of $\vv_\delta$ is contained in $\re^d\setminus\dom$, we find that $\oL_\delta(\vv_\delta)=\ve{0}$ on $\dom^-_{2\delta}$. Thus
\begin{equation}\label{E:I2Equal}
    I_{2,\delta}
    =
    \int_{\bnd^-_{2\delta}}\ls\oL_\delta\vv_\delta(\vx)\rs\bdot\ov{\vv}_\delta(\vx)\dd\vx.
\end{equation}
We now work to bound $\|\oL_\delta\vv_\delta\|_{\bnd^-_{2\delta}}$. Assume (A2) and (A3) hold. From its definition,
\begin{equation}\label{E:LIneq1}
    \|\oL_\delta\vv_\delta\|_{\bnd^-_{2\delta}}
    \le
    \|\oA_\delta\vv_\delta\|_{\bnd^-_{2\delta}}+\|\oB_\delta\vv_\delta\|_{\bnd^-_{2\delta}}.
\end{equation}

For a.e. $\vx\in\bnd^-_\delta$ and a.e. $\vz_1,\vz_2\in B_\delta$, we find
\begin{align*}
    |\vv_\delta(\vx+\vz_1+\vz+2)|
    \le&
    \dist(\vx+\vz_1+\vz_2,\partial\dom)^{1+\beta}
    \le
    \lp \dist(\partial\dom,\vx)+|\vz_1|+|\vz_2|\rp^{1+\beta}\\
    \le&
    3\delta^{1+\beta}.
\end{align*}
Thus, for a.e. $\vx\in\bnd^-$
\begin{align*}
    |\oA_\delta\vv_\delta(\vx)|
    \le&
    \frac{C_{\oA}}{m(\delta)}\int_{B_\delta}K_\delta(|\vz_1|)|\vz_1|
        |\theta_\delta(\vx+\vz_1);\vv_\delta)|\dd\vz_1\\
    \le&
    \frac{C'}{m(\delta)^2}
    \int_{B_\delta}\int_{B_\delta}K_\delta(|\vz_1|)|\vz_1|
        K_\delta(|\vz_2|)|\vz_2|\vv_\delta(\vx+\vz_1+\vz_2)\dd\vz_2\dd\vz_1\\
    \le&
    \frac{C'\delta^{1+\beta}}{m(\delta)^2}\int_{B_\delta}\int_{B_\delta}
        K_\delta(|\vz_1|)|\vz_1|
        K_\delta(|\vz_2|)|\vz_2|\dd\vz_2\dd\vz_1\\
    \le&
    C'\delta^{\beta-1}.
\end{align*}
It follows that
\[
    \|\oA_\delta\vv_\delta\|_{\bnd^-_{2\delta}}
    \le
    C'\delta^{\beta-1}|\bnd^-_{2\delta}|^\frac{1}{2}\le C'\delta^{\beta-\frac{1}{2}}.
\]
We similarly find $\|\oB_\delta\vv_\delta\|^2_{\bnd^-_{2\delta}}\le C'\delta^{\beta-\frac{1}{2}}$. From these bounds,~\eqref{E:LIneq1}, and the bounds for $I_{1,\delta}$ in Theorem~\ref{T:LConv}, we obtain
\[
    \|\ov{\vv}_\delta\|^2_{L^2(\dom)}
    \le
    C'
    \lp\|\oL_\delta\ov{\vu}_0-\oL_0\ov{\vu}_0\|_{L^2(\dom)}
        \|\ov{\vv}_\delta\|_{L^2(\dom)}
    +\|\oL_\delta\vv_\delta\|_{L^2(\bnd^-_{2\delta})}
        \|\ov{\vv}_\delta\|_{L^2(\bnd^-_{2\delta})}\rp,
\]
and part (a) of the Theorem follows. For part (b), we also assume (A4). From the inequality above, we get
\[
    \|\ov{\vv}_\delta\|^2_{L^2(\dom)}
    \le
    C'
    \lp\|\oL_\delta\ov{\vu}_0-\oL_0\ov{\vu}_0\|_{L^2(\dom)}
        \|\ov{\vv}_\delta\|_{L^2(\dom)}
    +\delta^{\beta+\gamma-\frac{1}{2}}\rp.
\]
This proves the result, else we reach a contradiction.
\qed
\end{proof}
\begin{remark}\label{rmk1}
The results of Theorem \ref{thm:conv} above provide three ``knobs" for identifying lower bounds for rates of convergence, as follows: 
\begin{itemize}
\item[$\bullet$] The exponent $\alpha$ is limited by the degree of smoothness for the extended local solution and is independent of the collar values for the nonlocal problem. For a local solution with a smooth extension $\alpha=2$, which corresponds to quadratic convergence away from the boundary.
\item[$\bullet$] The exponent $\beta$ quantifies the order at which the prescribed collar values converge to the extension of the local solution as the boundary $\partial\dom$ is approached. For constant extensions of the boundary values on $\partial\dom$, we have $\beta\ge0$; for linear extensions, $\beta\ge1$. If the prescribed collar data is provided by the extension $\ov{\vu}_0$ of the local solution, then $\beta$ is effectively infinite and the convergence is of order $\delta^\alpha$.
\item[$\bullet$] Part (a) of the theorem only assumes $\{\vu_\delta\}_{\delta>0}\subseteq L^2(\dom)$. The family of nonlocal solutions need not even be uniformly bounded in $L^2(\dom)$. If information about the rate of convergence of $\vu_\delta\to\vu_0$ on the interior collar $\bnd^-_{2\delta}$ is available, this can be captured with the parameter $\gamma$. For example, if there is a uniform bound for the nonlocal solutions, then $\gamma=\frac{1}{2}$.
\end{itemize}
\end{remark}

We next provide a corollary, which provides an alternative to (A4). Given $\vv\in L^2(\re^d;\re^d)$, define
\[
    \oG\vv(\vx)
    :=
    \frac{1}{m(\delta)}\int_{B_\delta}|\vK_\delta(\vz)||\vv(\vx+\vz)-\vv(\vx)|\dd\vz.
\]
Consider the following assumption, where $\vv_\delta$ is defined in~\eqref{E:FunDiff}.
\begin{itemize}
\item[(A4')] There exists an $M_3'<\infty$ and a $\gamma'\ge 0$ such that
\[
    \|\oG\ov{\vv}_\delta\|_{L^\infty(\bnd_{2\delta}^-)}\le M_3'\delta^{\gamma'},
    \frl0<\delta\le\ov{\delta}.
\]
\end{itemize}

\begin{corollary}\label{cor:conv}
Suppose that $\dom$ satisfies the uniform exterior sphere condition. In addition to (A1), (A2), and (A3), assume (A4'), then part (b) of Theorem~\ref{thm:conv} holds with $\gamma=\gamma'+\frac{3}{2}$.
\end{corollary}
\begin{proof}
We use an argument similar to the one used in~\cite{foss2019nonlocal}. Given a set $E\subseteq\re^d$ and $\delta>0$, define $E_\delta:=\{\rho\vx\in\re^d:0\le\rho<\delta\text{ and }\vx\in E\}$. The uniform exterior sphere condition ensures there is a $\ov{\delta}>0$, $j_0\in\nats$, and open sets $\{U_k\}_{k=1}^{k_0}\subseteq\bnd^-_{\ov{\delta}}$ and $\{S_k\}_{k=1}^{k_0}\subseteq \partial B_1$ with the following properties:
\begin{itemize}
    \item[$\bullet$] $\bnd^+_{2\ov{\delta}}\subseteq\bigcup_{k=1}^{k_0}U_k$ and there exists a $\sigma>0$ and $0<\rho_0<\ov{\delta}$ such that
    \[
        \frac{1}{m(\delta)}\int_{Z_{\rho,k}} K_\delta(|\vz|)|\vz|^2\dd\vz>\sigma,
        \frl 0<\rho\le\ov{\delta},
    \]
    where $Z_{\rho,k}:=\rho Z_k\setminus B_{\rho_0}$.
    \item[$\bullet$] For each $\delta>0$, $\vx\in U_{\delta,k}:=\{\vy\in U_k:\dist(\vy,\partial\dom)<\delta\}$ and $\vom\in Z_k$, we find $\vx+\rho\vom\in\bnd^-_{2\delta}\cup\re^d\setminus\dom$ for all $0<\rho\le\delta$.
    \item[$\bullet$] For each $0<\delta\le\ov{\delta}$ and $\vx\in U_{\delta,k}$ there exists a $1\le j_{\vx}\le j_0$ such that $\vx+j\vz\in\re^d\setminus\dom$ for all $\vz\in Z_{\delta,k}$.
\end{itemize}
Roughly speaking, we require a finite collection of subsets of the unit sphere such that for each $\vx\in U_{\delta,k}$ and each $\vz$ in the associated annular sector $Z_{\delta,k}$ of $B_\delta$, the sequence $\{\vx+j\vz\}_{j=0}^{j_{\vx}}$ stays within $2\delta$ units of the $\partial\dom$ and terminates at $\vx+j_{\vx}\vz\in\re^d\setminus \dom$. An example is a square, in which case we can use four trapezoidal regions $\{U_k\}_{k=1}^4$ that border the boundary of the square and four corresponding quarter sectors $\{S_k\}_{k=1}^4$ of the unit circle.

Let $k=1,\dots,k_0$ and $\delta>0$ be given. Since $\ov{\vv}_\delta=\ve{0}$ on $\re^d\setminus\dom$, for each $\vx\in U_{\delta,k}$ and $\vz\in Z_{\delta,k}$, we have
\[
    |\ov{\vv}_\delta(\vx)|\le\sum_{j=1}^{j_{\vx}}
        |\ov{\vv}_\delta(\vx+j\vz)-\ov{\vv}(\vx+(j-1)\vz)|,
\]
with $\vx_j=\vx+j\vz$. Multiplying both sides by $|\vK_\delta(\vz)||\vz|m(\delta)^{-1}$ and integrating with respect to $\vz$ yields
\[
    |\ov{\vv}_\delta(\vx)|
    \le
    \frac{1}{\sigma m(\delta)}\sum_{j=0}^{j_{\vx}-1}\int_{Z_k}|\vK_\delta(\vz)||\vz|
        |\ov{\vv}_\delta(\vx_{j}+\vz)-\ov{\vv}(\vx_{j})|\dd\vz
    \le
    \frac{M_3'j_0}{\sigma}\delta^{1+\gamma'}.
\]
Taking the $L^2$-norm over $\bnd^-_{2\delta}$, we see that assumption (A4) is satisfied with $\gamma=\gamma'+\frac{3}{2}$. Thus, we may apply Theorem 5.
\qed
\end{proof}
\begin{remark}
\begin{itemize}
    \item[$\bullet$] The uniform exterior sphere condition is satisfied, for example, by any $C^2$-domain or a convex set.
    \item[$\bullet$] A straightforward modification of the argument can be used to replace $\|\oG\ov{\vv}_\delta\|_{L^\infty(\bnd^-_{2\delta})}$ with $\|\oG\ov{\vv}_\delta\|_{L^2(\bnd^-_{2\delta})}$ in assumption (A4'). In this case, we obtain $\gamma=\gamma'+1$
\end{itemize}
\end{remark}

\section{Convergence of Solutions for the Dynamic System}\label{sec:dynamic}

The results obtained in the steady-state case problem extend to the dynamic system, with similar bounds, but constants that may grow exponentially in time. Suppose that $\vu_\delta\in C^2(0, T;L^2(\dom^+_{\delta}))$ and $\vu_0\in C^2(0,T; W^{1,2}(\dom))$ satisfy the initial boundary value problems
\begin{equation}\label{eqn:dynamic}
    \lc\begin{array}{ll}
        \rho\vu_{{\delta}_{tt}}-\set{L}_\delta\vu_\delta=\ve{f}(t,\vx),\, & \vx\in\dom, \, t>0\\
        \vu_\delta(t, \vx)=\vu_D(t,\vx), & \vx\in\bnd^+_{2\delta}, \, t>0\\
        \vu_{\delta}(0,\vx)=\phib(\vx), & \vx\in \dom^+_{2\delta}\\
        \vu_{{\delta}_t}(0,\vx)=\psib(\vx), & \vx\in \dom^+_{2\delta}
    \end{array}\rper
    \nd
     \lc\begin{array}{ll}
        \rho\vu_{0_{tt}}-\set{L}_0\vu_0=\ve{f}(t,\vx), & \vx\in\dom, t>0\\
        \vu_0(t,\vx)=\vu_D(t, \vx), & \vx\in\partial\dom, \, t>0\\
        \vu_{0}(0,\vx)=\phib(\vx), & \vx\in \dom^+_{2\delta}\\
        \vu_{0_t}(0,\vx)=\psib(\vx), & \vx\in \dom^+_{2\delta}.
    \end{array}\rper
\end{equation}
Here we consider the material density $\rho=1$ without loss of generality, and denote
\[
    \vu_{\delta_t}=\pder{}{t}\vu_{\delta}
    \nd
    \vu_{\delta_{tt}}=\pder{{}^2}{t^2}\vu_{\delta}.
\]

The precise statements mirror Theorem \ref{thm:conv} and are given by
\begin{theorem}\label{thm:dynamicconv}
Let $\ve{f}\in L^2(0, T; L^2(\Omega)), \,  \vu_D\in C(0,T; L^2(\Omega_{2\delta}^+)\cap L^2(\partial\Omega))$ and $\phib, \psib \in L^2(\dom^+_{2\delta})$ be given. Let $\alpha,\beta,\gamma\ge0$ be given, and suppose that
$\vu_0\in C^2(0,T; W^{2+\alpha,2}(\dom;\re^d))$. 
Then the estimates in parts (a1), (a2) of Theorem 5 hold under same assumptions (A1)-(A3) for the time dependent solution, with all statements appropriately adjusted for the dynamic case, and the bounds $M_1, M_2, C$ time dependent. 

For the second part, we additionally assume  that 
\begin{itemize}
\item[(A4'')] There is an $M_3'<\infty$ and $\gamma_0\in \re$ such that
\begin{equation}\label{E:Aiv1}
    \|\vu_{0_t}-\vu_{\delta_t}\|_{L^2(\bnd^-_{2\delta})}\le M_3'\delta^{\gamma_0},
    \frl0<\delta\le\ov{\delta}.
\end{equation}
\end{itemize}
Then 
then there exists a $C=C(t)<\infty$ such that
\begin{itemize}
    \item[(b1')] If $\alpha=0,1$, then for each $\vep>0$:
\[
    \|\vu_\delta(t)-\vu_0(t)\|_{L^2(\dom;\re^d)}
    \le
    \vep\delta^{\alpha}+C\delta^{\frac{2\beta+2\gamma_{_0}-1}{4}}.\]
    \item[(b2')] If $0<\alpha\le2$ and $\alpha\neq1$, then
\[
    \|\vu_\delta(t)-\vu_0(t)\|_{L^2(\dom;\re^d)}
    \le C\lp\delta^\alpha+\delta^{\frac{2\beta+2\gamma_{_0}-1}{4}}\rp.
\]
\end{itemize}

\end{theorem}
\begin{proof}
We use the same approach as in the proof of Theorem \ref{thm:conv}, with the same notation and extensions to the collar $\Gamma_{2\delta}^+$). Subtracting the equations from \eqref{eqn:dynamic} that hold on $\Omega$, multiplying by the difference $\ov{\vv}_{\delta}:=\ov{\vu}_{\delta}-\ov{\vu}_0$, integrating over the domain, and using the symmetry of the kernel (nonlocal integration by parts) yield:
\begin{equation}\label{eq:diffineq}
\frac{d}{2dt}\int_{\Omega} (\ov{\vv}_{{\delta}_t})^2\, \dd\vx 
+ W_\delta(\ov{\vv}_{\delta},\ov{\vv}_{\delta_t}) 
=  \int_\dom\ls\set{L}_\delta\ov{\vu}_0(\vx)-\set{L}_0\ov{\vu}_0(\vx)\rs
        \bdot\ov{\vv}_{\delta_t}(\vx)\dd\vx
    +\int_\dom\ls\set{L}_\delta\lp\vu_\delta-\ov{\vu}_\delta\rp(\vx)\rs
            \bdot\ov{\vv}_{\delta_t}(\vx)\dd\vx
\end{equation}

Using Young's inequality, bounds on $\bnd^+_{2\delta}$ and rearranging the terms gives the estimate:
\begin{equation}\label{ineq1}
\frac{d}{dt}\left( \int_{\Omega} (\ov{\vv}_{{\delta}_t})^2\, \dd\vx + W_\delta(\ov{\vv}_{\delta},\ov{\vv}_{\delta})\right) 
\leq 8\|\ov{\vv}_{\delta_t}\|_{L^2(\Omega)}^2+ 4\|\set{L}_{\delta} \ov{\vu}_0-\oL_0\ov{\vu}_0 \|_{L^2(\Omega)}^2+4\|\set{L}_{\delta} (\vu_\delta-\ov{\vu}_\delta) \|_{L^2(\Gamma^-_{2\delta})}^2,     \end{equation}
where we used the fact that $\vu_\delta-\ov{\vu}_\delta=0$ in $\Omega$.
Let 
\[
\zeta(t):=\int_{\Omega} (\vv_{{\delta}_t})^2\, \dd\vx + W_\delta(\ov{\vv}_{\delta},\ov{\vv}_{\delta})
\]
then \eqref{ineq1} together with the bounds 
\begin{equation}\label{eq:boundsL}
\|\set{L}_{\delta} \ov{\vu}_0-\oL_0\ov{\vu}_0 \|_{L^2(\Omega)}^2 \leq C_1\delta^{2\alpha},\quad\quad 
\|\set{L}_{\delta} (\vu_\delta-\ov{\vu}_\delta) \|_{L^2(\Gamma^-_{2\delta})}^2\leq C_2\delta^{2\beta-1}
\end{equation}
which are established exactly in the same way as in the proof of Theorem \ref{thm:conv}, we obtain
\[
\zeta'(t)\leq 8\zeta(t)+ C_1\delta^{2\alpha}+C_2\delta^{2\beta-1}.
\]
We obtain the bound
\begin{equation}\label{eq:zeta}
\zeta(t)\leq C(e^t-1)(\delta^{2\alpha}+\delta^{2\beta-1} ),
\end{equation}
so $W_\delta(\ov{\vv}_{\delta},\ov{\vv}_{\delta}) \leq C(e^t-1)(\delta^{\alpha}+\delta^{2\beta-1} )$. By the Poincar\'e-Korn inequality we obtain that 
\[
\|\ov{\vv}_{\delta}\|_{L^2(\Omega)} \leq C(e^t-1)(\delta^{\alpha}+\delta^{\beta-\frac{1}{2}} ),
\]
where $C$ depends also on the Poincar\'e-Korn constant. 

For the second part, we go back to inequality \eqref{ineq1} and use the bounds \eqref{eq:boundsL} and assumption (A4'') to obtain:
\begin{equation}\label{eq:diffineq2}
\zeta'(t)\leq 8\zeta+C\delta^{2\alpha}+M_3'\delta^{\beta-\frac{1}{2}+\gamma_0}.
\end{equation}
Appying Gronwall's inequality again yields the bounds (b1') and (b2').
\qed
\end{proof}


\section{Numerical Tests and Empirical Convergence Rates}\label{sec:numerics}

 \begin{figure}[!htb]\centering
  \subfigure{\includegraphics[width=0.35\textwidth]{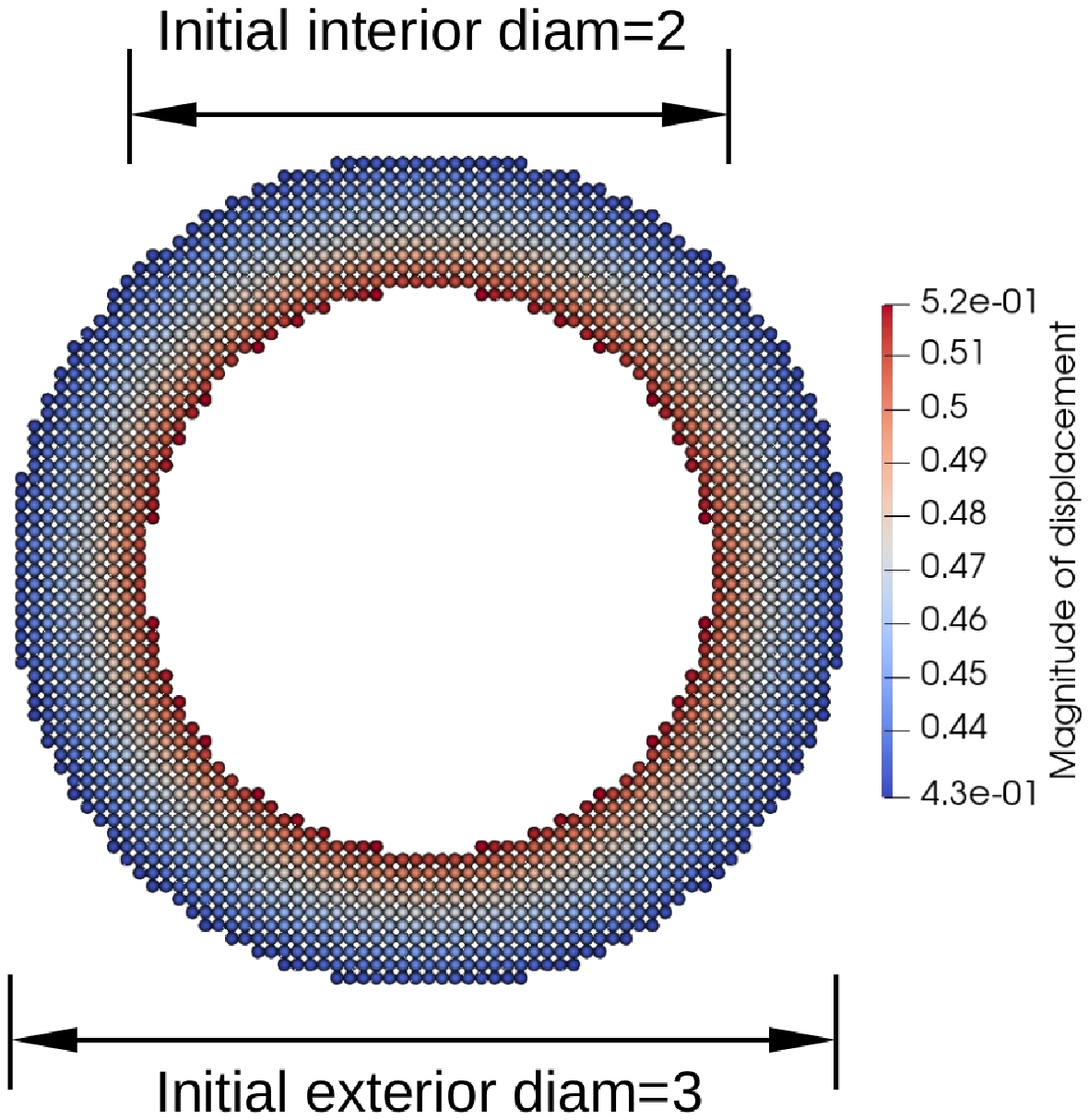}}
 \subfigure{\includegraphics[width=0.3\textwidth]{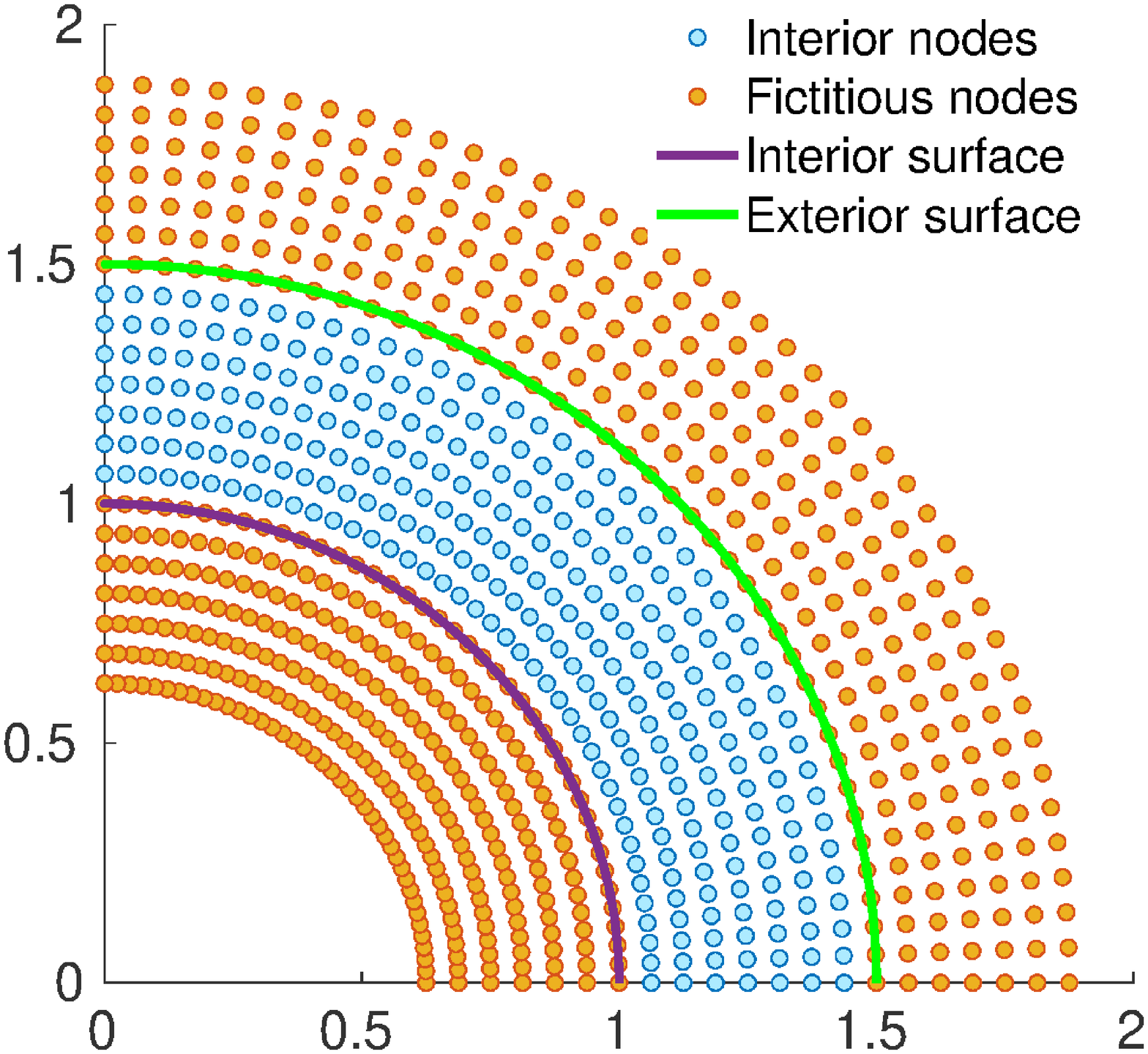}}
  \subfigure{\includegraphics[width=0.3\textwidth]{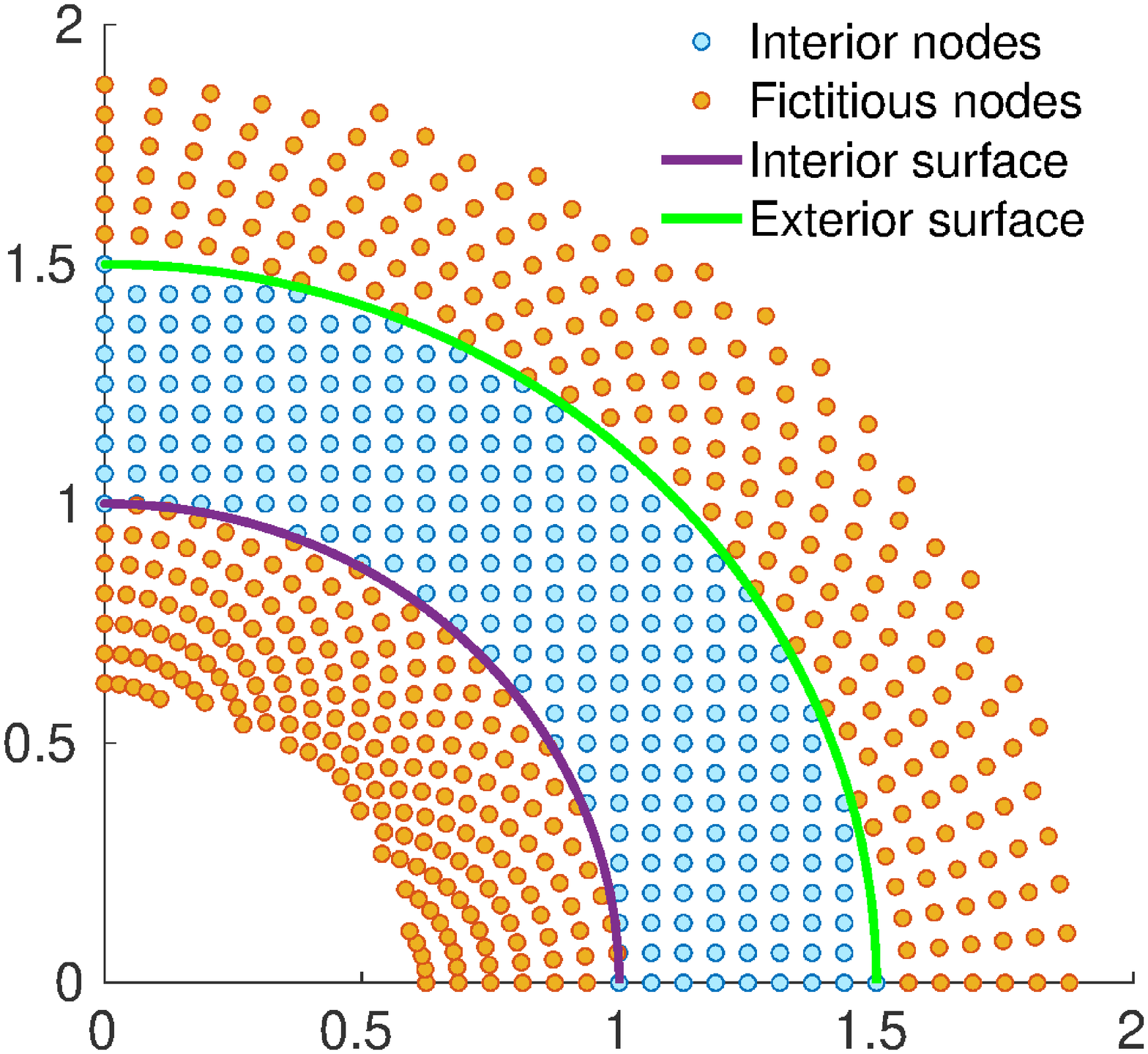}} 
 \caption{Problem settings for a hollow disk and two types of grids (in blue) with mirror-based fictitious nodes (in orange) for the hollow cylinder example. Left: polar grids where the interior nodes are generated as $X_h\cap\Omega:=\{(p_1h\cos(\pi p_2h/5),p_1h\sin(\pi p_2h/5))|\bm{p}=(p_1,p_2)\in\mathbb{Z}^2\}\cap\Omega$. Right: Cartesian grids where the interior nodes are generated as $X_h\cap\Omega:=\{(p_1h,p_2h|\bm{p}=(p_1,p_2)\in\mathbb{Z}^2\}\cap\Omega$.
 }
 \label{fig:diskgrids}
\end{figure}

In this section, we aim to numerically verify the analysis developed in Sections \ref{sec:conv}-\ref{sec:dynamic}. In particular, with nonlocal Dirichlet-type boundary conditions define by extending the surface (local) data into volumetric data, we denote the solution of \eqref{eqn:probd} as $\ub_\delta$ and its local limit, the solution of \eqref{eqn:local}, as $\ub_0$. We now investigate the convergence of numerical approximations $\ub_\delta$, by taking $\delta\rightarrow 0$ and the spatial refinement $h\rightarrow 0$ simultaneously. Three types of boundary conditions will be investigated: the ``smooth extension'' where a smooth extended local solution is provided as $\ub_D(\xb)$ for $\xb\in\Gamma^+_{2\delta}$, the ``constant extension'' which employs the naive fictitious node methods such that $\ub_D(\xb)$ is defined by the surface (local) data at its corresponding projection on $\partial\Omega$, and the ``linear extension'' which employs the mirror-based fictitious node methods and extends the surface (local) and the interior data into the exterior layer $\Gamma^+_{2\delta}$ linearly. As an illustration for the mirror grids, two example meshes for a hollow cylinder problem are provided in Figure \ref{fig:diskgrids}, with further details to be provided later on in the section.

To observe the convergence rates of $\|\vu_\delta-\vu_0\|_{L^2(\dom;\re^d)}$, an important feature of the discretization would be to preserve this asymptotic limit as $\delta,h\rightarrow 0$. Discretizations which preserve the correct local limit under spatial refinement $h \rightarrow 0$ and $\delta\rightarrow 0$ are termed asymptotically compatible (AC) \cite{tian2014asymptotically}. For further discussions and an incomplete list of AC methods see  \cite{d2020numerical,du2016local,hillman2020generalized,leng2019asymptotically,pasetto2018reproducing,seleson2016convergence,tao2017nonlocal,trask2019asymptotically,You_2019,you2020asymptotically}.  Here we numerically solve \eqref{eqn:probd} by employing a meshfree, particle discretization method introduced in \cite{yu2021asymptotically} and analyzed in \cite{fan2021asymptotically}. This optimization-Based meshfree method features asymptotic compatibility in the $\delta-$convergence tests \cite{bobaru2009convergence}, i.e., when $\ub_\delta\overset{\delta\rightarrow 0}{\rightarrow}\ub_0$ and one refines both $\delta$ and $h$ at the same rate, the numerical nonlocal solution converges to the local limit. In this work we will also inverstigate the convergence of numerical nonlocal solutions under the $\delta-$convergence setting, since banded stiffness matrices are obtained in such a setting and scalable implementations are allowed. Therefore, we always choose $h$ such that the ratio $\frac{h}{\delta}$ is bound by a constant $M$ as $\delta \rightarrow 0$.

Discretizing the whole interaction region $\Omega^+_{2\delta}$ by a collection of points $X_{h} = \{\vx_i\}_{\{i=1,2,\cdots,N_p\}} \subset \Omega^+_{2\delta}$, we aim to solve for the displacement $\ub_i\approx\ub(\xb_i)$ and the nonlocal dilatation $\theta_i\approx \theta(\xb_i)$ on all $\xb_i\in X_h\cap \omg$. We first characterize the distribution of collocation points as follows. Recall the definitions \cite{wendland2004scattered} of fill distance $h_{X_h,\Omega} = \underset{\xb_i \in X_h}{\sup}\, \underset{\xb_i \in X_h}{\min}||\xb_i - \xb_j||_2$ and separation distance $q_{X_h} = \frac12 \underset{i \neq j}{\min} ||\xb_i - \xb_j||_2$, 
we assume that $X_h$ is quasi-uniform, namely that there exists $c_{qu} > 0$ such that $q_{\chi_h} \leq h_{\chi_h,\Omega} \leq c_{qu} q_{\chi_h}$. 

We first consider the spatial discretization for the static LPS model \eqref{eqn:probd} through the following one point quadrature rule at $X_h$ \cite{silling_2010}:
\begin{equation}\label{eq:discreteNonlocElasticity}
\left\{\begin{array}{ll}
    (\mathcal{L}_\delta^h \ub)_i:=-\frac{C_{\oA}}{m(\delta)} \sum_{\xb_j \in B_\delta (\vx_i)} \left(\lambda - \mu\right) K_\delta(\verti{\xb_j-\xb_i}) \left(\vx_j-\vx_i\right)\left(\theta_i + \theta_j \right) \omega_{j,i}&\\
  \qquad-  \frac{C_\beta}{m(\delta)}\sum_{\xb_j \in B_\delta (\vx_i)} \mu K_\delta(\verti{\xb_j-\xb_i})\frac{\left(\vx_j-\vx_i\right)\otimes\left(\vx_j-\vx_i\right)}{\left|\vx_j-\vx_i\right|^2} \cdot \left(\vu_j- \vu_i \right) \omega_{j,i} = \vf_i,&\, \text{for }\xb_i \text{ in }\Omega,\\
\theta_i =  \frac{d}{m(\delta)} \sum_{\xb_j \in B_\delta (\vx_i)}  K_\delta(\verti{\xb_j-\xb_i}) \left(\vx_j-\vx_i\right) \cdot  \left(\vu_j - \vu_i \right) \omega_{j,i},& \, \text{for }\xb_i \text{ in }\omg^+_{\delta},\\
\ub_i=\ub_D(\xb_i), & \, \text{for }\xb_i \text{ in }\Gamma^+_{2\delta},\\
  \end{array}\right.
\end{equation}  
where we specify $\left\{\omega_{j,i}\right\}$ as a to-be-determined collection of quadrature weights admitting interpretation as a measure associated with each collocation point $\xb_i$. Particularly, we employ an optimization-based approach to define these weights \cite{trask2019asymptotically}, by seeking $\omega_{j,i}$ for integrals supported on balls of the form
\begin{equation}
I[f] := \int_{B_\delta (\vx_i)} f(\xb,\yb) d\yb \approx I_h[f] := \sum_{\xb_j \in B_\delta (\vx_i)} f(\xb_i,\xb_j) \omega_{j,i}
\end{equation}
where we include the subscript $i$ in $\left\{\omega_{j,i}\right\}$ to denote that we seek a different family of quadrature weights for different subdomains $B_\delta(\vx_i)$. We obtain these weights from the following optimization problem
\begin{align}\label{eq:quadQP}
  \underset{\left\{\omega_{j,i}\right\}}{\text{argmin}} \sum_{\xb_j \in B_\delta (\vx_i)} \omega_{j,i}^2 \quad
  \text{such that}, \quad
  I_h[p] = I[p] \quad \forall p \in \mathbf{V}_h
\end{align}
where $\mathbf{V}_h$ denotes a Banach space of functions which should be integrated exactly. For the LPS model we take $\mathbf{V}_h:=\left\{ q = \frac{p(\yb)}{|\yb-\xb|^3} \,|\, p \in P_5(\mathbb{R}^d) \text{ such that } \int_{B_\delta(\xb)} q(\yb) d \yb < \infty \right\}$ where $p \in P_5(\mathbb{R}^d)$ is the space of quintic polynomials. As shown in \cite{yu2021asymptotically}, the above $V_h$ provides a reproducing space which is sufficient to integrate \eqref{eq:discreteNonlocElasticity} exactly in the case where $\ub$ and $\theta$ are quadratic polynomials. We refer to previous work \cite{trask2019asymptotically,you2020asymptotically,yu2021asymptotically} for further information, analysis, and implementation details.

For the dynamic LPS model \eqref{eqn:dynamic}, to discretize in time we apply the Newmark scheme. With time step size $\Delta t$, at the $(n+1)-$th time step we solve for the displacement $\ub_i^{n+1}\approx\ub((n+1)\Delta t,\xb_i)$ and the nonlocal dilatation $\theta_i^{n+1}\approx \theta((n+1)\Delta t,\xb_i)$ following:
\begin{equation}\label{eqn:probdis}
\left\{\begin{array}{ll}
    \frac{4\rho}{\Delta t^2} \ddot{\ub}_i^{n+1}+(\mcL^h_\delta
    \ub)_i^{n+1}= \fb((n+1)\Delta t,\xb_i)
    +\frac{4\rho}{\Delta t^2}(\ub_{i}^n + \Delta t \dot{\ub}_i^n + \frac{\Delta t^2}{4}
    \ddot{\ub}_i^n)
    ,& \quad \text{for }\xb_i \text{ in }\Omega,\\
\theta_i^{n+1}=\dfrac{d}{m(\delta)} \sum\limits_{\xb_j \in B_{\delta}(\xb_i)} K_{ij} (\xb_j-\xb_i)^T \left(\ub_j^{n+1} - \ub_i^{n+1} \right)\omega_{j,i},&\quad \text{for }\xb_i \text{ in }\omg^+_\delta,\\
\ub_i^{n+1}=\ub_D(\xb_i), & \quad \text{for }\xb_i \text{ in }\Gamma^+_{2\delta},\\
\ub_i^0=\phib(\xb_i),\;\dot{\ub}_i^0=\psib(\xb_i), &\quad \text{ for }\xb_i \in\omg^+_{2\delta},\\
\end{array}\right.
\end{equation}
where $\mcL^h_{\delta}$ is the discretized nonlocal operators as defined in \eqref{eq:discreteNonlocElasticity}. The acceleration and velocity at the $n+1$-th time step are then calculated as follows:
$$\ddot{\ub}_i^{n+1} 
:= \frac{4}{\Delta t^2}(\ub_i^{n+1}
-\ub_i^{n}-\Delta t \ub_i^n) - \ddot{\ub}_i^n, \quad
\dot{\ub}_i^{n+1} := \dot{\ub}_i^n
+ \frac{\Delta t}{2} (\ddot{\ub}_i^n
+ \ddot{\ub}_i^{n+1}).$$
Note that although the Newmark scheme is unconditionally stable, in all numerical tests we take a sufficiently small time step size $\Delta t\ll \delta$ so as to study the convergence rates with respect to $\delta$.

In the following we consider static and dynamic problems on three numerical cases: linear patch tests, smooth manufactured nonliear solutions, and analytical solutions to curvilinear surface problems. For the first two cases we consider square domains with Cartesian grids and mesh spacing $h$. For the analytical solutions to curvilinear surface problems, we consider the deformation of a hollow cylinder under an internal pressure $p_0$, with the setting plotted in Figure \ref{fig:diskgrids}. Two types of grids are generated and tested: polar grids where the interior nodes are generated as $X_h\cap\Omega:=\{(p_1h\cos(\pi p_2h/5),p_1h\sin(\pi p_2h/5))|\bm{p}=(p_1,p_2)\in\mathbb{Z}^2\}\cap\Omega$, and Cartesian grids where the interior nodes are generated as $X_h\cap\Omega:=\{(p_1h,p_2h|\bm{p}=(p_1,p_2)\in\mathbb{Z}^2\}\cap\Omega$. In all cases we adopt material parameters under plane strain assumptions:
$$E=1,\lambda=E\nu/((1+\nu)(1-2\nu)),\mu=E/(2(1+\nu)).$$
Two values of Poisson ratio $\nu=0.3$ and $0.49$ are investigated which correspond to compressible and nearly-incompressible materials, respectively.

We summarize the setup and report the formal convergence study for the static LPS problem in Section \ref{sec:staticnonlocal} and for the dynamic problem in Section \ref{sec:dynamicnonlocal}. To investigate whether the theoretical convergence rates in Theorem \ref{thm:conv} and Theorem \ref{thm:dynamicconv} are realized as $\delta \rightarrow 0$, our numerical results particularly focus on identifying the convergence rates of $\vertii{\ub_\delta-\ub_0}_{L^2(\Omega;\re^d)}$ and numerically evaluating the value of $\gamma'$ by calculating $\|\oG\ov{\vv}_\delta\|_{L^\infty(\bnd_{2\delta}^-)}=O(\delta^{\gamma'})$.

\subsection{Static Problem with Dirichlet-type Boundary Conditions}\label{sec:staticnonlocal}

 \begin{figure}[!htb]\centering
 \subfigure{\includegraphics[width=0.49\textwidth]{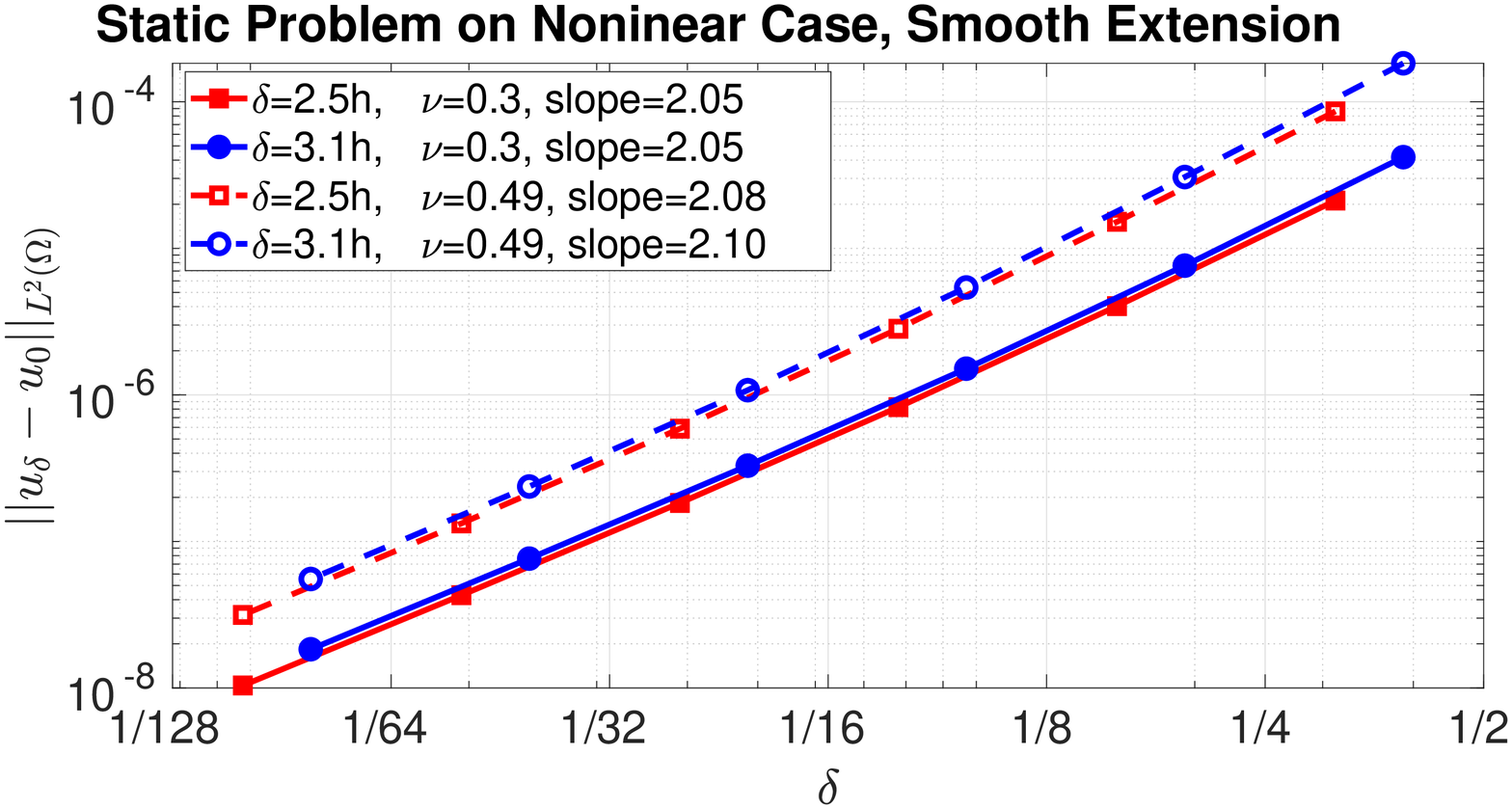}}
 \subfigure{\includegraphics[width=0.49\textwidth]{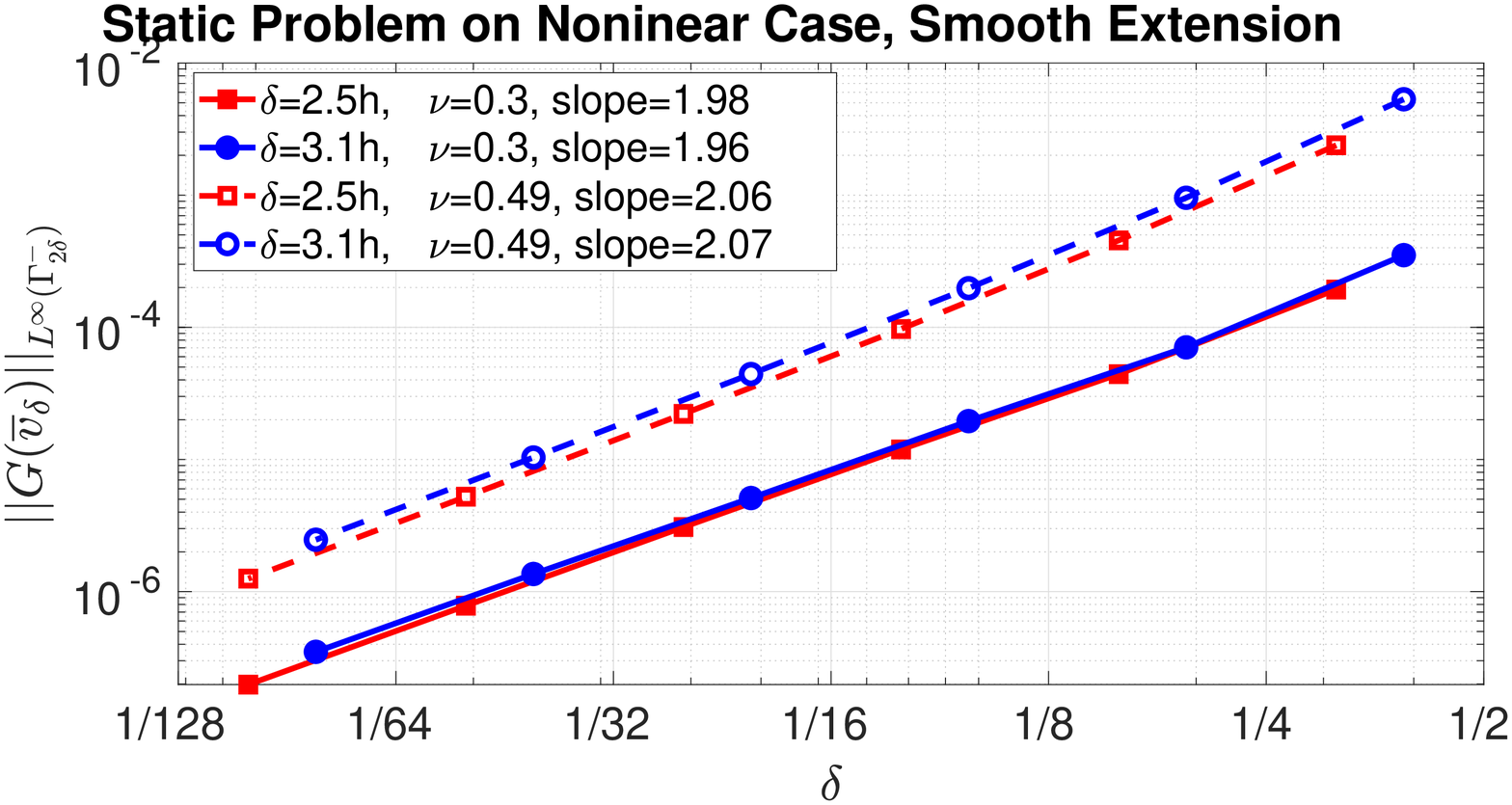}}
\subfigure{\includegraphics[width=0.49\textwidth]{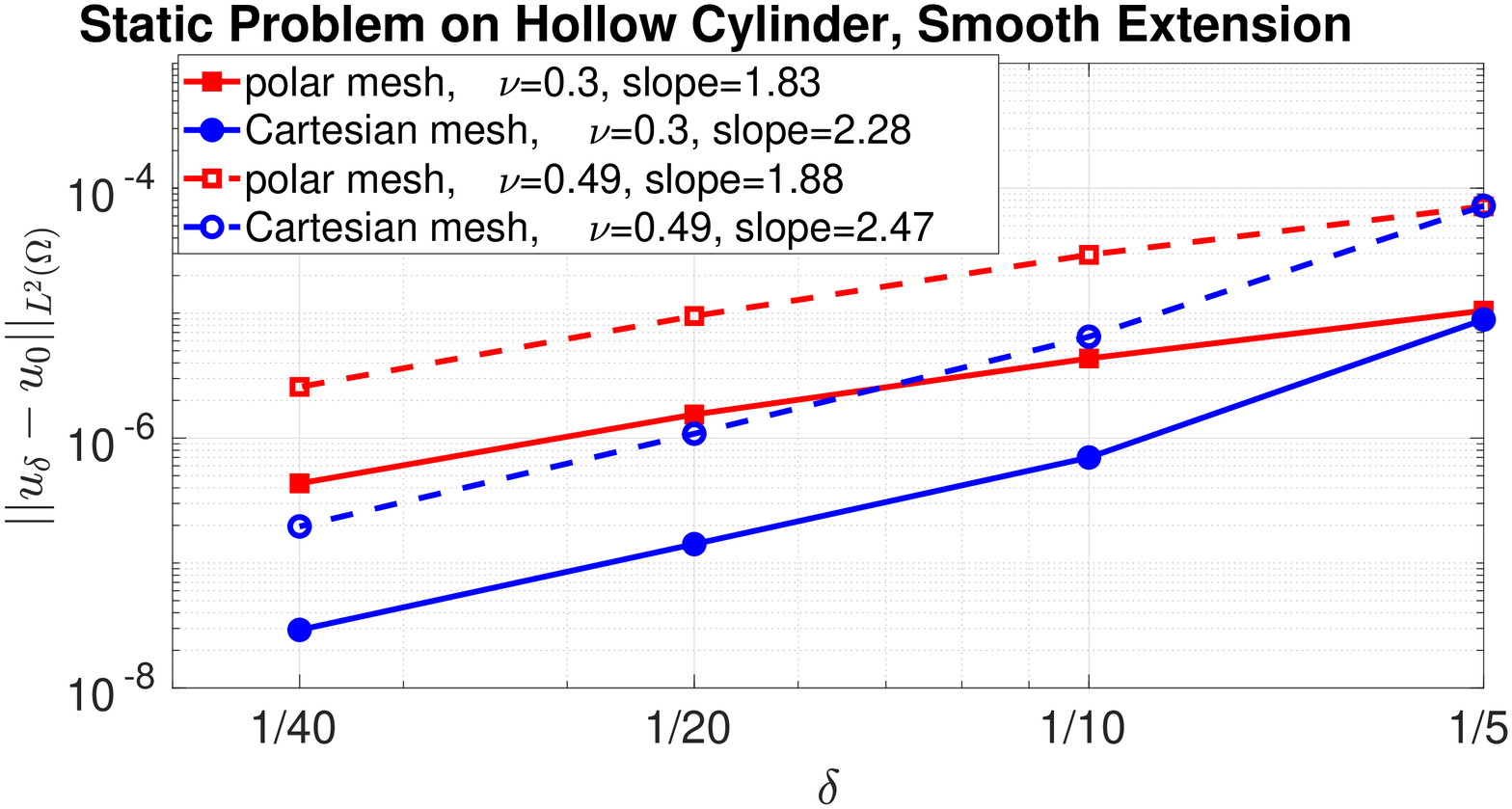}}
 \subfigure{\includegraphics[width=0.49\textwidth]{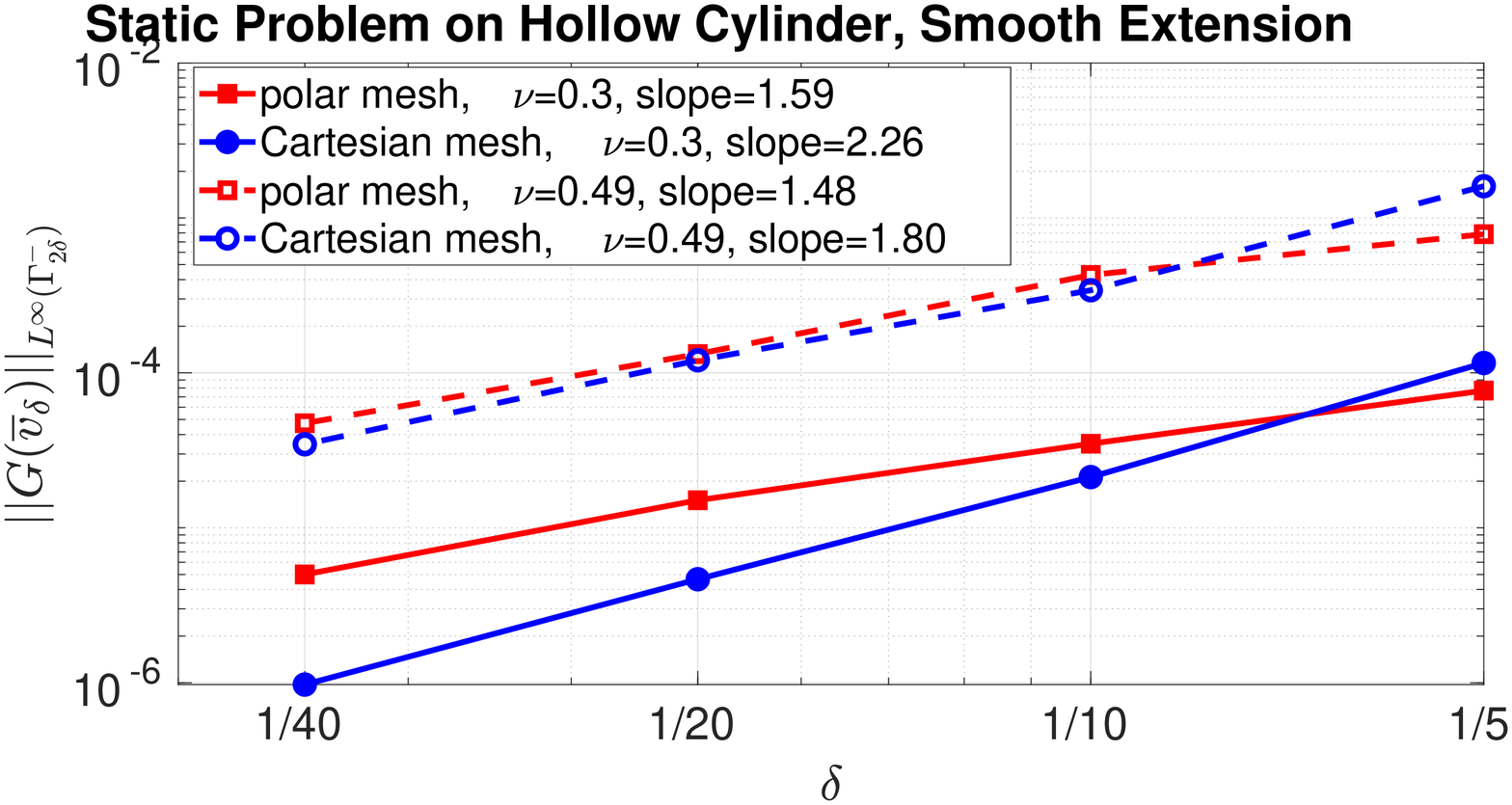}} 
 \caption{{\bf Static problem with smooth extension boundary condition.} Left: the $L^2(\omg;\re^d)$ difference between displacement $\ub_\delta$ and its local limit $\ub_0$. Right: the convergence of $\|\oG\ov{\vv}_\delta\|_{L^\infty(\bnd_{2\delta}^-)}$ in condition (A4').
 }
 \label{fig:sin_Diri_s}
\end{figure}

 \begin{figure}[!htb]\centering
 \subfigure{\includegraphics[width=0.49\textwidth]{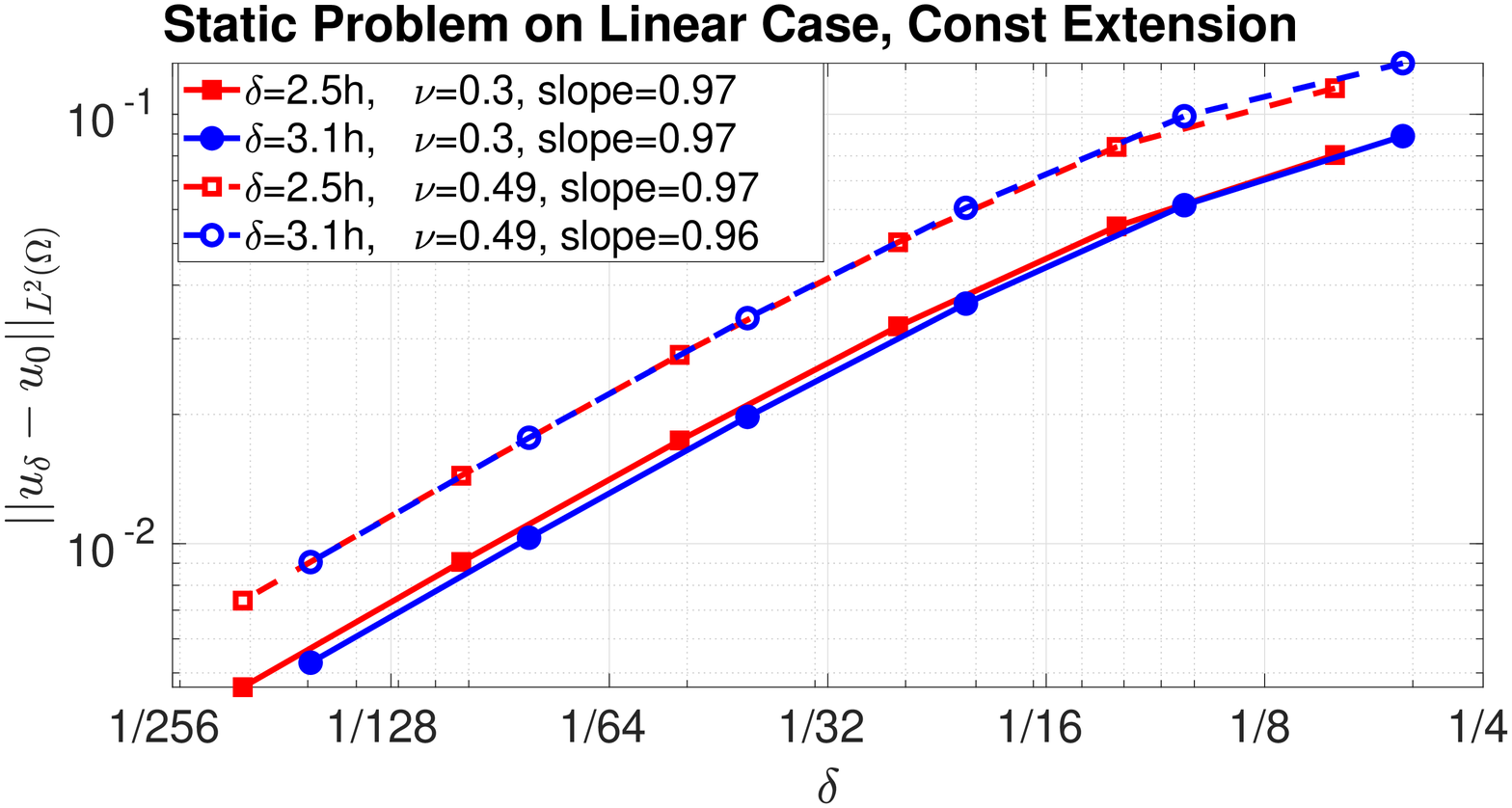}}
   \subfigure{\includegraphics[width=0.49\textwidth]{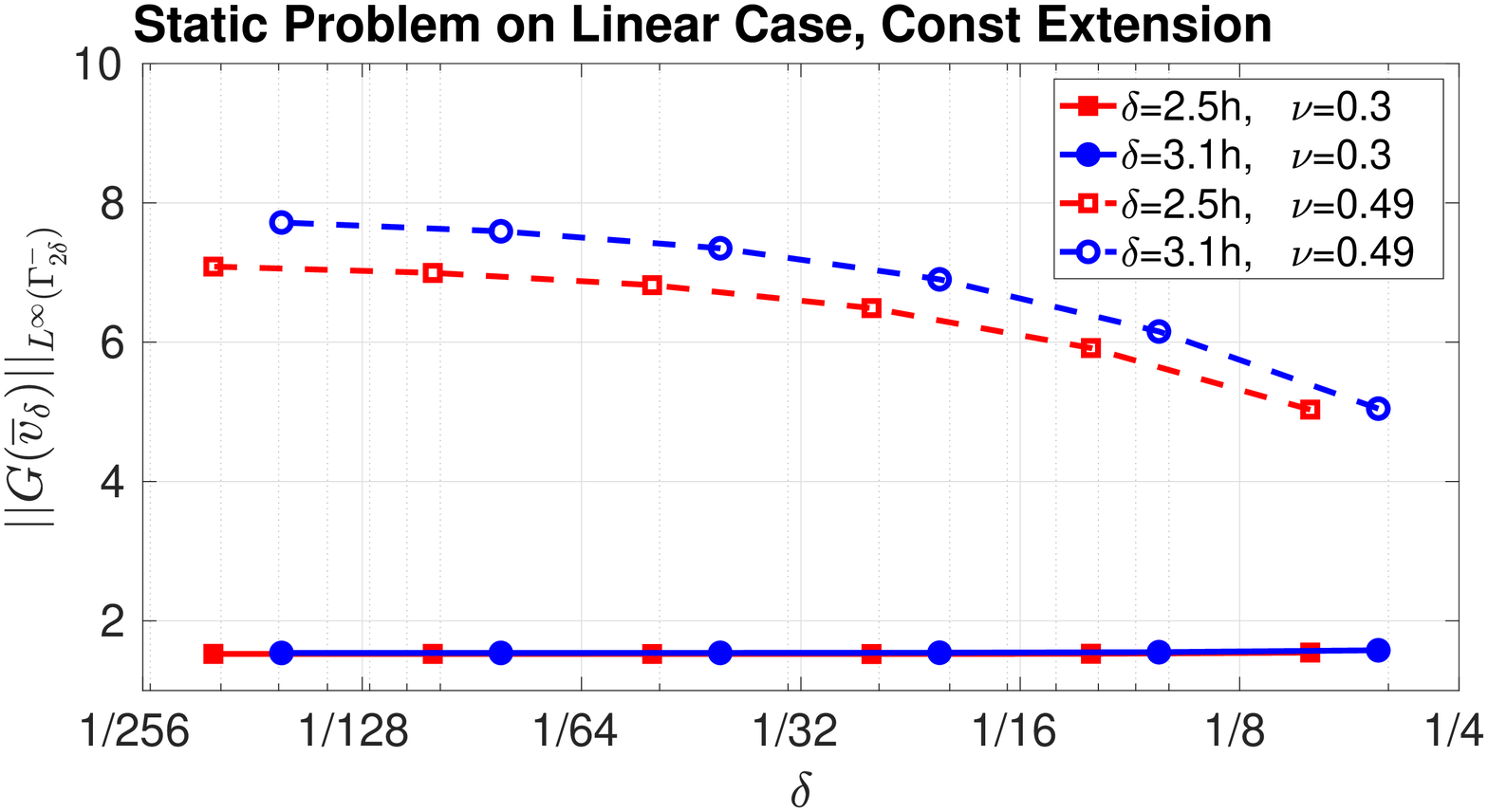}}
 \subfigure{\includegraphics[width=0.49\textwidth]{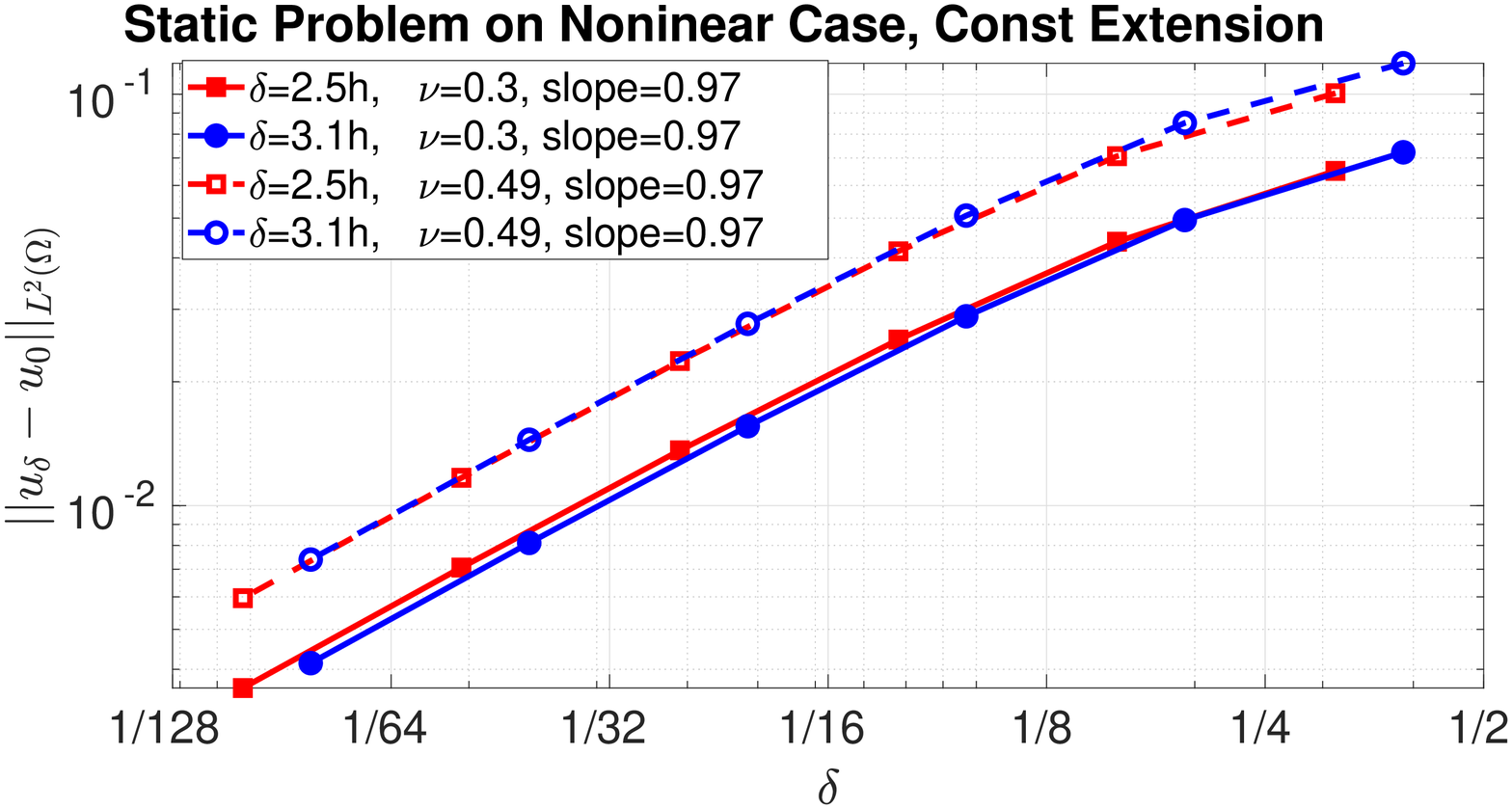}}
 \subfigure{\includegraphics[width=0.49\textwidth]{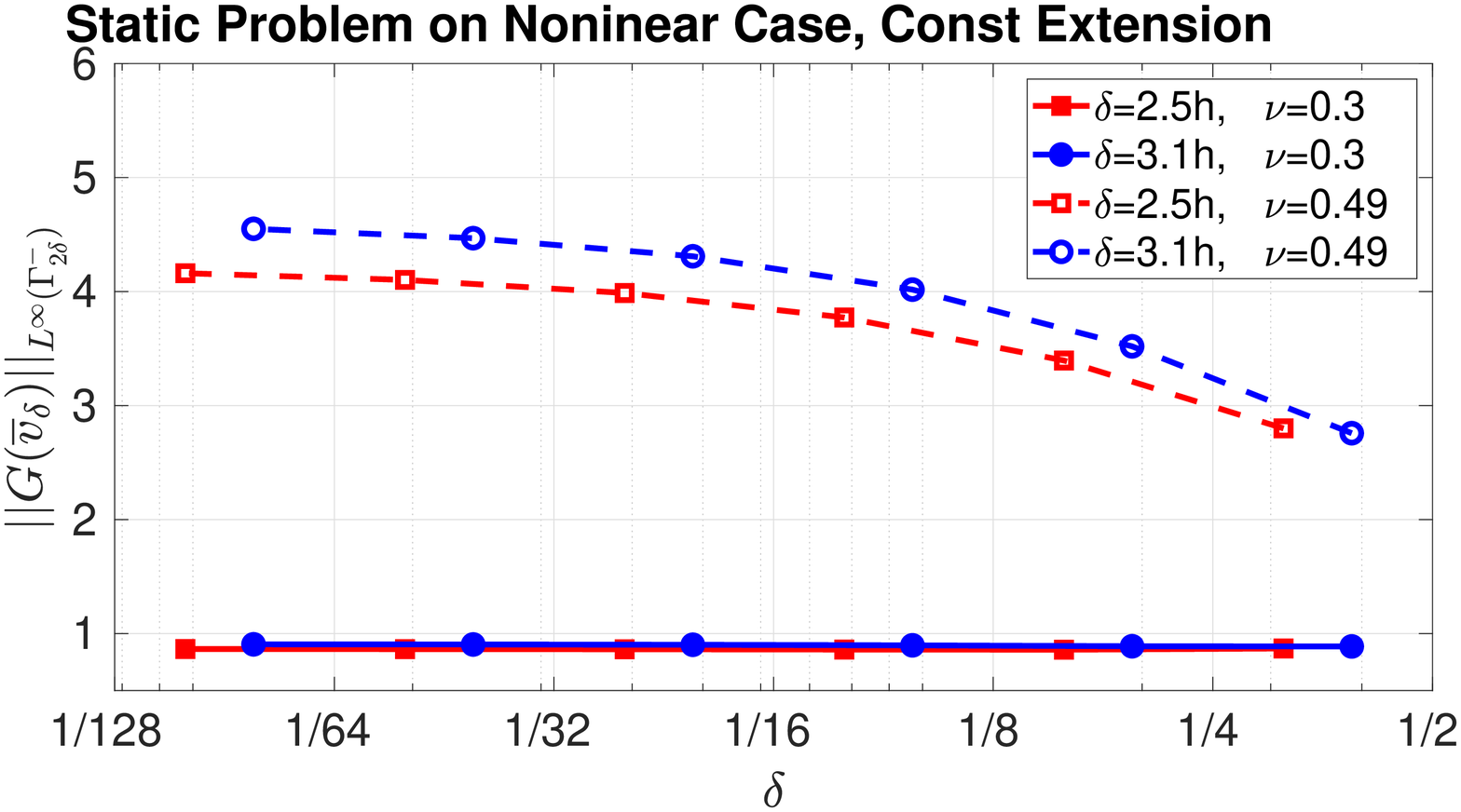}} 
 \subfigure{\includegraphics[width=0.49\textwidth]{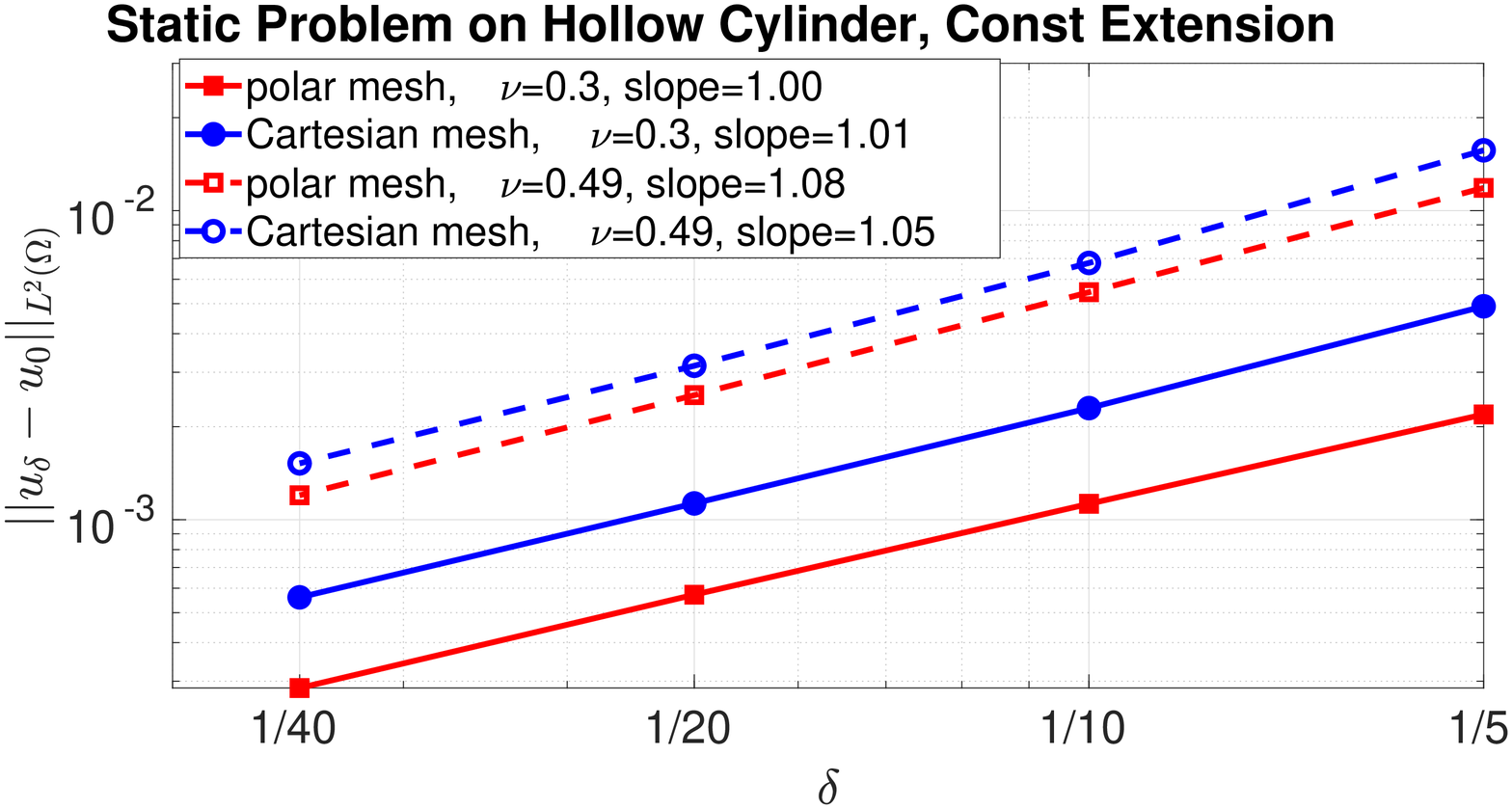}}
 \subfigure{\includegraphics[width=0.49\textwidth]{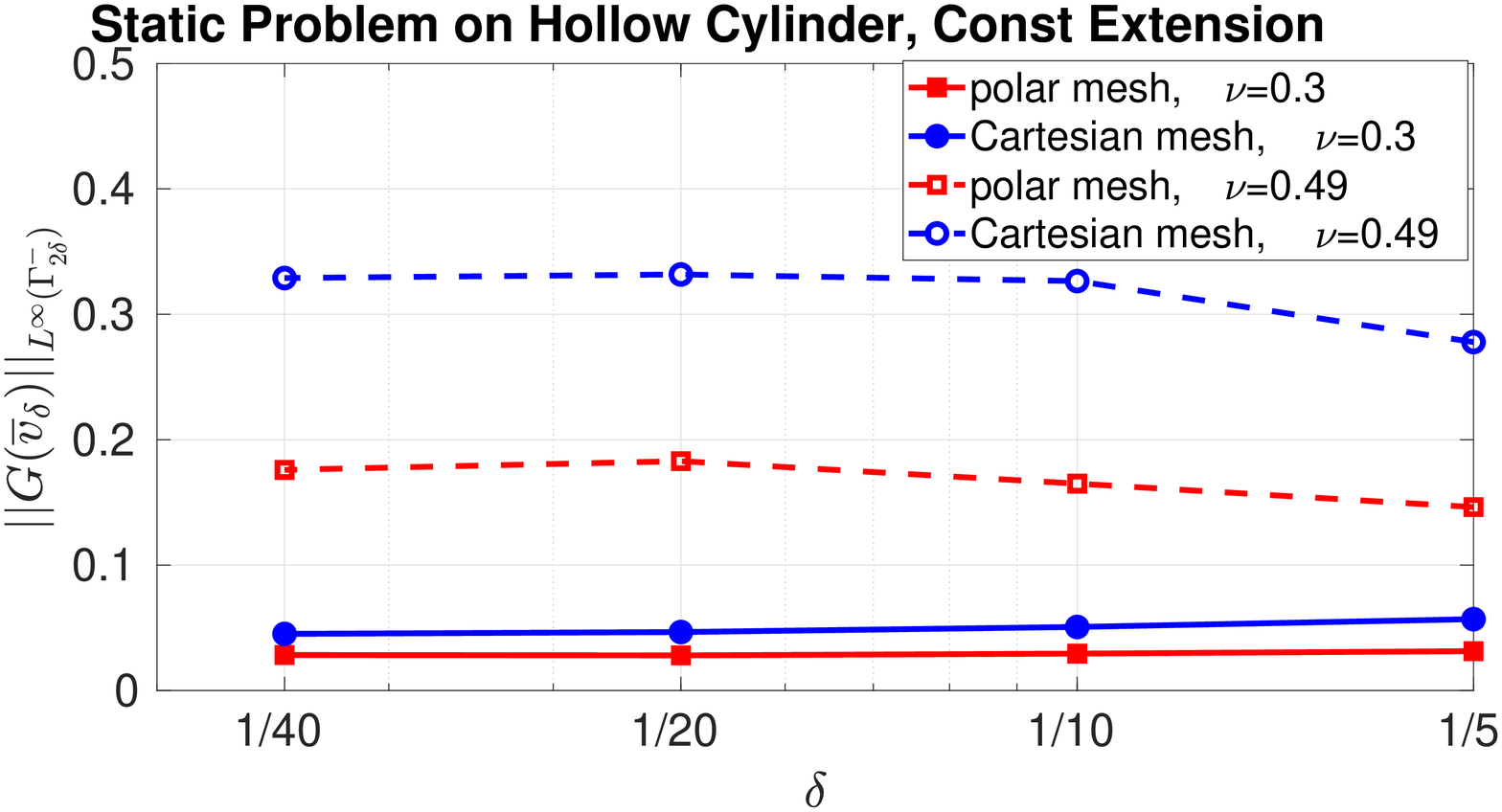}}
 \caption{{\bf Static problem with constant extension boundary condition.} Left: the $L^2(\omg;\re^d)$ difference between displacement $\ub_\delta$ and its local limit $\ub_0$. Right: the convergence of $\|\oG\ov{\vv}_\delta\|_{L^\infty(\bnd_{2\delta}^-)}$ in condition (A4').
 }
 \label{fig:linear_const_s}
\end{figure}

 \begin{figure}[!htb]\centering
 \subfigure{\includegraphics[width=0.49\textwidth]{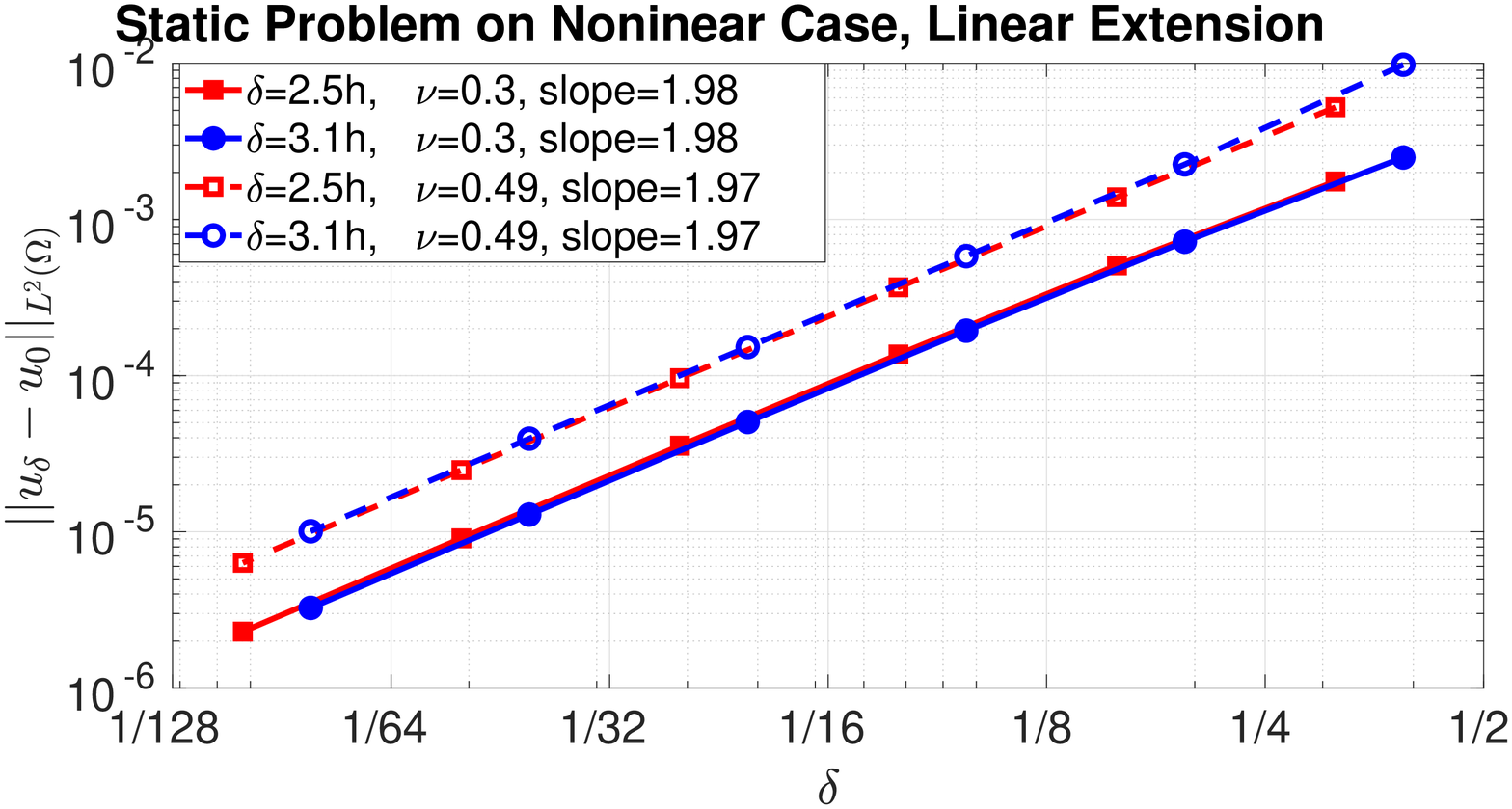}}
 \subfigure{\includegraphics[width=0.49\textwidth]{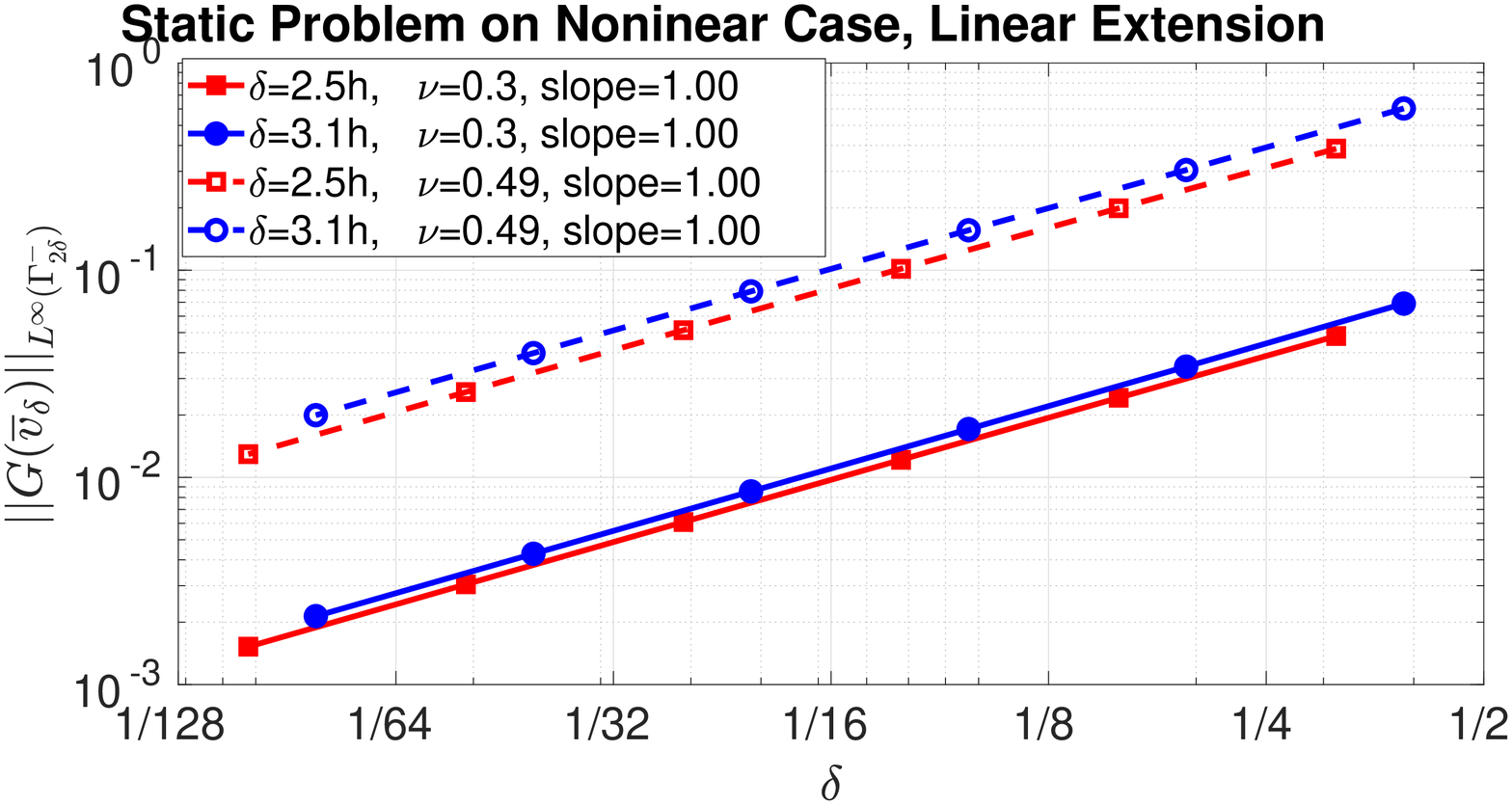}}
 \subfigure{\includegraphics[width=0.49\textwidth]{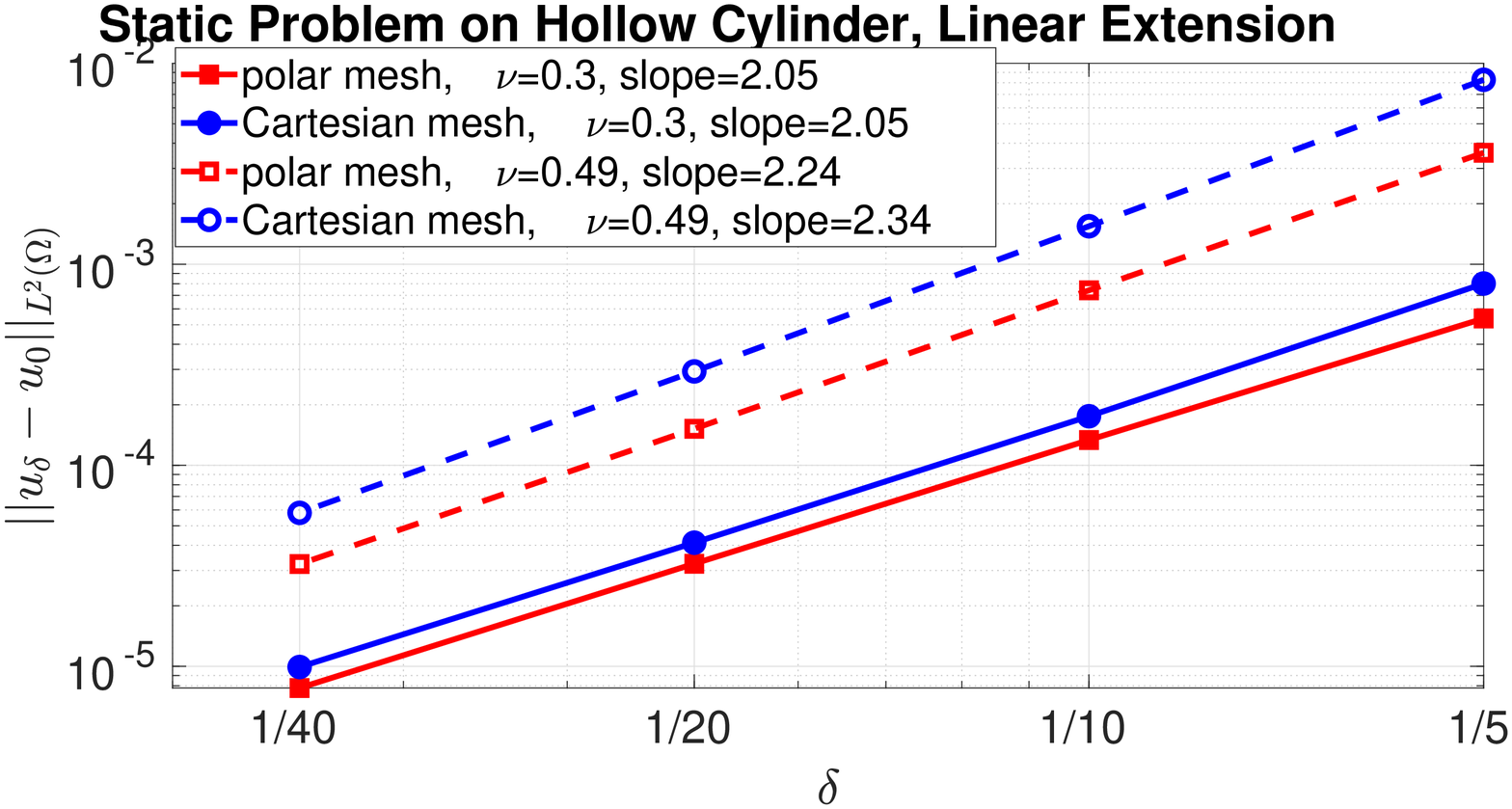}}
 \subfigure{\includegraphics[width=0.49\textwidth]{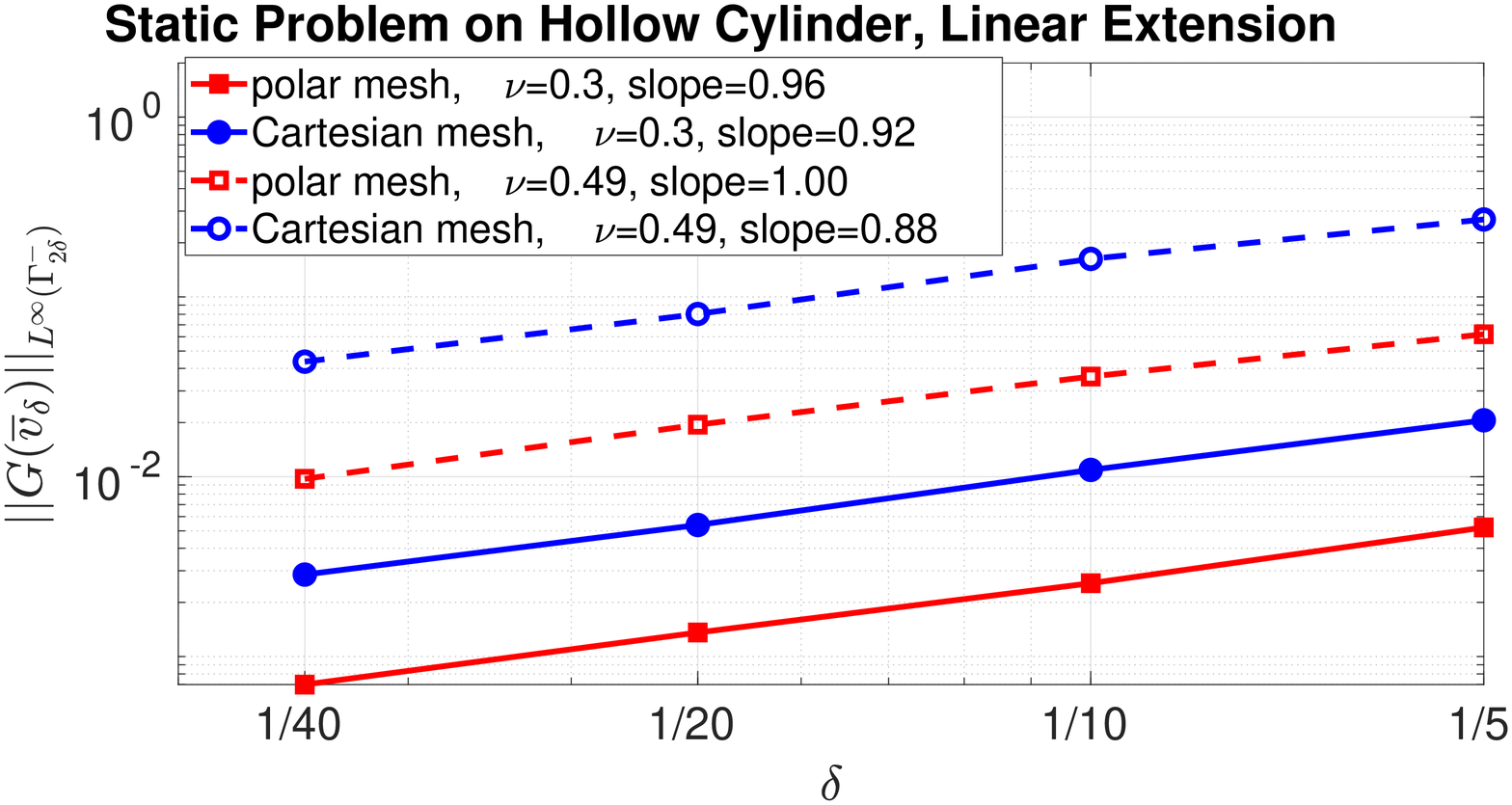}}
 \caption{{\bf Static problem with linear extension boundary condition.} Left: the $L^2(\omg;\re^d)$ difference between displacement $\ub_\delta$ and its local limit $\ub_0$. Right: the convergence of $\|\oG\ov{\vv}_\delta\|_{L^\infty(\bnd_{2\delta}^-)}$ in condition (A4').
 }
 \label{fig:sin_linear_s}
\end{figure}

To demonstrate the convergence in static LPS problems and the impacts of Dirichlet boundary conditions, three different settings are considered:
\begin{enumerate}
    \item {\bf Linear case}: We consider as linear patch test the displacement
$$\ub_0(x,y)=(3x+2y,-x+2y)$$
on a square domain $\Omega=[-1/4,1/4] \times [-1/4,1/4]$. Note that when the analytical $\ub_0$ is provided on $\bnd^+_{2\delta}$, in linear patch tests the local and nonlocal solutions coincide.
    \item {\bf Nonlinear case}: We consider a manufactured local solution  adopted from \cite{Yu_2012}:
    $$\ub_0(x,y) = [Ax+C\sin(bx),0],\qquad \fb(x,y)=[(\lambda+2\mu)b^2C\sin(bx),0].$$
    on a square domain $\Omega=[-1/2,-1/2] \times [-1/2,1/2]$. The parameters are taken as $A=0.9$, $C=1.4$ and $b=1.6$. 
\item {\bf Hollow cylinder case}: We consider the expansion of a hollow cylinder under an internal pressure $p_0=0.1$. The classical linear elasticity model predicts displacements given by
    $$\ub_0(x,y) = \left[Ax+\dfrac{Bx}{x^2+y^2},Ay+\dfrac{By}{x^2+y^2}\right]$$
    where
    $$A=\dfrac{(1+\nu)(1-2\nu)p_0R_0^2}{K(R_1^2-R_0^2)},\;B=\dfrac{(1+\nu)p_0R_0^2R_1^2}{K(R_1^2-R_0^2)},$$
    $R_0=1$ and $R_1=1.5$ are the interior and exterior radius of the hollow disk. An illustration of the hollow cylinder problem setting and displacement magnitudes are plotted in Figure \ref{fig:diskgrids}. 
\end{enumerate}

To establish the numerical convergence rates for the difference between $\ub_\delta$ and $\ub_0$ with given extensions of different degrees of regularity, for each setting three boundary extension strategies are considered to provide data on the exterior collar $\Gamma^+_{2\delta}$:
\begin{enumerate}
\item {\bf Smooth extension}: A smooth extended local solution is provided as $\ub_D(\xb)$ for $\xb\in\Gamma^+_{2\delta}$. In particular, we set:
$$\ub_D(\xb)=\ub_0(\xb),\quad \forall\xb\text{ in }\Gamma^+_{2\delta}.$$

Results of solution convergence with smooth extension are plotted in Figure \ref{fig:sin_Diri_s}. We observe that the numerical solution passes the patch test within machine precision. For the nonlinear case and the hollow disk case, we observe second-order $L^2(\omg;\re^d)$-norm convergence for displacements in both compressible and nearly-incompressible materials, which is consistent with the theoretical bound discussed in Theorem \ref{thm:conv} (a) and Remark \ref{rmk1}.

\item {\bf Constant extension}: The naive fictitious node method is employed such that $\ub_D(\xb)$ is defined by the surface (local) data at its corresponding projection on $\partial\Omega$. In particular,
$$\vu_D(\vx)=\vu_0(\overline{\vx}),\quad \forall \xb \text{ in }\bnd^+_{2\delta},$$
where $\overline{\vx}:=\underset{\yb\in\partial\omg}{\text{argmin}}\,\text{dist}(\yb,\xb)$. We note that the constant extension on the boundary region doesn't pass the linear patch test theoretically. {Results for the three numerical cases are plotted in Figure \ref{fig:linear_const_s}, where we can see that $\|\oG\ov{\vv}_\delta\|_{L^\infty(\bnd_{2\delta}^-)}$ is uniformly bounded but not converging in all cases, i.e., $\gamma'=0$. Therefore, Theorem \ref{thm:conv} (b) and Corollary \ref{cor:conv} provides an estimated convergence bound for $\vertii{\ub_\delta-\ub_0}_{L^2(\omg;\re^d)}$ as $O(\delta^{1/2})$. Numerically, we observe an $O(\delta)$ convergence for displacements, which is slightly better than the $O(\delta^{1/2})$ theoretical convergence rate and suggests that the theoretical rate of convergence in Theorem \ref{thm:conv} (b) might be sub-optimal.} 

\item {\bf Linear extension}: We consider the mirror-based fictitious node methods which extends the surface (local) and the interior data into the exterior layer $\Gamma^+_{2\delta}$ linearly. Particularly, we set
$$\ub_D(\xb)=2\ub_0(\overline{\xb})-\ub_\delta(2\overline{\xb}-\xb),\quad \forall \xb \text{ in }\bnd^+_{2\delta},$$
where $\overline{\vx}:=\underset{\yb\in\partial\omg}{\text{argmin}}\,\text{dist}(\yb,\xb)$. 
For the linear patch test case, machine precision accuracy is again observed. For the other two cases, in Figures \ref{fig:sin_linear_s} we present convergence results, where we observe $O(\delta)$ convergence for $\|\oG\ov{\vv}_\delta\|_{L^\infty(\bnd_{2\delta}^-)}$ as $\delta\rightarrow 0$,  i.e., $\gamma'=1$. Therefore, an $O(\delta^{3/2})$ theoretical convergence rate is provided by Theorem \ref{thm:conv} (b) and Corollary \ref{cor:conv}. On the other hand, numerically quadratic convergence is observed for $\ub_\delta-\ub_0$ in the $L^2$ norm, which is again slightly better than the $O(\delta^{3/2})$ theoretical convergence rate.
\end{enumerate}

\subsection{Dynamic problem with Dirichlet boundary conditions}\label{sec:dynamicnonlocal}

 \begin{figure}[!htb]\centering
 \subfigure{\includegraphics[width=0.49\textwidth]{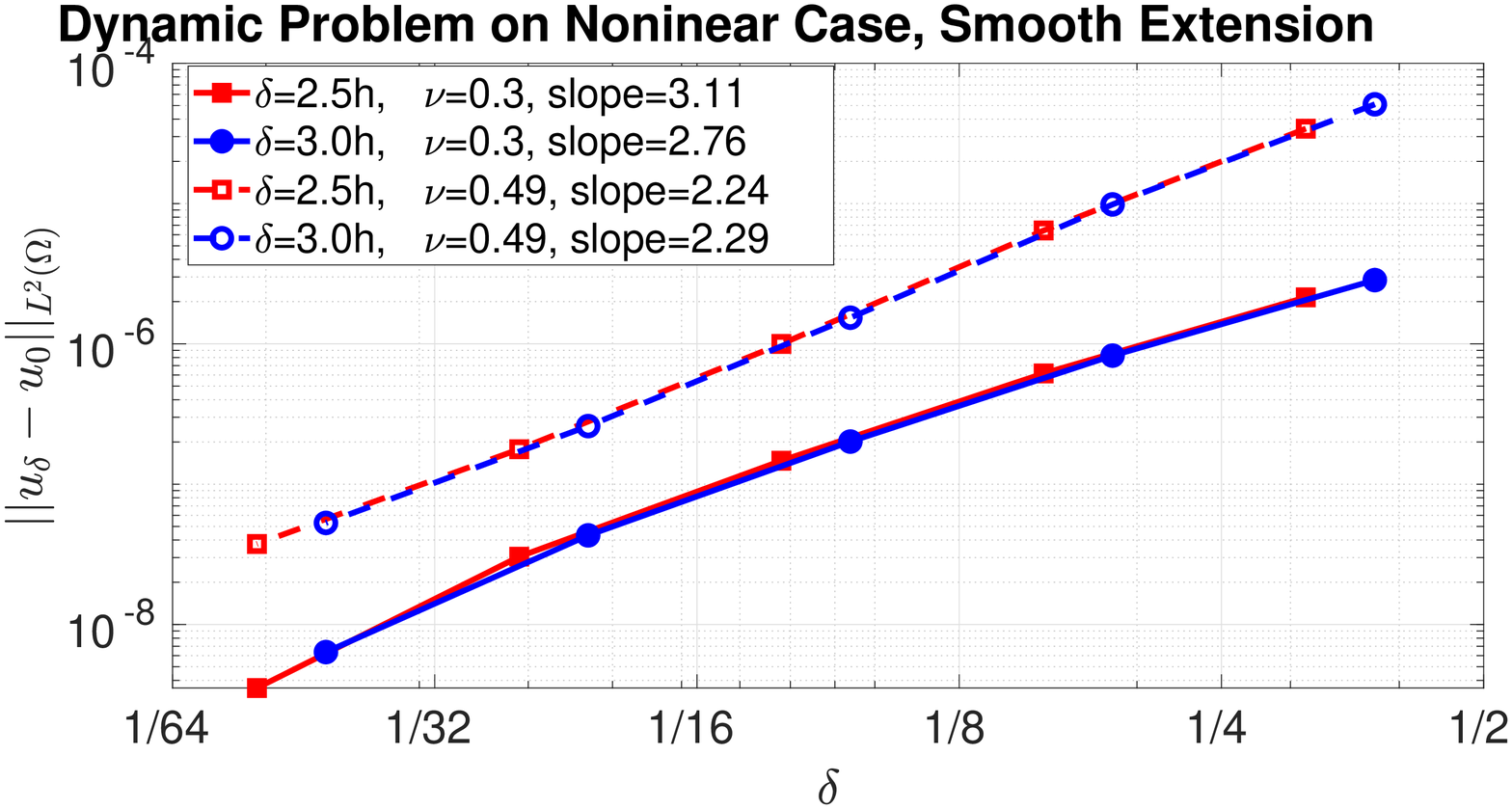}}
 \subfigure{\includegraphics[width=0.49\textwidth]{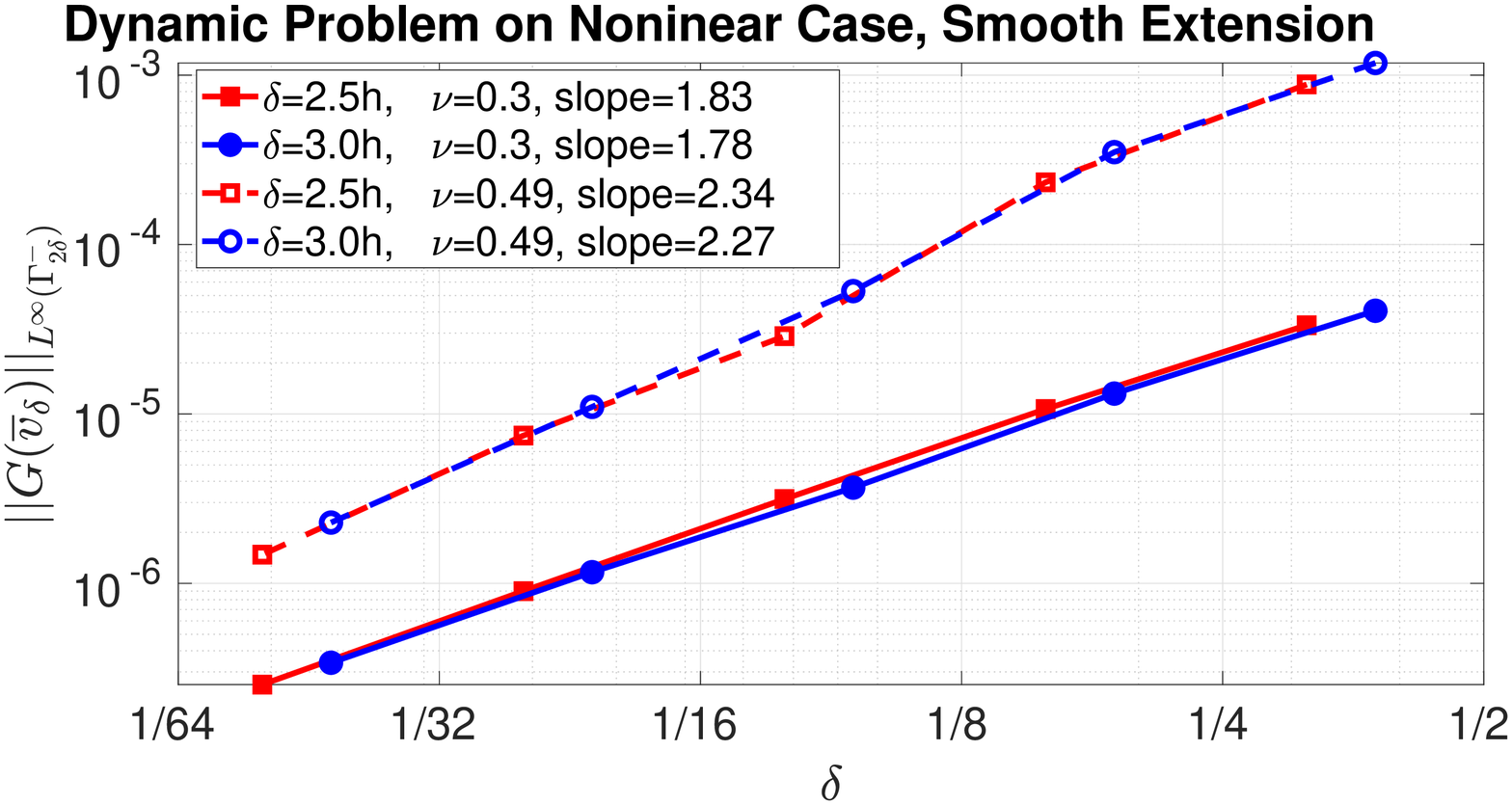}}
 \subfigure{\includegraphics[width=0.49\textwidth]{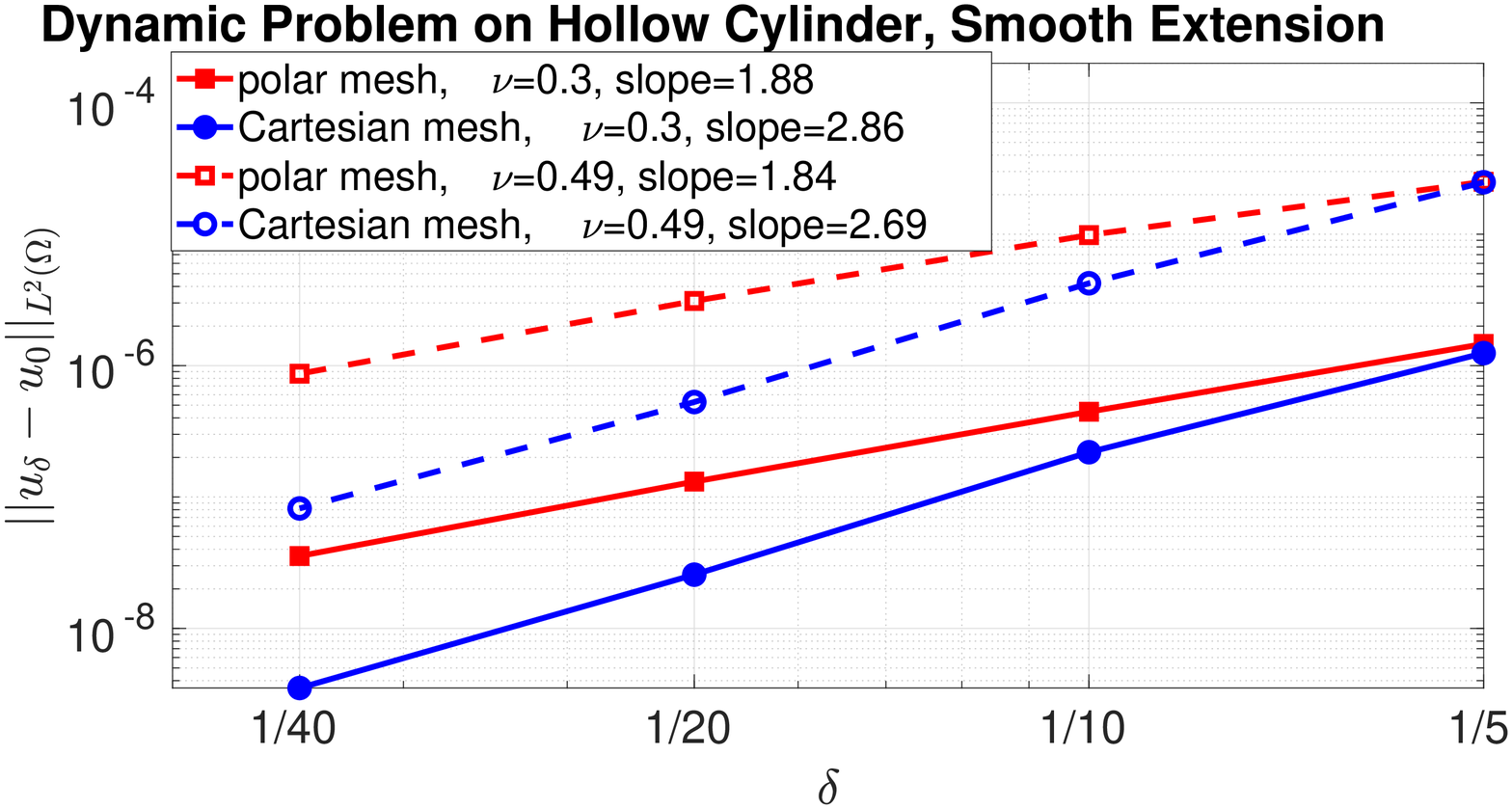}}
 \subfigure{\includegraphics[width=0.49\textwidth]{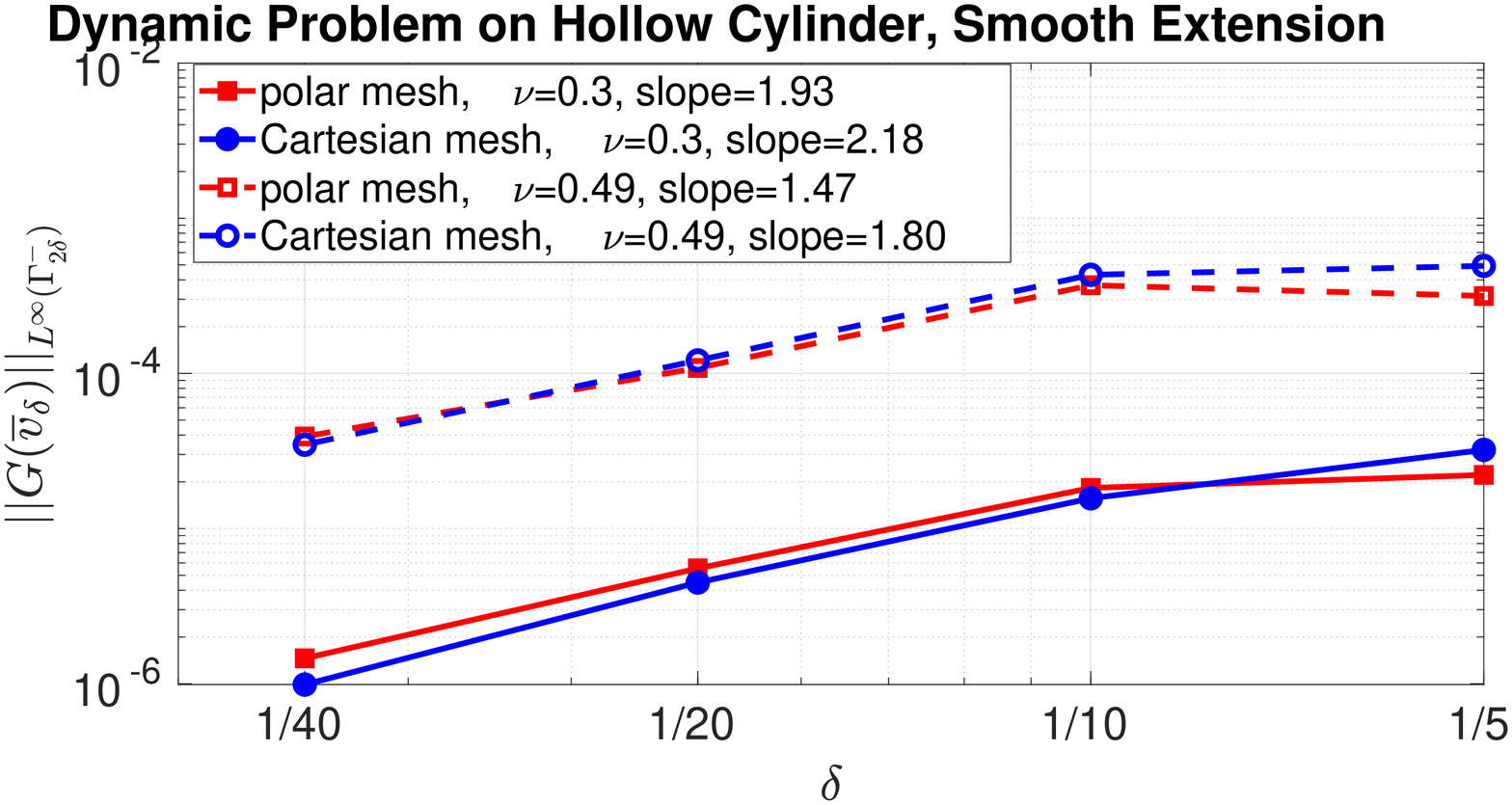}}
 \caption{{\bf Dynamic problem with smooth extension boundary condition.} Left: the $L^2(\omg;\re^d)$ difference between displacement $\ub_\delta(T,\cdot)$ and its local limit $\ub_0(T,\cdot)$. Right: the convergence of $\|\oG\ov{\vv}_\delta\|_{L^\infty(\bnd_{2\delta}^-)}$ in condition (A4').
 }
\label{fig:sin_Diri_d}
\end{figure}

 \begin{figure}[!htb]\centering
 \subfigure{\includegraphics[width=0.49\textwidth]{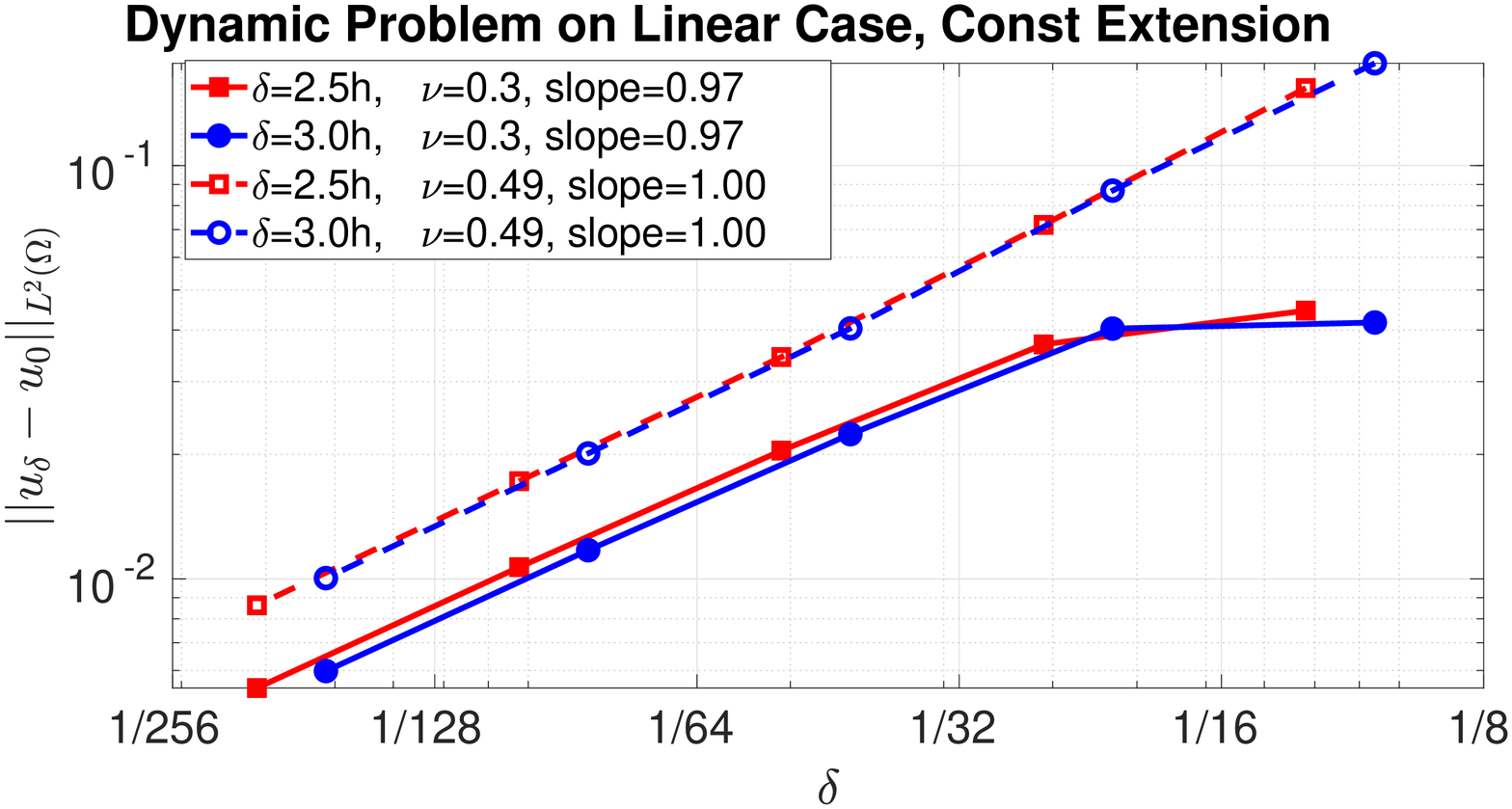}}
 \subfigure{\includegraphics[width=0.49\textwidth]{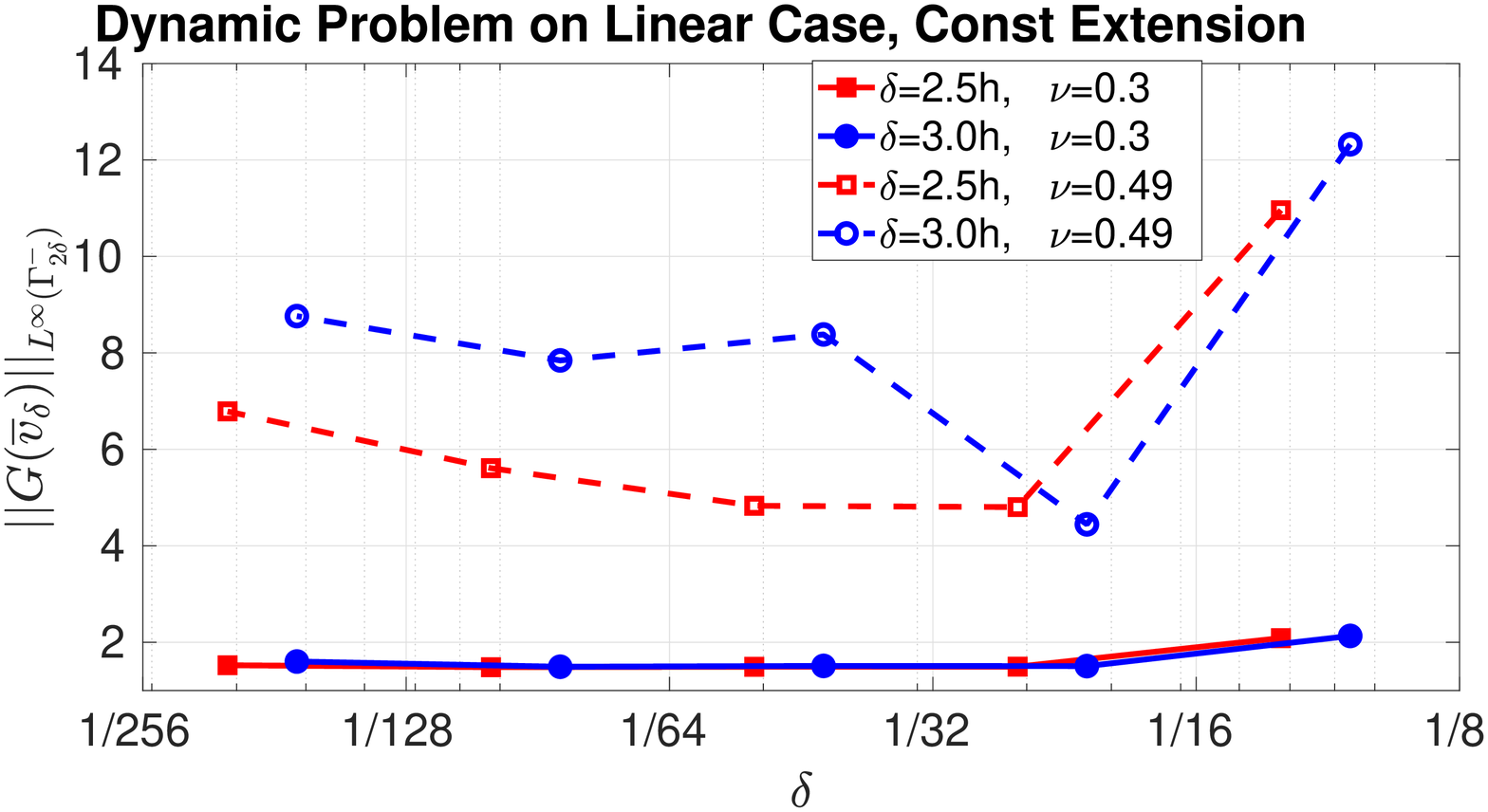}}
 \subfigure{\includegraphics[width=0.49\textwidth]{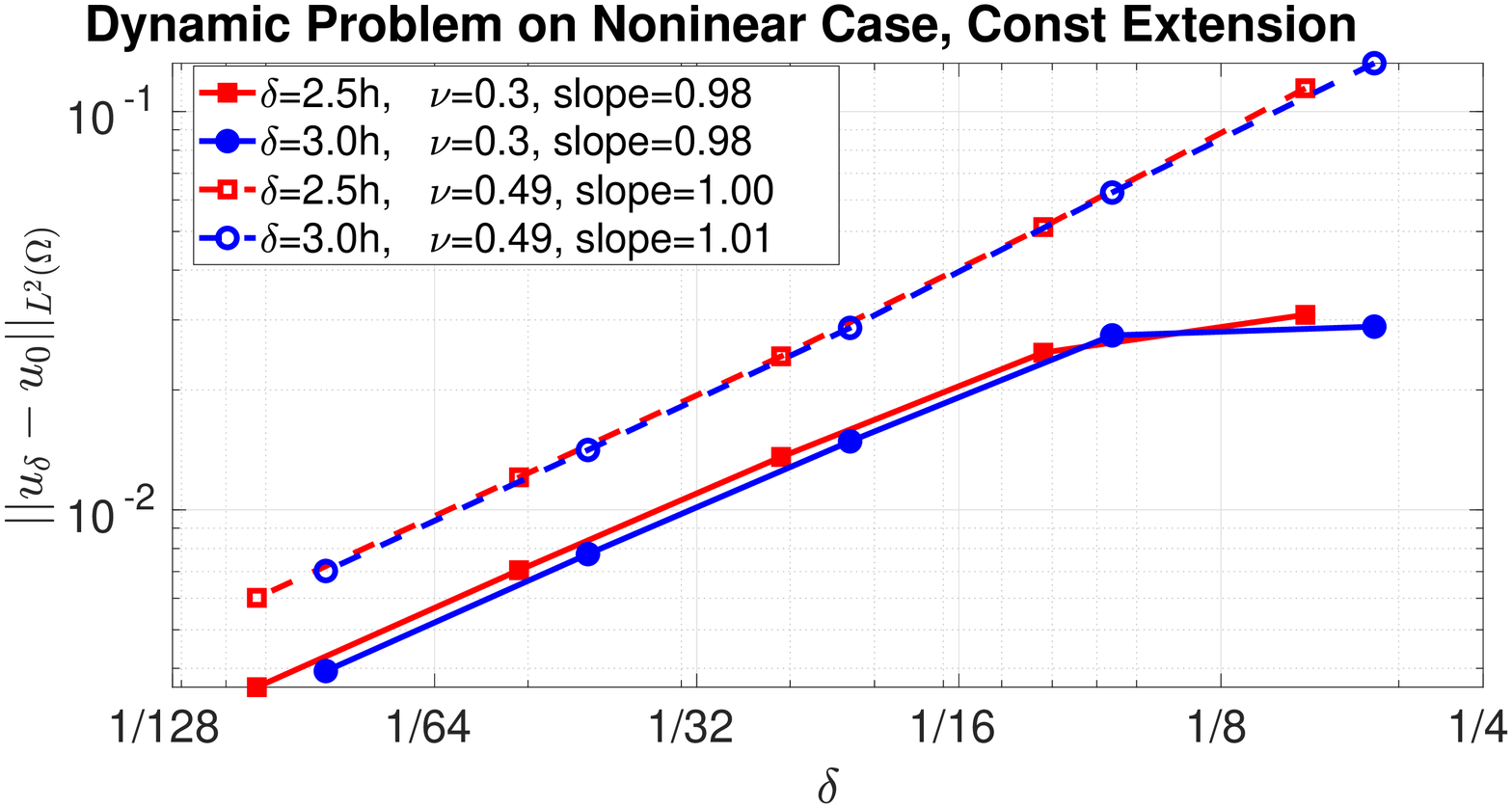}}
 \subfigure{\includegraphics[width=0.49\textwidth]{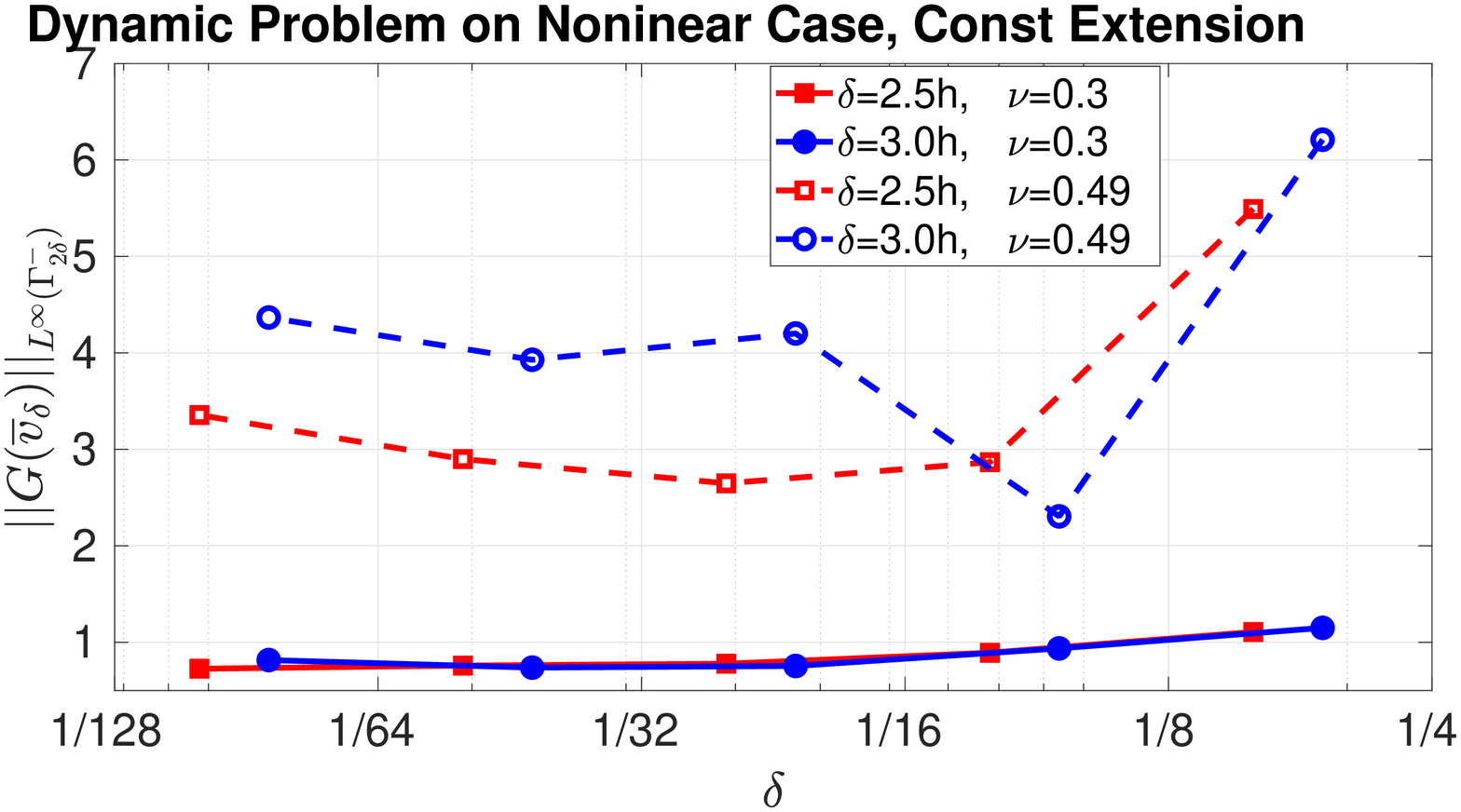}}
 \subfigure{\includegraphics[width=0.49\textwidth]{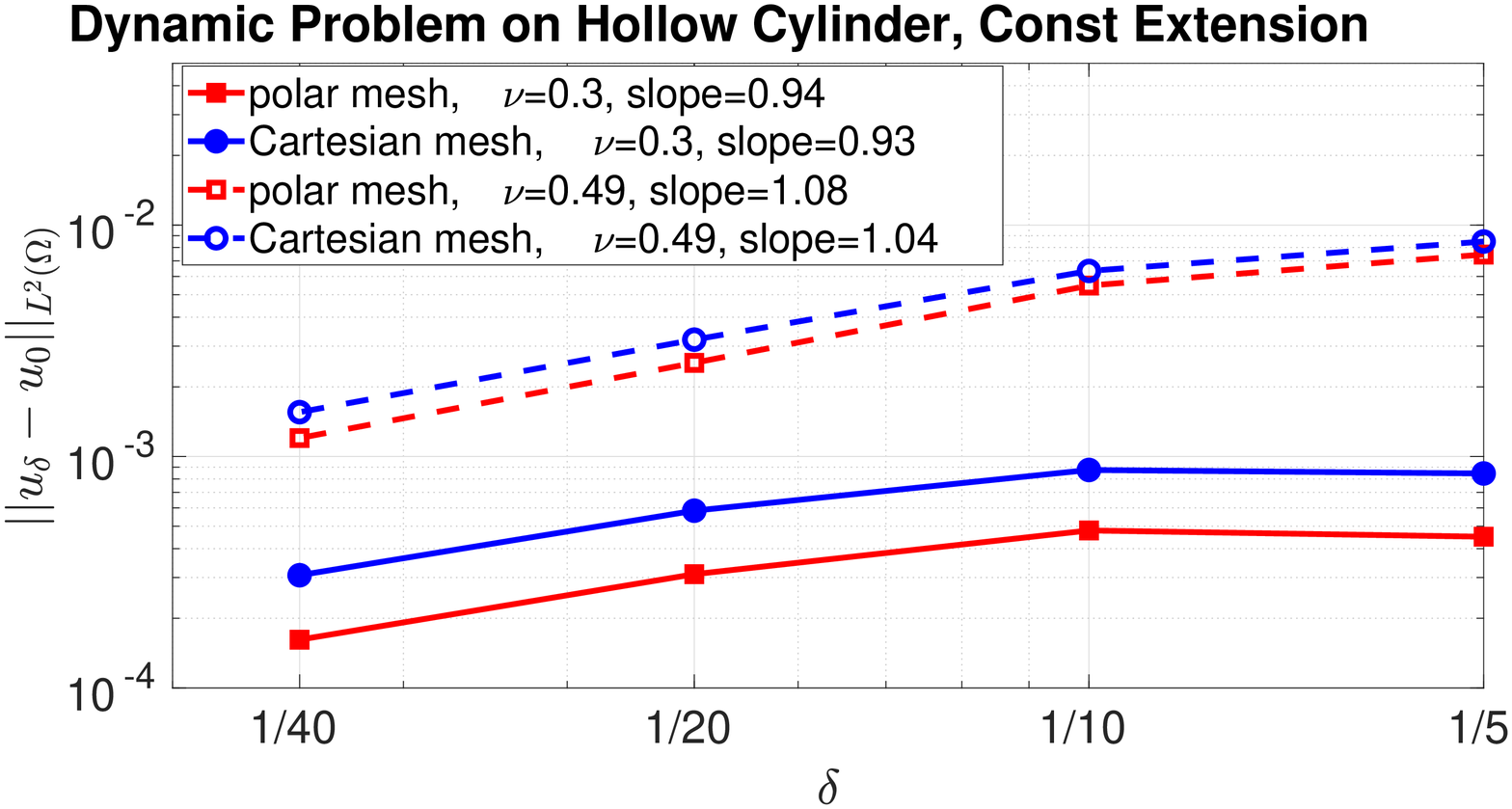}}
 \subfigure{\includegraphics[width=0.49\textwidth]{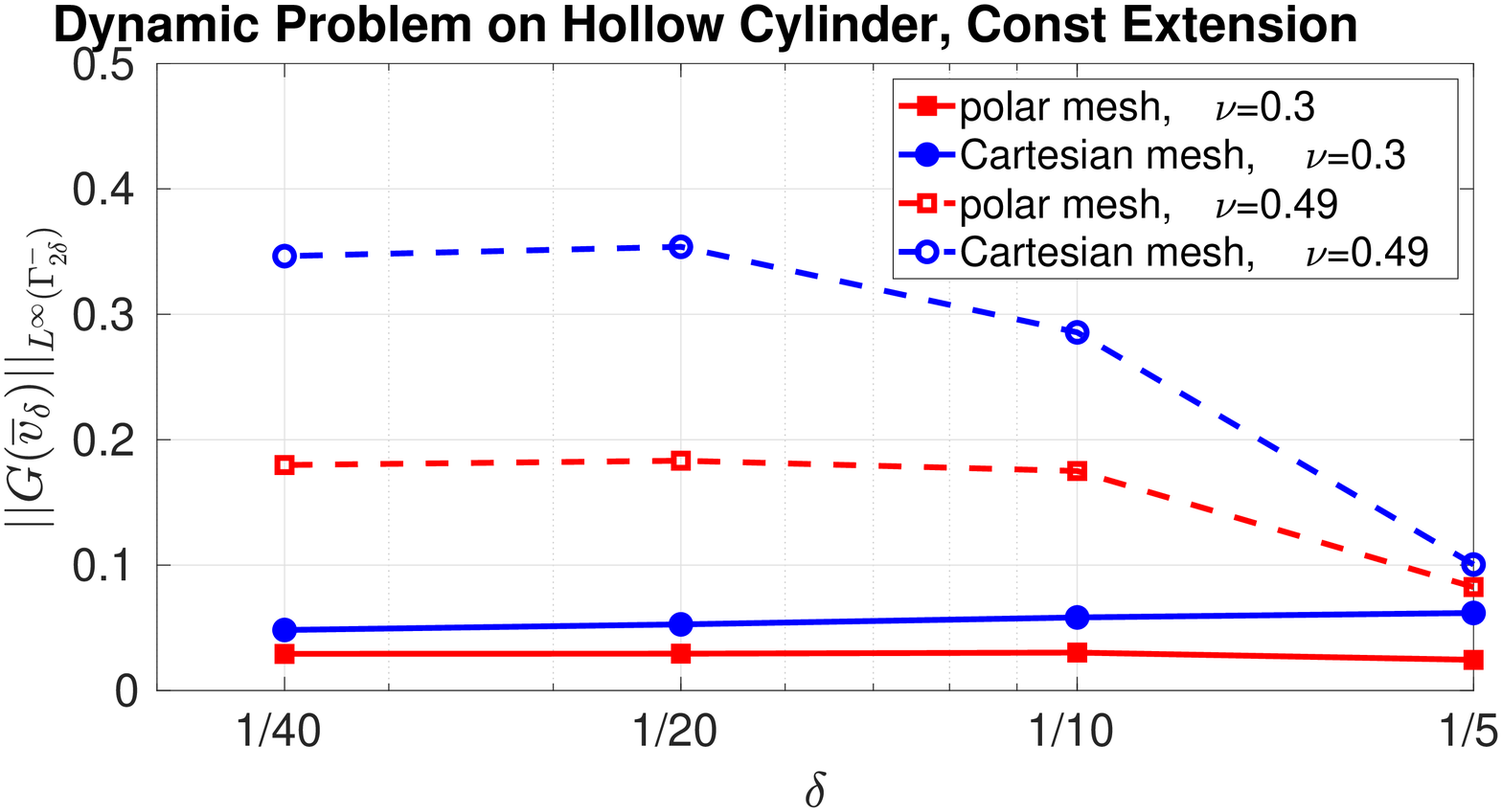}} 
 \caption{{\bf Dynamic problem with constant extension boundary condition.} Left: the $L^2(\omg;\re^d)$ difference between displacement $\ub_\delta(T,\cdot)$ and its local limit $\ub_0(T,\cdot)$. Right: the convergence of $\|\oG\ov{\vv}_\delta\|_{L^\infty(\bnd_{2\delta}^-)}$ in condition (A4').
 }
 \label{fig:linear_const_d}
\end{figure}

 \begin{figure}[!htb]\centering
 \subfigure{\includegraphics[width=0.49\textwidth]{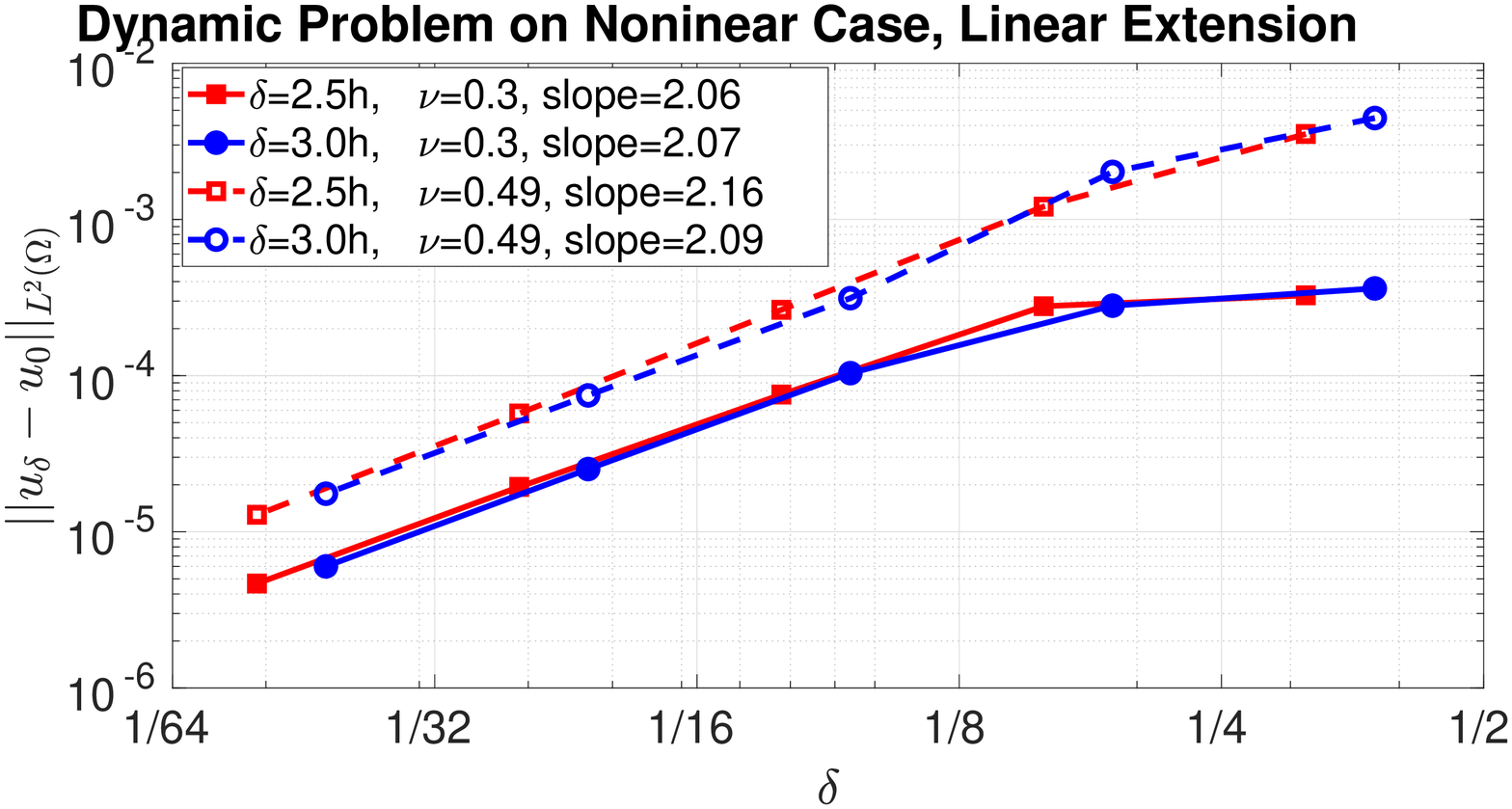}}
 \subfigure{\includegraphics[width=0.49\textwidth]{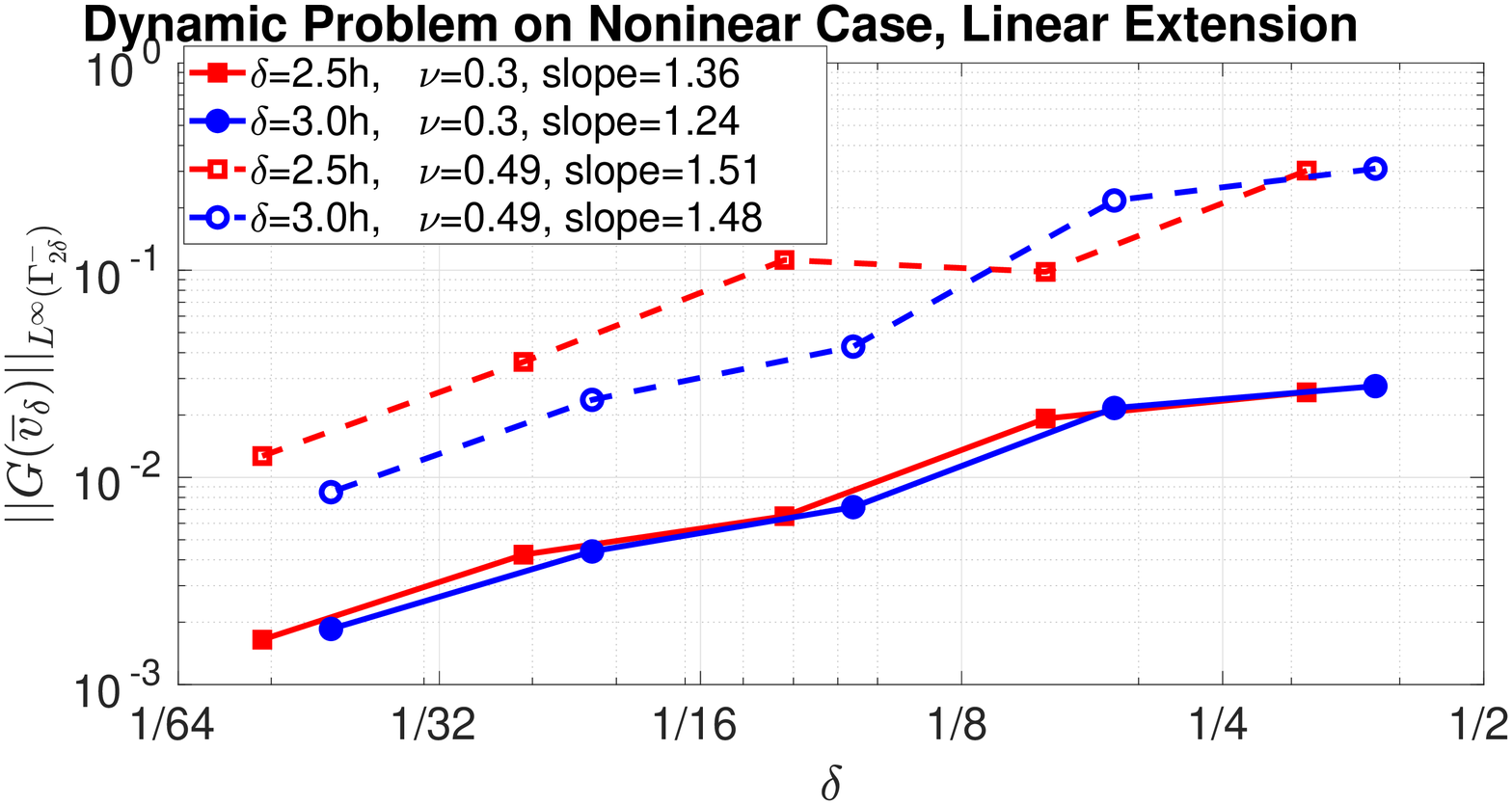}}
 \subfigure{\includegraphics[width=0.49\textwidth]{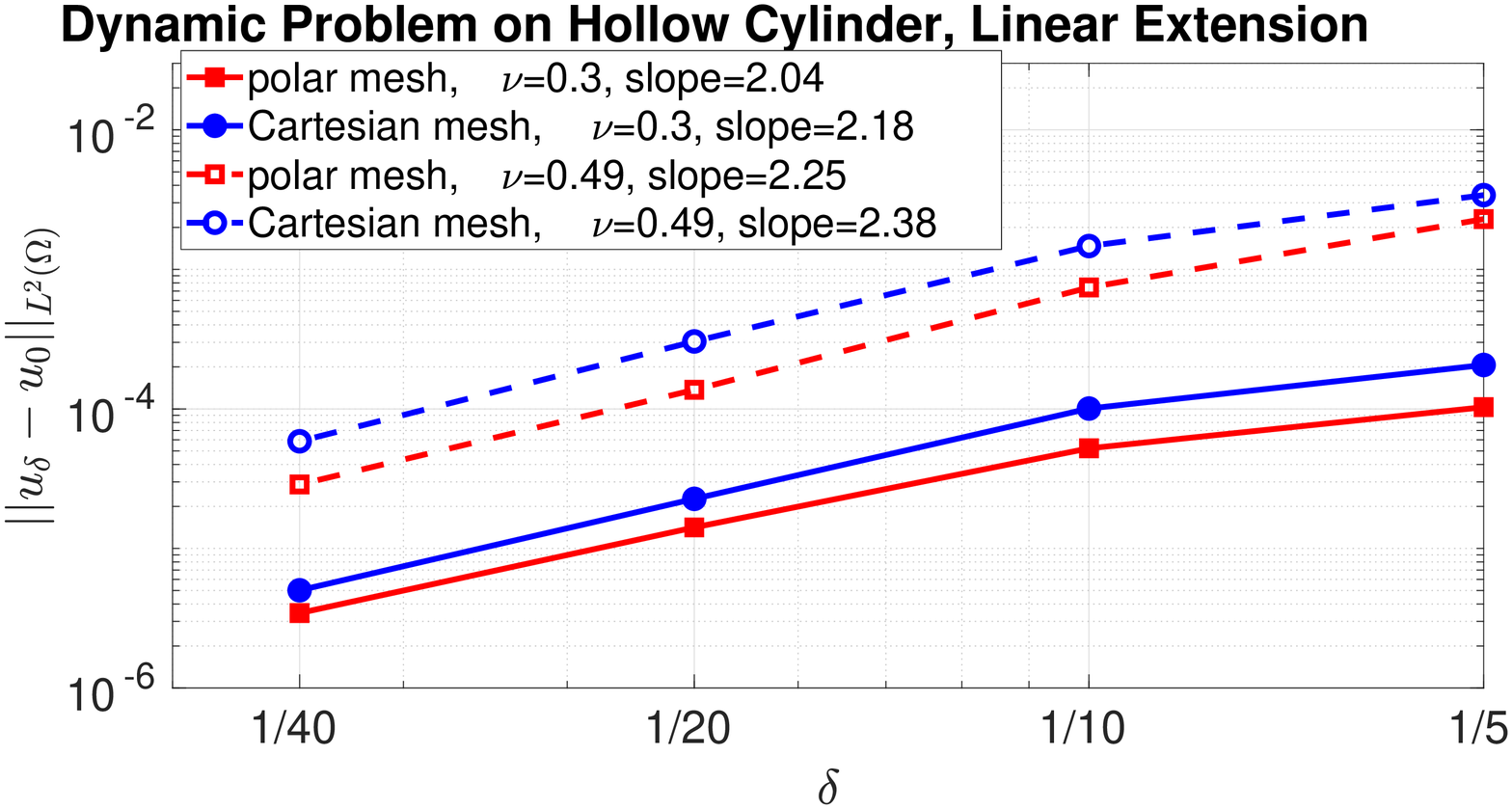}}
   \subfigure{\includegraphics[width=0.49\textwidth]{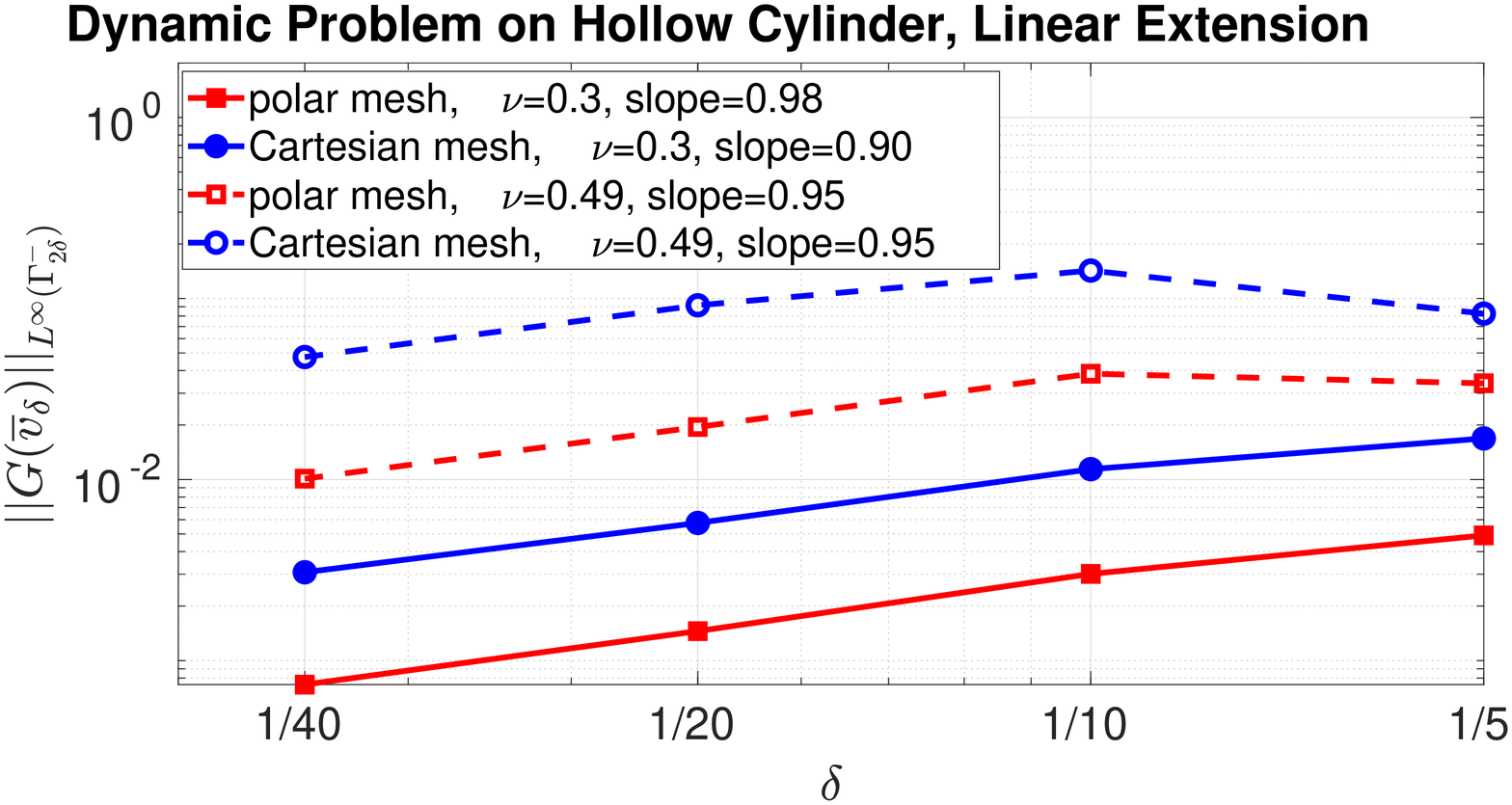}}  
 \caption{{\bf Dynamic problem with linear extension boundary condition.} Left: the $L^2(\omg;\re^d)$ difference between displacement $\ub_\delta(T,\cdot)$ and its local limit $\ub_0(T,\cdot)$. Right: the convergence of $\|\oG\ov{\vv}_\delta\|_{L^\infty(\bnd_{2\delta}^-)}$ in condition (A4').
 }
 \label{fig:sin_linear_d}
\end{figure}

With density $\rho=1.0$, in dynamic LPS problems we also consider three different settings. To observe the convergence with respect to $\delta$, in all cases we take a sufficiently small time step size $\Delta t=0.01$ and simulate until $T=0.1$.
\begin{enumerate}
    \item {\bf Linear case}: We consider as linear patch test the displacement
$$\ub_0(t,x,y)=(t+3x+2y,t-x+2y)$$
on a square domain $\Omega=[-1/4,1/4] \times [-1/4,1/4]$. 
    \item {\bf Nonlinear case}: We consider a manufactured local solution adopted from \cite{Yu_2012}:
    $$\ub_0(t,x,y) = [Ax+B\sin(at)x+C\sin(bx),0]$$
    on a square domain $\Omega=[-1/2,1/2] \times [-1/2,1/2]$ with the external loading 
    $$\fb(t,x,y)=[(\lambda+2\mu)b^2C\sin(bx)-\rho a^2Bx\sin(at),0].$$
    The parameters are taken as $A=0.9$, $B=0.1$, $C=1.4$, $a=1.0$ and $b=1.2$.
    \item {\bf Hollow cylinder case}: We consider the expansion of a hollow cylinder under an increasing internal pressure $p_0t$. Under this setting the classical solution $\ub_0$ is given by
    $$\ub_0(t,x,y) = \left[Axt+\dfrac{Bxt}{x^2+y^2},Ayt+\dfrac{Byt}{x^2+y^2}\right]$$
    where $p_0$, $A$, $B$ and $R_0$, $R_1$ are as defined in the static hollow cylinder case.
\end{enumerate}

Similar as in static LPS problems,  three  boundary  extension strategies  are  considered to study the numerical convergence rates with given extensions of different degrees of regularity:
\begin{enumerate}
\item {\bf Smooth extension}: We set 
$$\ub_D(t,\xb)=\ub_0(t,\xb),\quad \forall\xb\text{ in }\Gamma^+_{2\delta},t>0.$$
It is observed that the numerical solution passes the linear patch test to within machine precision. In Figures \ref{fig:sin_Diri_d} we plot the numerical results for the nonlinear and hollow cylinder cases, where second-order $L^2(\omg;\re^d)$-norm convergence for displacements is again observed, which is consistent with the analysis provided in Theorem \ref{thm:dynamicconv}.

\item {\bf Constant extension}: We set
$$\vu_D(t,\vx)=\vu_0(t,\overline{\vx}),\quad \forall \xb \text{ in }\bnd^+_{2\delta},t>0,$$
where $\overline{\vx}:=\underset{\yb\in\partial\omg}{\text{argmin}}\,\text{dist}(\yb,\xb)$. Numerical convergence results for all three numerical cases are provided in Figure \ref{fig:linear_const_d}, where we can see that the convergence rates are very similar to the static cases: $\|\oG\ov{\vv}_\delta\|_{L^\infty(\bnd_{2\delta}^-)}$ is uniformly bounded with $\gamma'=0$. Therefore, the theoretical analysis in Theorem \ref{thm:dynamicconv} provides an $O(\delta^{1/2})$ bound for $\vertii{\ub_\delta-\ub_0}_{L^2(\omg;\re^d)}$, while numerically we observe first order convergence.

\item {\bf Linear extension}: We set
$$\ub_D(t,\xb)=2\ub_0(t,\overline{\xb})-\ub_\delta(t,2\overline{\xb}-\xb),\quad \forall \xb \text{ in }\bnd^+_{2\delta},t>0$$
where $\overline{\vx}:=\underset{\yb\in\partial\omg}{\text{argmin}}\,\text{dist}(\yb,\xb)$. 
For the linear patch test case, machine precision accuracy is obtained. For the other two cases, in Figures \ref{fig:sin_linear_d} we plot the numerical results for $\vertii{\ub_\delta-\ub_0}_{L^2(\omg;\re^d)}$ and $\|\oG\ov{\vv}_\delta\|_{L^\infty(\bnd_{2\delta}^-)}$. The convergence rates are very similar to the rates observed in static cases: second order convergence rate is observed for $\ub_\delta-\ub_0$ in the $L^2$ norm and first order convergence rate for $\|\oG\ov{\vv}_\delta\|_{L^\infty(\bnd_{2\delta}^-)}$. Therefore, the numerical rate of convergence for $\vertii{\ub_\delta-\ub_0}_{L^2(\omg;\re^d)}$ is half order higher than the theoretical rate provided by Theorem \ref{thm:dynamicconv}.
\end{enumerate}

\section{Conclusions and Future Directions}\label{sec:concl}

The paper illustrates several aspects behind convergence of vector-valued solutions to  nonlocal systems of equations to classical counterparts. As a particular aspect, we investigate the relationship between the smoothness of the extension of a classical solution and the degree of convergence as the horizon of interaction shrinks to zero. The conclusion is that in order to increase the rate of convergence one must align the nonlocal boundary data with smoother extensions of the classical solution. A better agreement with the classical extension and a higher degree of regularity would guarantee a faster convergence rate.

The optimality of the theoretical results regarding the rates of convergence is an important question for future explorations, especially as the theoretical results and numerical experiments exhibit an $O(
\delta^{1/2})$ gap in the order of convergence. It remains to be shown if one can produce examples of nonlocal solutions which match the theoretical rates of convergence. For this direction, a sensible approach would be to explore more irregular domains, or problems with more irregular data. These studies can further be carried out also in terms of other types of regularity to nonlinear problems, e.g. nonlinear diffusion or wave propagation in the nonlocal setting. While it is expected that more general kernels can be considered (due to the availability of results in \cite{fosstraces,foss2019nonlocal}), a much more complex investigation is needed for dynamic fracture, where the limit for the nonlocal solutions would first need to be identified. Nonlinear problems in the vectorial setting are additionally beset with well-posedness issues, especially as nonlocal Korn-type inequalities do not have the same strength as for classical operators, so higher-order nonlinearities need to be tackled with more sophisticated analysis tools.

Finally, these results have been obtained for Dirichlet-type boundary conditions, so Neumann-type boundary problems such as \cite{you2020asymptotically,yu2021asymptotically} will be investigated elsewhere. An additional open problem would be the investigation of convergence rates for higher-order problems, such as for the nonlocal biharmonic with clamped or hinged boundary conditions, as introduced in \cite{RTY}.

\section*{Conflict of Interest Statement} 
 
On behalf of all authors, the corresponding author states that there is no conflict of interest.  
 
\bibliographystyle{spmpsci} 
\bibliography{yyu}

\end{document}

%% file: CommPreamble.tex
\usepackage{anyfontsize}
\usepackage{enumerate}
\usepackage{esint}
\usepackage{etoolbox}
\usepackage{isomath}
\usepackage{mathrsfs}
\usepackage{mathtools}
\usepackage{multicol}
\usepackage{pgfplots}
\usepackage{scalerel}
\usepackage{stackengine}
\usepackage{tabulary}
\usepackage{tikz}

\newcommand{\ds}[1]{\displaystyle{#1}}

\newcommand{\qqquad}{\qquad\qquad}
\newcommand{\qqqquad}{\qqquad\qqquad}

\newcommand{\set}[1]{\mathcal{#1}}

\newcommand{\nats}{\mathbb{N}}

\newcommand{\re}{\mathbb{R}}

\newcommand{\ren}{\re^n}

\newcommand{\bll}{\set{B}}

\newcommand{\dd}{\mathrm{d}}

\newcommand{\pder}[2]{\frac{\partial #1}{\partial #2}}

\makeatletter
\newcommand*\bdot{\mathpalette\bdot@{.65}}
\newcommand*\bdot@[2]{\mathbin{\vcenter{\hbox{\scalebox{#2}{$\m@th#1\bullet$}}}}}
\makeatother

\makeatletter
\newcommand*\bddot{\mathpalette\bddot@{.65}}
\newcommand*\bddot@[2]{\mathbin{\vcenter{\hbox{\scalebox{#2}
    {$\m@th#1\smash{{}_{\bullet}^{\bullet}}$}}}}}
\makeatother

\newcommand{\circled}[2][]{%
  \tikz[baseline=(char.base)]{%
    \node[shape = circle, draw, inner sep = .5pt]
    (char) {\phantom{\ifblank{#1}{#2}{#1}}};%
    \node at (char.center) {\makebox[0pt][c]{#2}};}}
\robustify{\circled}

\newcommand{\ve}[1]{\vectorsym{#1}}

\newcommand{\vom}{\ve{\omega}}

\newcommand{\vf}{\ve{f}}

\newcommand{\vK}{\ve{K}}

\newcommand{\vu}{\ve{u}}
\newcommand{\vv}{\ve{v}}
\newcommand{\vw}{\ve{w}}
\newcommand{\vx}{\ve{x}}
\newcommand{\vy}{\ve{y}}

\newcommand{\vz}{\ve{z}}

\newcommand{\lp}{\left(}
\newcommand{\rp}{\right)}
\newcommand{\ls}{\left[}
\newcommand{\rs}{\right]}
\newcommand{\lc}{\left\{}

\newcommand{\rper}{\right.}

\newcommand{\te}[1]{\tensorsym{#1}}
\newcommand{\tA}{\te{A}}
\newcommand{\tB}{\te{B}}
\newcommand{\tC}{\te{C}}

\newcommand{\tI}{\te{I}}

\newcommand{\tv}{\te{v}}
\newcommand{\tw}{\te{w}}

\def\dist{\operatorname{dist}}
\def\diverge{\operatorname{div}}

\def\sym{\operatorname{sym}}

\newcommand{\nd}{\quad\text{ and }\quad}
\newcommand{\frl}{\quad\text{ for all }}

\newcommand{\ov}[1]{\overline{#1}}

\newcommand{\vep}{\varepsilon}

\newcommand{\dom}{\Omega}
\newcommand{\bnd}{\Gamma}

\newcommand{\oA}{\mathcal{A}}
\newcommand{\oB}{\mathcal{B}}

\newcommand{\oG}{\mathcal{G}}
\newcommand{\oL}{\mathcal{L}}

\newcommand{\oR}{\mathcal{R}}

\stackMath
\newcommand\reallywidecheck[1]{%
\savestack{\tmpbox}{\stretchto{%
  \scaleto{%
    \scalerel*[\widthof{\ensuremath{#1}}]{\kern-.6pt\bigwedge\kern-.6pt}%
    {\rule[-\textheight/2]{1ex}{\textheight}}
  }{\textheight}%
}{0.5ex}}%
\stackon[1pt]{#1}{\scalebox{-1}{\tmpbox}}%
}